\documentclass[12pt]{article}
\usepackage{subfigure}
\usepackage{tikz}
\usepackage{pstricks, amssymb, xr, multind, multicol,textcomp,hyperref}
\def\hfl#1{\smash{\mathop{\hbox to 10mm{\rightarrowfill}}\limits^{\textstyle
#1}}}
\def\hookfl#1{\smash{\mathop{\hookrightarrow}\limits^{\textstyle
#1}}}
\def\norm#1{\parallel #1 \parallel}
\newtheorem{proposition}[equation]{Proposition} 
\newtheorem{corollary}[equation]{Corollary} 
\newtheorem{theorem}[equation]{Theorem} 
\newtheorem{exa}[equation]{Example} 
\newtheorem{ex}[equation]{Exercise} 
\newtheorem{s-ex}[equation]{Side-exercise} 
\newtheorem{exas}[equation]{Examples} 
\newtheorem{lemma}[equation]{Lemma} 
 
\newtheorem{remar}[equation]{Remark} 
\newtheorem{remars}[equation]{Remarks} 
\newtheorem{nota}[equation]{Notation} 
\newtheorem{sremar}[equation]{Side-remark} 
\newtheorem{definitio}[equation]{Definition}

\newenvironment{remark}{\begin{remar} \rm }{\end{remar}}

\newenvironment{exo}{\begin{ex} \rm }{\end{ex}} 
 
\newenvironment{examples}{\begin{exas} \rm }{\end{exas}} 
\newenvironment{example}{\begin{exa} \rm }{\end{exa}} 
\newenvironment{definition}{\begin{definitio} \rm }{\end{definitio}} 
\newcommand{\blowup}[2]{B\!\ell(#1,#2)}

\newcommand{\HH}{{\mathbb{H}}} 
\newcommand{\KK}{\mathbb{K}}

\newcommand{\CA}{{\cal A}}

\newcommand{\ansothree}{\beta}

\newcommand{\CD}{{\cal D}} 
 
\newcommand{\CF}{{\cal F}} 
 
\newcommand{\handlebody}{{\cal H}} 
\newcommand{\CI}{{\cal I}} 
 
\newcommand{\CK}{{\cal K}}
\newcommand{\CL}{{\cal L}} 
\newcommand{\linearmapanv}{{\cal L}_v} 
\newcommand{\CM}{{\cal M}} 

\newcommand{\normbun}{\mathfrak{V}} 
\newcommand{\normnu}{\nu} 
\newcommand{\neightub}{N}
\newcommand{\neigh}{N}
\newcommand{\nfinset}{N}
\newcommand{\CO}{{\cal O}}
\newcommand{\CP}{{\cal P}} 
\newcommand{\prop}{{\cal P}}
 
\newcommand{\tang}{{T}} 
\newcommand{\taust}{{\tau_s}}

\newcommand{\ZZ}{\mathbb{Z}}

\newcommand{\RR}{\mathbb{R}} 
\newcommand{\QQ}{\mathbb{Q}} 
\newcommand{\CC}{\mathbb{C}} 
\newcommand{\NN}{\mathbb{N}}

\newcommand{\bp}{\noindent {\sc Proof: }} 
\newcommand{\bsp}{\noindent {\sc Sketch of proof: }} 
\newcommand{\eop}{\nopagebreak \hspace*{\fill}{$\Box$} \medskip} 
\newcommand{\tata}{ \begin{tikzpicture} \useasboundingbox (-.3,-.1) rectangle (.3,.2);
\draw (-.2,0) -- (.2,0) (0,0) circle (.2);
\fill (-.2,0) circle (1.5pt) (.2,0) circle (1.5pt);
\end{tikzpicture}} 
\newcommand{\smally}{\begin{tikzpicture} \useasboundingbox (-.4,-.2) rectangle (.4,.2);
\draw (150:.15) -- (0,0) -- (-90:.15) (0,0) -- (30:.15) (30:.25) circle (.1) (150:.25) circle (.1) (-90:.25) circle (.1);
\fill (0,0) circle (1pt);
\end{tikzpicture}}
\newcommand{\pcbg}{\begin{tikzpicture}
\useasboundingbox (-.6,.1) rectangle (.6,.55);
\draw [->] (.2,0) -- (-.2,.4) node[left]{\tiny $J$};
\draw [->,draw=white,double=gray,very thick]  (-.2,0) -- (.2,.4) node[right]{\tiny $K$};
\draw [->,draw=gray]  (-.2,0) -- (.2,.4);
\end{tikzpicture}} 
\newcommand{\pcbgkj}{\begin{tikzpicture}
\useasboundingbox (-.6,.1) rectangle (.6,.55);
\draw [->,draw=gray] (.2,0) -- (-.2,.4) node[left]{\tiny $K$};
\draw [->,draw=white,double=black,very thick]  (-.2,0) -- (.2,.4) node[right]{\tiny $J$};
\draw [->]  (-.2,0) -- (.2,.4);
\end{tikzpicture}} 
\newcommand{\ncbg}{\begin{tikzpicture}
\useasboundingbox (-.6,.1) rectangle (.6,.55);
\draw [->] (-.2,0) -- (.2,.4) node[right]{\tiny $J$};
\draw [->,draw=white,double=gray,very thick]  (.2,0) -- (-.2,.4) node[left]{\tiny $K$};
\draw [->,draw=gray]  (.2,0) -- (-.2,.4);
\end{tikzpicture}} 
\newcommand{\ncbgkj}{\begin{tikzpicture}
\useasboundingbox (-.6,.1) rectangle (.6,.55);
\draw [->,draw=gray] (-.2,0) -- (.2,.4) node[right]{\tiny $K$};
\draw [->,draw=white,double=black,very thick]  (.2,0) -- (-.2,.4) node[left]{\tiny $J$};
\draw [->]  (.2,0) -- (-.2,.4);
\end{tikzpicture}} 
\newcommand{\pc}{\begin{tikzpicture} \useasboundingbox (-.5,.1) rectangle (.5,.55);
\draw [->] (.2,0) -- (-.2,.4);
\draw [->,draw=white,double=black,very thick]  (-.2,0) -- (.2,.4);
\draw [->]  (-.2,0) -- (.2,.4);
\end{tikzpicture}}
\newcommand{\pccirc}{\begin{tikzpicture} \useasboundingbox (-.5,.1) rectangle (.5,.55);
\draw [->] (.2,0) -- (-.2,.4);
\draw [->,draw=white,double=black,very thick]  (-.2,0) -- (.2,.4);
\draw [->]  (-.2,0) -- (.2,.4);
\draw [dash pattern=on 2pt off 2pt] (0,.2) circle (.283);
\end{tikzpicture}}
\newcommand{\pcortrig}{\begin{tikzpicture} \useasboundingbox (-.5,.1) rectangle (.5,.55);
\draw [->] (.2,0) -- (-.2,.4);
\draw [->,draw=white,double=black,very thick]  (-.2,0) -- (.2,.4);
\draw [->]  (-.2,0) -- (.2,.4);
\draw [draw=gray,->] (.1,.45) .. controls (0,.5) .. (-.1,.45);
\end{tikzpicture}}
\newcommand{\nc}{\begin{tikzpicture}
\useasboundingbox (-.5,.1) rectangle (.5,.55);
\draw [->] (-.2,0) -- (.2,.4);
\draw [->,draw=white,double=black,very thick]  (.2,0) -- (-.2,.4);
\draw [->]  (.2,0) -- (-.2,.4);
\end{tikzpicture}}
\newcommand{\nccirc}{\begin{tikzpicture}
\useasboundingbox (-.5,.1) rectangle (.5,.55);
\draw [->] (-.2,0) -- (.2,.4);
\draw [->,draw=white,double=black,very thick]  (.2,0) -- (-.2,.4);
\draw [->]  (.2,0) -- (-.2,.4);
\draw [dash pattern=on 2pt off 2pt] (0,.2) circle (.283);
\end{tikzpicture}}
\newcommand{\doublep}{\begin{tikzpicture} \useasboundingbox (-.5,.1) rectangle (.5,.55);
\draw [->] (.2,0) -- (-.2,.4);
\draw [->]  (-.2,0) -- (.2,.4);
\fill  (0,.2) circle (1.5pt);
\end{tikzpicture}}
\newcommand{\zerochord}{\begin{tikzpicture} \useasboundingbox (-.4,-.2) rectangle (.4,.3);
\draw [->,dash pattern=on 2pt off 2pt] (.3,0) arc (0:360:.3);
\end{tikzpicture}}
\newcommand{\onechord}{\begin{tikzpicture} \useasboundingbox (-.4,-.2) rectangle (.4,.3);
\draw [->,dash pattern=on 2pt off 2pt] (.3,0) arc (0:360:.3);
\draw (0,-.3) -- (0,.3);
\fill (0,-.3) circle (1.5pt) (0,.3) circle (1.5pt);
\end{tikzpicture}}
\newcommand{\onechordR}{\begin{tikzpicture} \useasboundingbox (-.25,-.2) rectangle (.2,.3);
\draw [->,dash pattern=on 2pt off 2pt] (0,-.3) -- (0,.4);
\draw (0,.15) arc (90:270:.15);
\fill (0,-.15) circle (1pt) (0,.15) circle (1pt);
\end{tikzpicture}}
\newcommand{\twoischord}{\begin{tikzpicture} \useasboundingbox (-.4,-.2) rectangle (.4,.3);
\draw [->,dash pattern=on 2pt off 2pt] (.3,0) arc (0:360:.3);
\draw (-45:.3) .. controls (-90:.1) ..  (-135:.3);
\draw (45:.3) .. controls (90:.1) ..  (135:.3);
\fill (-45:.3) circle (1.5pt) (45:.3) circle (1.5pt) (-135:.3) circle (1.5pt) (135:.3) circle (1.5pt);
\end{tikzpicture}}
\newcommand{\twoxchord}{\begin{tikzpicture} \useasboundingbox (-.4,-.2) rectangle (.4,.3);
\draw [->,dash pattern=on 2pt off 2pt] (.3,0) arc (0:360:.3);

\draw (-45:.3) -- (135:.3);
\draw (-135:.3) -- (45:.3);
\fill (-45:.3) circle (1.5pt) (45:.3) circle (1.5pt) (-135:.3) circle (1.5pt) (135:.3) circle (1.5pt);
\end{tikzpicture}}
\newcommand{\threeischord}{\begin{tikzpicture} \useasboundingbox (-.4,-.2) rectangle (.4,.5);
\draw [->,dash pattern=on 2pt off 2pt] (.3,0) arc (0:360:.3);
\draw (-30:.3) .. controls (-60:.1) ..  (-90:.3);
\draw (30:.3) .. controls (60:.1) ..  (90:.3);
\draw  (-150:.3) .. controls (180:.1) ..  (150:.3);
\fill (-30:.3) circle (1.5pt) (30:.3) circle (1.5pt) (-90:.3) circle (1.5pt) (90:.3) circle (1.5pt) (-150:.3) circle (1.5pt) (150:.3) circle (1.5pt);
\end{tikzpicture}}
\newcommand{\threeparchord}{\begin{tikzpicture} \useasboundingbox (-.4,-.2) rectangle (.4,.5);
\draw [->,dash pattern=on 2pt off 2pt] (.3,0) arc (0:360:.3);
\draw (-30:.3) -- (150:.3);
\draw (30:.3) .. controls (60:.1) ..  (90:.3);
\draw  (-150:.3) .. controls (-120:.1) ..  (-90:.3);
\fill (-30:.3) circle (1.5pt) (30:.3) circle (1.5pt) (-90:.3) circle (1.5pt) (90:.3) circle (1.5pt) (-150:.3) circle (1.5pt) (150:.3) circle (1.5pt);
\end{tikzpicture}}
\newcommand{\threexchord}{\begin{tikzpicture} \useasboundingbox (-.4,-.2) rectangle (.4,.5);
\draw [->,dash pattern=on 2pt off 2pt] (.3,0) arc (0:360:.3);
\draw (-150:.3) .. controls (180:.1) ..  (150:.3);
\draw (30:.3) -- (-90:.3);
\draw  (90:.3) -- (-30:.3);
\fill (-30:.3) circle (1.5pt) (30:.3) circle (1.5pt) (-90:.3) circle (1.5pt) (90:.3) circle (1.5pt) (-150:.3) circle (1.5pt) (150:.3) circle (1.5pt);
\end{tikzpicture}}
\newcommand{\threexparchord}{\begin{tikzpicture} \useasboundingbox (-.4,-.2) rectangle (.4,.5);
\draw [->,dash pattern=on 2pt off 2pt] (.3,0) arc (0:360:.3);
\draw (-150:.3) -- (-30:.3);
\draw (90:.3) -- (-90:.3);
\draw  (150:.3) -- (30:.3);
\fill (-30:.3) circle (1.5pt) (30:.3) circle (1.5pt) (-90:.3) circle (1.5pt) (90:.3) circle (1.5pt) (-150:.3) circle (1.5pt) (150:.3) circle (1.5pt);
\end{tikzpicture}}
\newcommand{\threexxchord}{\begin{tikzpicture} \useasboundingbox (-.4,-.2) rectangle (.4,.5);
\draw [->,dash pattern=on 2pt off 2pt] (.3,0) arc (0:360:.3);
\draw (-150:.3) -- (30:.3);
\draw (-30:.3) -- (150:.3);
\draw  (90:.3) -- (-90:.3);
\fill (-30:.3) circle (1.5pt) (30:.3) circle (1.5pt) (-90:.3) circle (1.5pt) (90:.3) circle (1.5pt) (-150:.3) circle (1.5pt) (150:.3) circle (1.5pt);
\end{tikzpicture}}

\newcommand{\exjactwosonebis}{\begin{tikzpicture} \useasboundingbox (-1.4,-.2) rectangle (1.4,.5);
\draw [->,dash pattern=on 2pt off 2pt] (0,-.3) arc (-90:270:.3);
\draw [->,dash pattern=on 2pt off 2pt] (-1,-.3) arc (-90:270:.3);
\draw (-1.3,0) -- (30:.3);
\draw (0,0) -- (-140:.3) (0,0) -- (-40:.3) (0,0) -- (0,.3);
\fill (-40:.3) circle (1.5pt) (-140:.3) circle (1.5pt) (90:.3) circle (1.5pt) (30:.3) circle (1.5pt) (-1.3,0) circle (1.5pt) (0,0) circle (1.5pt) (1.2,0) circle (1.5pt) (.8,0) circle (1.5pt);
\draw (1,0) circle (.2) (.8,0) -- (1.2,0);
\end{tikzpicture}}
\newcommand{\lolli}{\begin{tikzpicture} \useasboundingbox (-.4,-.2) rectangle (.4,.3);
\draw [->,dash pattern=on 2pt off 2pt] (.3,0) arc (0:360:.3);
\draw (0,-.3) -- (0,-.1) (0,.05) circle (.15);
\fill (0,-.3) circle (1.5pt) (0,-.1) circle (1.5pt);
\end{tikzpicture}}
\newcommand{\tatasone}{\begin{tikzpicture} \useasboundingbox (-.4,-.2) rectangle (.4,.3);
\draw [->,dash pattern=on 2pt off 2pt] (.3,0) arc (0:360:.3);
\draw (-.2,0) -- (.2,0) (0,0) circle (.2);
\fill (-.2,0) circle (1.5pt) (.2,0) circle (1.5pt);
\end{tikzpicture}}
\newcommand{\tatasonetw}{\begin{tikzpicture} \useasboundingbox (-.4,-.2) rectangle (.4,.3);
\draw [->,dash pattern=on 2pt off 2pt] (.3,0) arc (0:360:.3);
\draw (.2,0) arc (0:180:.2) (-.2,0) -- (-.05,0) .. controls (.05,0) and (.2,-.3) .. (.2,0) 
(-.2,0) arc (-180:-90:.2) (0,-.2) .. controls (.15,-.2) and (0,.1) .. (.2,0);
\fill (-.2,0) circle (1.5pt) (.2,0) circle (1.5pt);
\end{tikzpicture}}
\newcommand{\haltsone}{\begin{tikzpicture} \useasboundingbox (-.4,-.2) rectangle (.4,.3);
\draw [->,dash pattern=on 2pt off 2pt] (.35,0) arc (0:360:.35);
\draw (0,-.07) -- (0,.07) (0,.19) circle (.12) (0,-.19) circle (.12);
\fill (0,-.07) circle (1.5pt) (0,.07) circle (1.5pt);
\end{tikzpicture}}
\newcommand{\trivAS}{\begin{tikzpicture} \useasboundingbox (-.4,-.3) rectangle (.4,.4);
\draw (0,0) -- (0,-.4) (0,0) -- (50:.4) (0,0) -- (130:.4);
\fill (0,0) circle (1.5pt);
\end{tikzpicture}}
\newcommand{\trivASop}{\begin{tikzpicture} \useasboundingbox (-.4,-.2) rectangle (.4,.3);
\draw (0,0) -- (0,-.2);
\draw (0,0) .. controls (50:.3) ..  (130:.5);
\draw [draw=white,double=black,very thick] (0,0) .. controls (130:.3) ..  (50:.5);
\fill (0,0) circle (1.5pt);
\end{tikzpicture}}
\newcommand{\stuy}{\begin{tikzpicture} \useasboundingbox (-.1,-.25) rectangle (.7,.4);
\draw [->,dash pattern=on 2pt off 2pt] (0,-.2) -- (.6,-.2);
\draw (.3,-.2) -- (.3,0) -- (.15,.2) (.3,0) -- (.45,.2);
\fill (.3,-.2) circle (1.5pt) (.3,0) circle (1.5pt);
\end{tikzpicture}}
\newcommand{\stui}{\begin{tikzpicture} \useasboundingbox (-.1,-.25) rectangle (.7,.4);
\draw [->,dash pattern=on 2pt off 2pt] (0,-.2) -- (.6,-.2);
\draw (.15,-.2) -- (.15,.2)  (.45,.2) -- (.45,-.2);
\fill (.15,-.2) circle (1.5pt) (.45,-.2) circle (1.5pt);
\end{tikzpicture}}
\newcommand{\stux}{\begin{tikzpicture} \useasboundingbox (-.1,-.25) rectangle (.7,.4);
\draw [->,dash pattern=on 2pt off 2pt] (0,-.2) -- (.6,-.2);
\draw (.15,.2) -- (.45,-.2);
\draw [draw=white,double=black,very thick] (.15,-.2) -- (.45,.2);
\fill (.15,-.2) circle (1.5pt) (.45,-.2) circle (1.5pt);
\end{tikzpicture}}
\newcommand{\asunor}{\begin{tikzpicture} \useasboundingbox (-.1,-.25) rectangle (.7,.4);
\draw [->,dash pattern=on 2pt off 2pt] (0,-.2) -- (.6,-.2);
\draw (.3,.2) -- (.3,-.2);
\fill (.3,-.2) circle (1.5pt);
\end{tikzpicture}}
\newcommand{\asunorc}{\begin{tikzpicture} \useasboundingbox (-.1,-.25) rectangle (.7,.4);
\draw [->,dash pattern=on 2pt off 2pt] (0,-.2) -- (.6,-.2);
\draw (.3,.2) .. controls (.1,0) and  (.05,-.35)  .. (.15,-.35)  .. controls (.25,-.35) .. (.3,-.2);
\fill (.3,-.2) circle (1.5pt);
\end{tikzpicture}}
\newcommand{\stubone}{\begin{tikzpicture} \useasboundingbox (-.1,-.25) rectangle (1.1,.4);
\draw [->,dash pattern=on 2pt off 2pt] (0,-.2) -- (1,-.2);
\draw (.3,.3) -- (.3,-.2) (.7,.3) -- (.7,-.2);
\fill (.3,-.2) circle (1.5pt) (.7,-.2) circle (1.5pt);
\end{tikzpicture}}
\newcommand{\stubtwo}{\begin{tikzpicture} \useasboundingbox (-.1,-.25) rectangle (1.1,.4);
\draw [->,dash pattern=on 2pt off 2pt] (0,-.2) -- (1,-.2);
\draw (.4,.3)  .. controls (.1,-.1) and  (.3,-.35)  .. (.4,-.35)  .. controls (.5,-.35) .. (.7,-.2) (.7,.3) -- (.4,-.2);
\fill (.4,-.2) circle (1.5pt) (.7,-.2) circle (1.5pt);
\end{tikzpicture}}
\newcommand{\stubthree}{\begin{tikzpicture} \useasboundingbox (-.1,-.25) rectangle (1.1,.4);
\draw [->,dash pattern=on 2pt off 2pt] (0,-.2) -- (1,-.2);
\draw (.2,.3)  .. controls (.2,-.2) and  (.7,-.1)  .. (.7,0)  .. controls (.7,.1) .. (.5,.15) (.5,.3) -- (.5,-.2);
\fill (.5,-.2) circle (1.5pt) (.5,.15) circle (1.5pt);
\end{tikzpicture}}
\newcommand{\ihxone}{\begin{tikzpicture} \useasboundingbox (-.4,-.3) rectangle (.4,.4);
\draw (-90:.4) -- (0,0) -- (30:.4) (0,0) -- (150:.4);
\draw (-.3,-.25) -- (0,-.25);
\fill (0,-.25) circle (1.5pt) (0,0) circle (1.5pt);
\end{tikzpicture}}
\newcommand{\ihxtwo}{\begin{tikzpicture} \useasboundingbox (-.4,-.3) rectangle (.4,.4);
\draw (-90:.4) -- (0,0) -- (30:.4) (0,0) -- (150:.4);
\draw [out=0,in=-60] (-.3,-.25) to (30:.25);
\fill (30:.25) circle (1.5pt) (0,0) circle (1.5pt);
\end{tikzpicture}}
\newcommand{\ihxthree}{\begin{tikzpicture} \useasboundingbox (-.4,-.3) rectangle (.4,.4);
\draw (-90:.4) -- (0,0) -- (30:.4) (0,0) -- (150:.4);
\draw (-.3,-.25) .. controls (-.2,-.25) and  (.25,-.2)  .. (.25,0)  .. controls (.25,.2) and (120:.35) .. (150:.25);
\fill (150:.25) circle (1.5pt) (0,0) circle (1.5pt);
\end{tikzpicture}}
\newcommand{\ihxi}{\begin{tikzpicture}[scale=.7] \useasboundingbox (-.4,-.3) rectangle (.4,.4);
\draw (-.3,-.3) -- (.3,-.3) (-.3,.3) -- (.3,.3) (0,-.3) -- (0,.3);
\fill (0,-.3) circle (2pt) (0,.3) circle (2pt);
\end{tikzpicture}}
\newcommand{\ihxh}{\begin{tikzpicture}[scale=.7] \useasboundingbox (-.4,-.3) rectangle (.4,.4);
\draw (-.3,-.3) -- (-.3,.3) (.3,-.3) -- (.3,.3) (-.3,0) -- (.3,0);
\fill (-.3,0) circle (2pt) (.3,0) circle (2pt);
\end{tikzpicture}}
\newcommand{\ihxx}{\begin{tikzpicture}[scale=.7] \useasboundingbox (-.4,-.3) rectangle (.4,.4);
\draw (-.3,-.3) -- (.3,.3) (.3,-.3) -- (-.3,.3) (-.15,-.15) -- (.15,-.15);
\fill (-.15,-.15) circle (2pt) (.15,-.15) circle (2pt);
\end{tikzpicture}}
\newcommand{\singtrefoil}{\begin{tikzpicture}[scale=.5] \useasboundingbox (-1,-.5) rectangle (1,.7);
\draw [out=120,in=60] (30:.7) to (150:.3);
\draw [out=240,in=180] (150:.3) to (-90:.7);
\draw [out=0,in=-60] (-90:.7) to (30:.3);
\draw [out=120,in=60] (30:.3) to (150:.7);
\draw [out=240,in=180] (150:.7) to (-90:.3);
\draw [out=0,in=-60,->] (-90:.3) to (30:.7);
\fill (-30:.39) circle (2.5pt) (90:.39) circle (2.5pt) (-150:.39) circle (2.5pt);
\end{tikzpicture}}
\newcommand{\righttrefoil}{\begin{tikzpicture} \useasboundingbox (-1,-.5) rectangle (1,.7);
\draw [out=120,in=60] (30:.7) to (150:.3);
\draw [out=0,in=-60] (-90:.7) to (30:.3);
\draw [out=240,in=180] (150:.7) to (-90:.3);
\draw [out=240,in=180,draw=white,double=black,very thick] (150:.3) to (-90:.7);
\draw [out=120,in=60,draw=white,double=black,very thick] (30:.3) to (150:.7);
\draw [out=0,in=-60,draw=white,double=black,very thick] (-90:.3) to (30:.7);
\draw [out=0,in=-60,->] (-90:.3) to (30:.7);
\end{tikzpicture}}

\newcommand{\dquatSTU}{\begin{tikzpicture} \useasboundingbox (-.8,-.4) rectangle (.8,.6);
\draw [->,dash pattern=on 2pt off 2pt] (70:.5) arc (70:110:.5);
\draw [->,dash pattern=on 2pt off 2pt] (170:.5) arc (170:240:.5);
\draw [->,dash pattern=on 2pt off 2pt] (-60:.5) arc (-60:10:.5);
\draw (0,0) -- (-25:.5) (0,0) -- (90:.5) (0,0) -- (205:.5);
\fill (90:.5) circle (1.5pt) (-25:.5) circle (1.5pt) (205:.5) circle (1.5pt) (0,0) circle (1.5pt);
\end{tikzpicture}}
\newcommand{\dquatTun}{\begin{tikzpicture} \useasboundingbox (-.8,-.4) rectangle (.8,.6);
\draw [->,dash pattern=on 2pt off 2pt] (70:.5) arc (70:110:.5);
\draw [->,dash pattern=on 2pt off 2pt] (170:.5) arc (170:240:.5);
\draw [->,dash pattern=on 2pt off 2pt] (-60:.5) arc (-60:10:.5);
\draw (90:.5) --(-40:.5) (-10:.5) -- (205:.5);
\fill (90:.5) circle (1.5pt) (205:.5) circle (1.5pt) (-40:.5) circle (1.5pt) (-10:.5) circle (1.5pt);
\end{tikzpicture}}
\newcommand{\dquatTdeux}{\begin{tikzpicture} \useasboundingbox (-.8,-.4) rectangle (.8,.6);
\draw [->,dash pattern=on 2pt off 2pt] (70:.5) arc (70:110:.5);
\draw [->,dash pattern=on 2pt off 2pt] (170:.5) arc (170:240:.5);
\draw [->,dash pattern=on 2pt off 2pt] (-60:.5) arc (-60:10:.5);
\draw (90:.5) -- (220:.5) (190:.5) -- (-25:.5);
\fill (90:.5) circle (1.5pt) (220:.5) circle (1.5pt) (-25:.5) circle (1.5pt) (190:.5) circle (1.5pt);
\end{tikzpicture}}
\newcommand{\dquatTtrois}{\begin{tikzpicture} \useasboundingbox (-.8,-.4) rectangle (.8,.6);
\draw [->,dash pattern=on 2pt off 2pt] (70:.5) arc (70:110:.5);
\draw [->,dash pattern=on 2pt off 2pt] (170:.5) arc (170:240:.5);
\draw [->,dash pattern=on 2pt off 2pt] (-60:.5) arc (-60:10:.5);
\draw (90:.5) -- (-10:.5) (-40:.5) -- (205:.5);
\fill (90:.5) circle (1.5pt) (205:.5) circle (1.5pt) (-40:.5) circle (1.5pt) (-10:.5) circle (1.5pt);
\end{tikzpicture}}
\newcommand{\dquatTqua}{\begin{tikzpicture} \useasboundingbox (-.8,-.4) rectangle (.8,.6);
\draw [->,dash pattern=on 2pt off 2pt] (70:.5) arc (70:110:.5);
\draw [->,dash pattern=on 2pt off 2pt] (170:.5) arc (170:240:.5);
\draw [->,dash pattern=on 2pt off 2pt] (-60:.5) arc (-60:10:.5);
\draw (90:.5) --(190:.5) (220:.5) -- (-25:.5);
\fill (90:.5) circle (1.5pt) (220:.5) circle (1.5pt) (-25:.5) circle (1.5pt) (190:.5) circle (1.5pt);
\end{tikzpicture}}
\newcommand{\ischord}{\begin{tikzpicture} \useasboundingbox (-.2,-.2) rectangle (.4,.3);
\draw [->,dash pattern=on 2pt off 2pt] (-90:.3) arc (-90:90:.3);
\draw (-45:.3) -- (45:.3);
\fill (-45:.3) circle (1.5pt) (45:.3) circle (1.5pt);
\end{tikzpicture}}
\newcommand{\pkink}{\begin{tikzpicture} \useasboundingbox (-.1,-.1) rectangle (.5,.2);
\draw [out=-90,in=-90,->] (.4,0) to (0,.25);
\draw [out=90,in=90,draw=white,double=black,very thick] (0,-.25) to (.4,0);
\end{tikzpicture}}
\newcommand{\nkink}{\begin{tikzpicture} \useasboundingbox (-.1,-.1) rectangle (.5,.2);
\draw [out=90,in=90] (0,-.2) to (.4,0);
\draw [out=-90,in=-90,draw=white,double=black,very thick] (.4,0) to (0,.2);
\draw [out=-90,in=-90,->] (.4,0) to (0,.2);
\end{tikzpicture}}
\newcommand{\mkink}{\begin{tikzpicture} \useasboundingbox (-.1,-.1) rectangle (.5,.2);
\draw [out=90,in=90,<-] (0,-.2) to (.4,0);
\draw [out=-90,in=-90,draw=white,double=black,very thick] (.4,0) to (0,.2);
\draw [out=-90,in=-90] (.4,0) to (0,.2);
\end{tikzpicture}}

\newcommand{\noedgeuu}{\begin{tikzpicture} \useasboundingbox (-.6,-.2) rectangle (.2,.4);
\draw [->,dash pattern=on 2pt off 2pt] (0,-.3) -- (0,.3);
\draw [->,dash pattern=on 2pt off 2pt] (-.4,.3) -- (-.4,-.3);
\end{tikzpicture}}
\newcommand{\edgeuu}{\begin{tikzpicture} \useasboundingbox (-.7,-.2) rectangle (.2,.4);
\draw [->,dash pattern=on 2pt off 2pt] (0,-.3) -- (0,.3);
\draw [->,dash pattern=on 2pt off 2pt] (-.5,.3) -- (-.5,-.3);
\draw (-.5,0) -- (0,0);
\fill (-.5,0) circle (1.5pt) (0,0) circle (1.5pt);
\end{tikzpicture}}
\newcommand{\xchangeuu}{\begin{tikzpicture} \useasboundingbox (-.7,-.2) rectangle (.4,.4);
\draw [->,dash pattern=on 2pt off 2pt] (0,0) -- (0,.3);
\draw (0,-.3) -- (0,-.2);
\draw [->,dash pattern=on 2pt off 2pt] (.1,.1) .. controls (.15,.1) and (.25,.05) .. (.25,0) ..
controls (.25,-.05) and (.1,-.1) .. (0,-.1) .. controls (-.1,-.1) and (-.5,-.2) ..  (-.5,-.3);
\draw (0,-.3) -- (0,-.2);
\draw [dash pattern=on 2pt off 2pt] (-.5,.3) .. controls (-.5,.2) and (-.2,.1) .. (-.1,.1);
\end{tikzpicture}}

\newcommand{\edget}{\begin{tikzpicture} \useasboundingbox (-1.2,-.2) rectangle (.3,.4);
\draw (-1,.1) -- (-.9,0) -- (-1,-.1) (-.9,0) -- (.1,0) -- (.2,.1) (.1,0) -- (.2,-.1);
\fill (-.9,0) circle (1.5pt) (.1,0) circle (1.5pt);
\end{tikzpicture}}
\newcommand{\edgehopf}{\begin{tikzpicture} \useasboundingbox (-1.2,-.2) rectangle (.3,.4);
\draw [->] (-.3,0) arc (0:90:.2);
\draw [->] (-.1,0) arc (0:90:.2);
\draw (-.1,0) arc (0:-200:.2) (-.5,.2) arc (90:270:.2) (-.3,0) arc (0:-20:.2);
\draw (-1,.1) -- (-.9,0) -- (-1,-.1) (-.9,0) -- (-.7,0) (-.1,0) -- (.1,0) -- (.2,.1) (.1,0) -- (.2,-.1);
\fill (-.9,0) circle (1.5pt) (.1,0) circle (1.5pt);
\end{tikzpicture}}
\newcommand{\annulone}{\begin{tikzpicture} \useasboundingbox (-1.3,-.8) rectangle (1.3,1.2);
\draw [fill=gray!20,draw=white] (0,0) circle (.9);
\draw [fill=white,draw=white] (0,0) circle (.5) (.6,-.2) rectangle (.8,.2);
\draw [->,dash pattern=on 2pt off 2pt] (0,-1) -- (0,1);
\draw (.7,0) node{\tiny $A$} (0,-.2) -- (-.3,0) -- (0,0) (-.3,0) -- (0,.2) (-.5,-1) -- (0,-.7) (.1,-.6) node{\tiny $v$} (.2,.75) node{\tiny $\alpha_2$} (.2,-.8) node{\tiny $\alpha_1$} (.25,0) node{\tiny $D$} (-.3,-.7) node{\tiny $\beta$} (.8,-.8) node{$\Gamma_1$};
\fill (0,-.7) circle (1.5pt) (0,-.2) circle (1.5pt) (0,0) circle (1.5pt) (0,.2) circle (1.5pt) (-.3,0) circle (1.5pt);
\end{tikzpicture}}
\newcommand{\annultwo}{\begin{tikzpicture} \useasboundingbox (-1.3,-.8) rectangle (1.3,1.2);
\draw [fill=gray!20,draw=white] (0,0) circle (.9);
\draw [fill=white,draw=white] (0,0) circle (.5) (.6,-.2) rectangle (.8,.2);
\draw [->,dash pattern=on 2pt off 2pt] (0,-1) -- (0,1);
\draw (.7,0) node{\tiny $A$} (0,-.2) -- (-.3,0) -- (0,0) (-.3,0) -- (0,.2) (-.5,-1) .. controls (.1,-.7) and (.55,-.4) .. (.55,0) .. controls  (.55,.4) and  (.2,.7) ..  (0,.7) (-.1,.8) node{\tiny $v$} (.2,.75) node{\tiny $\alpha_2$} (.25,0) node{\tiny $D$} (.2,-.8) node{\tiny $\alpha_1$} (-.3,-.7) node{\tiny $\beta$} (.8,-.8) node{$\Gamma_2$};
\fill (0,.7) circle (1.5pt) (0,-.2) circle (1.5pt) (0,0) circle (1.5pt) (0,.2) circle (1.5pt) (-.3,0) circle (1.5pt);
\end{tikzpicture}}
\newcommand{\annulthree}{\begin{tikzpicture} \useasboundingbox (-1.3,-.8) rectangle (1.3,1.2);
\draw [fill=gray!20,draw=white] (0,0) circle (.9);
\draw [fill=white,draw=white] (0,0) circle (.5) (.6,-.2) rectangle (.8,.2);
\draw [->,dash pattern=on 2pt off 2pt] (0,-1) -- (0,1);
\draw (.7,0) node{\tiny $A$} (-1,0) -- (0,0) (-.5,-1) -- (0,-.7) (.1,-.6) node{\tiny $v$} (.25,0) node{\tiny $D$} (.2,.75) node{\tiny $\alpha_2$} (.2,-.8) node{\tiny $\alpha_1$} (-.7,.15) node{\tiny $\alpha_3$} (-.3,-.7) node{\tiny $\beta$} (.8,-.8) node{$\Gamma_1$};
\fill (0,-.7) circle (1.5pt) (0,0) circle (1.5pt);
\end{tikzpicture}}
\newcommand{\annulfour}{\begin{tikzpicture} \useasboundingbox (-1.3,-.8) rectangle (1.3,1.2);
\draw [fill=gray!20,draw=white] (0,0) circle (.9);
\draw [fill=white,draw=white] (0,0) circle (.5) (.6,-.2) rectangle (.8,.2);
\draw [->,dash pattern=on 2pt off 2pt] (0,-1) -- (0,1);
\draw  (-.5,-1) .. controls (.1,-.7) and (.55,-.4) .. (.55,0) .. controls  (.55,.4) and  (.2,.7) ..  (0,.7) (.7,0) node{\tiny $A$} (-1,0) -- (0,0) (-.1,.8) node{\tiny $v$} (.25,0) node{\tiny $D$} (.2,.75) node{\tiny $\alpha_2$} (.2,-.8) node{\tiny $\alpha_1$} (-.7,.15) node{\tiny $\alpha_3$} (-.3,-.7) node{\tiny $\beta$} (.8,-.8) node{$\Gamma_2$};
\fill (0,.7) circle (1.5pt) (0,0) circle (1.5pt);
\end{tikzpicture}}
\newcommand{\annulfive}{\begin{tikzpicture} \useasboundingbox (-1.3,-.8) rectangle (1.3,1.2);
\draw [fill=gray!20,draw=white] (0,0) circle (.9);
\draw [fill=white,draw=white] (0,0) circle (.5) (.6,-.2) rectangle (.8,.2);
\draw [->,dash pattern=on 2pt off 2pt] (0,-1) -- (0,1);
\draw (-.5,-1) .. controls (.1,-.7) and (.55,-.4) .. (.55,0) .. controls  (.55,.4) and  (.2,.7) ..  (0,.7) .. controls  (-.25,.7) and  (-.7,.35) ..  (-.7,0)  (.7,0) node{\tiny $A$} (-1,0) -- (0,0) (-.8,.1) node{\tiny $v$} (.25,0) node{\tiny $D$} (.2,.75) node{\tiny $\alpha_2$} (.2,-.8) node{\tiny $\alpha_1$} (-.7,-.15) node{\tiny $\alpha_3$} (-.3,-.7) node{\tiny $\beta$} (.8,-.8) node{$\Gamma_3$};
\fill (-.7,0) circle (1.5pt) (0,0) circle (1.5pt);
\end{tikzpicture}}
\newcommand{\dmer}{\begin{tikzpicture} \useasboundingbox (-.5,-.8) rectangle (.2,1);
\draw [->,dash pattern=on 2pt off 2pt] (0,-.8) -- (0,.8);
\draw (0,-.2) -- (-.3,0) -- (0,0) (-.3,0) -- (0,.2);
\fill (0,-.2) circle (1.5pt) (0,0) circle (1.5pt) (0,.2) circle (1.5pt) (-.3,0) circle (1.5pt);
\end{tikzpicture}}
\newcommand{\cirt}{\begin{tikzpicture} \useasboundingbox (-1,-.8) rectangle (1,1);
\draw [->,dash pattern=on 2pt off 2pt]  (-.8,0) arc (-180:180:.8);
\draw (-135:.8) -- (0,-.15) -- (0,.15) -- (135:.8)  (0,.15) -- (45:.8) (0,-.15) --  (-45:.8);
\fill (-135:.8) circle (1.5pt) (135:.8) circle (1.5pt) (45:.8) circle (1.5pt) (-45:.8) circle (1.5pt) (0,-.15) circle (1.5pt) (0,.15) circle (1.5pt);
\end{tikzpicture}}
\newcommand{\cirtd}{\begin{tikzpicture} \useasboundingbox (-1,-.8) rectangle (1,1);
\draw [->,dash pattern=on 2pt off 2pt]  (-.8,0) arc (-180:180:.8);
\draw (-135:.8) -- (0,-.15) -- (0,.15) -- (135:.8)  (0,.15) -- (60:.8) (0,-.15) --  (-60:.8);
\draw (-30:.8) -- (.4,0) -- (30:.8) (.4,0) -- (.8,0);
\fill (-135:.8) circle (1.5pt) (135:.8) circle (1.5pt) (60:.8) circle (1.5pt) (-60:.8) circle (1.5pt) (30:.8) circle (1.5pt) (0:.8) circle (1.5pt) (-30:.8) circle (1.5pt) (.4,0) circle (1.5pt) (0,-.15) circle (1.5pt) (0,.15) circle (1.5pt);
\end{tikzpicture}}
\newcommand{\cirtdp}{\begin{tikzpicture} \useasboundingbox (-1,-.8) rectangle (1,1);
\draw [->,dash pattern=on 2pt off 2pt]  (-.8,0) arc (-180:180:.8);
\draw (-135:.8) -- (0,-.15) -- (0,.15) -- (150:.8)  (0,.15) -- (30:.8) (0,-.15) --  (-45:.8);
\draw (60:.8) -- (0,.45) -- (120:.8) (0,.45) -- (0,.8);
\fill (-135:.8) circle (1.5pt) (150:.8) circle (1.5pt) (30:.8) circle (1.5pt) (-45:.8) circle (1.5pt) (60:.8) circle (1.5pt) (90:.8) circle (1.5pt) (120:.8) circle (1.5pt) (0,.45) circle (1.5pt) (0,-.15) circle (1.5pt) (0,.15) circle (1.5pt);
\end{tikzpicture}}
\newcommand{\haltere}{\begin{tikzpicture} \useasboundingbox (-1.2,-.1) rectangle (1.2,.2);
\draw [->,dash pattern=on 2pt off 2pt]  (-.5,0) arc (-180:180:.15);
\draw [->,dash pattern=on 2pt off 2pt]  (.5,0) arc (0:360:.15);
\draw (-.2,0) -- (.2,0);
\fill (-.2,0) circle (1.5pt) (.2,0) circle (1.5pt) (.85,0) node{$S^1_j$} (-.85,0) node{$S^1_i$};
\end{tikzpicture}}

\newcommand{\onechordsmallj}{\begin{tikzpicture} \useasboundingbox (-.3,-.1) rectangle (.8,.2);
\draw [->,dash pattern=on 2pt off 2pt]  (0,.2) arc (-270:90:.2);
\draw (-.2,0) -- (.2,0) (.5,0) node{$S^1_j$};
\fill (-.2,0) circle (1.5pt) (.2,0) circle (1.5pt);
\end{tikzpicture}}
\newcommand{\onechordsmalljnum}{\begin{tikzpicture} \useasboundingbox (-.3,-.1) rectangle (.8,.2);
\draw [->,dash pattern=on 2pt off 2pt]  (0,-.2) arc (-90:270:.2);
\draw [->] (-.2,0) -- (.2,0);
\draw  (.5,0) node{$S^1_j$} (0,.1) node{\tiny $k$};
\fill (-.2,0) circle (1.5pt) (.2,0) circle (1.5pt);
\end{tikzpicture}}
\newcommand{\loopedge}{\begin{tikzpicture} \useasboundingbox (-.2,-.1) rectangle (.6,.2);
\draw  (.5,0) arc (0:360:.15);
\draw (-.1,0) -- (.2,0);
\fill (.2,0) circle (1.5pt);
\end{tikzpicture}}
\makeindex{N}
\begin{document} 
\title{An introduction to finite type invariants of knots and $3$-manifolds defined by counting graph configurations\\}

\author{Christine Lescop \thanks{Institut Fourier, UJF Grenoble, CNRS}}
\maketitle
\begin{abstract}
The finite type invariant concept for knots was introduced in the 90's in order to classify knot invariants, with the work of Vassiliev, Goussarov and Bar-Natan, shortly after the birth of numerous quantum knot invariants. This very useful concept was extended to $3$--manifold invariants by Ohtsuki.

These introductory lectures show how to define finite type invariants of links and $3$-manifolds by counting graph configurations in $3$-manifolds, following ideas of Witten and Kontsevich.

The linking number is the simplest finite type invariant for $2$--component links. It is defined in many equivalent ways in the first section. As an important example, we present it as the algebraic intersection of a torus and a $4$-chain called a {\em propagator\/} in a configuration space.

In the second section, we introduce the simplest finite type $3$--manifold invariant, which is the Casson invariant (or the $\Theta$--invariant) of integer homology $3$--spheres. It is defined as the algebraic intersection of three propagators in the same two-point configuration space.

In the third section, we explain the general notion of finite type invariants and introduce relevant spaces of Feynman Jacobi diagrams.

In Sections 4 and 5, we sketch an original construction based on configuration space integrals of universal finite type invariants for links in rational homology $3$--spheres and we state open problems. Our construction generalizes the known constructions for links in $\RR^3$ and for rational homology $3$--spheres, and it makes them more flexible.

In Section 6, we present the needed properties of parallelizations of $3$--manifolds and associated Pontrjagin classes, in details.

\end{abstract}


\maketitle

\tableofcontents
\section*{Foreword}
These notes contain some details about talks that were presented in the international conference ``Quantum Topology'' organized by Laboratory of Quantum Topology of Chelyabinsk State University in July 2014. They are based on the notes of five lectures presented in the ICPAM-ICTP research school of Mekn\`es in May 2012.
I thank the organizers of these two great events. I also thank Catherine Gille and K\'evin Corbineau for useful comments on these notes.

These notes have been written in an introductory way, in order to be understandable by graduate students. In particular, Sections~\ref{seclk}, \ref{secTheta} and \ref{secmorepar} provide an elementary self-contained presentation of the $\Theta$--invariant. The notes also contain original statements (Theorems~\ref{thmmain},\ref{thmunivone}, \ref{thmkeyuniv} and \ref{thmmainstraight}) together with sketches of proofs. Complete proofs of these statements, which generalize known statements, will be included in a monograph \cite{lesbookz}.

\section{Various aspects of the linking number}
\label{seclk}

\subsection{The Gauss linking number of two disjoint knots in \texorpdfstring{$\RR^3$}{the ambient space}}
\label{subdefgauss}

The modern powerful invariants of links and $3$--manifolds that will be defined in Section~\ref{secconstcsi} can be thought of as generalizations of the linking number. In this section, we warm up with several ways of defining this classical basic invariant. This allows us to introduce conventions and methods that will be useful througout the article.

Let $S^1$ denote the unit circle of $\CC$.
$$S^1=\{z;z \in \CC,|z|=1\}.$$
Consider two $C^{\infty}$ embeddings
$$J \colon S^1 \hookrightarrow \RR^3\;\;\;\mbox{and}\;\;\;K \colon S^1 \hookrightarrow \RR^3 \setminus J(S^1)$$

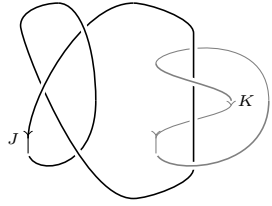
\begin{center}
\begin{tikzpicture}
\useasboundingbox (.5,.5) rectangle (4,3.5);
\draw [>-,draw=white,double=black,very thick] (.8,1.5) .. controls (.8,1) and (1.7,1) .. (1.7,2);
\draw [>-,draw=gray] (4,2) .. controls (4,3) and (2.5,2.7) .. (2.5,2.5) (3.5,2) .. controls (3.5,1.8) and (2.5,1.7) .. (2.5,1.5);
\draw [draw=white,double=black,very thick] (1.2,3.3) .. controls (1,3.3) and (.7,3.2) .. (.7,3) .. controls (.7,2.8) and (1.5,.7) .. (2.2,.7) .. controls (2.4,.7) and (3,.9) .. (3,1.1) -- (3,2.9) .. controls (3,3.1) and (2.4,3.3) .. (2.2,3.3);
\draw [>-,draw=white,double=gray,very thick] (2.5,2.5) .. controls (2.5,2.3) and (3.5,2.2) .. (3.5,2) (2.5,1.5) .. controls (2.5,1) and (4,1) .. (4,2);
\draw [draw=white,double=black,very thick] (2.2,3.3)  .. controls (1.5,3.3) and (.8,2) .. (.8,1.5);
\draw [draw=white,double=black,very thick] (1.7,2)  .. controls (1.7,2.7) and (1.5,3.3) .. (1.2,3.3);
\draw [>-] (.8,1.6)--(.8,1.3);
\draw [>-,draw=gray] (2.5,1.6)--(2.5,1.3);
\draw  (.6,1.5) node{\tiny $J$} (3.7,2) node{\tiny $K$};
\end{tikzpicture}
\end{center}
and the associated {\em Gauss map\/}
$$\begin{array}{llll}p_{JK}\colon & S^1 \times S^1 &\hookrightarrow & S^2\\
&(w,z) &\mapsto & \frac{1}{\parallel K(z)-J(w) \parallel}(K(z)-J(w)) \end{array}
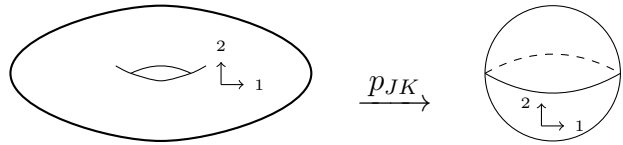
$$

\begin{center}
\begin{tikzpicture} \useasboundingbox (0,-.5) rectangle (4,1); 
\draw [thick] (0,0) .. controls (0,-.5) and (1.3,-.9) .. (2,-.9) .. controls (2.7,-.9) and (4,-.5) .. (4,0) .. controls (4,.5) and (2.7,.9) .. (2,.9) .. controls (1.3,.9) and (0,.5) .. (0,0); 
\draw plot[smooth] coordinates{(1.4,.1) (1.6,0) (2,-.1) (2.4,0) (2.6,.1)};
\draw [out=25,in=155] (1.6,0) to (2.4,0);
\draw [->] (2.8,-.15) -- (3.1,-.15) node[right]{\tiny $1$};
\draw [->] (2.8,-.15) -- (2.8,.15) node[above]{\tiny $2$};
\end{tikzpicture}
$\;\;\;\; \hfl{p_{JK}}\;\;\;\;$
\begin{tikzpicture} \useasboundingbox (-1,-.5) rectangle (1,1);
\draw (0,0) circle (.9);
\draw [out=-30,in=-150] (-.9,0) to (.9,0); 
\draw [out=29,in=151,dashed] (-.9,0) to (.9,0); 
\draw [->] (-.15,-.7) -- (.15,-.7) node[right]{\tiny $1$};
\draw [->] (-.15,-.7) -- (-.15,-.4) node[left]{\tiny $2$};
\end{tikzpicture}
\end{center}

\vspace{.5cm}

Denote the standard area form of $S^2$ by $4\pi \omega_{S^2}$ so that $\omega_{S^2}$ \index{N}{omegaS2@$\omega_{S^2}$} is the homogeneous volume form of $S^2$ such that
$\int_{S^2}\omega_{S^2}=1$.
In 1833, Gauss defined the {\em linking number\/} of the disjoint {\em knots\/} $J(S^1)$ and $K(S^1)$, simply denoted by $J$ and $K$, as
an integral \cite{gauss}.
 With modern notation, his definition reads
$$lk_G(J,K)=\int_{S^1 \times S^1}p_{JK}^{\ast}(\omega_{S^2}).$$
It can be rephrased as {\em $lk_G(J,K)$ is the degree of the Gauss map $p_{JK}$\/}.

\subsection{Some background material on manifolds without boundary, orientations, and degree}
\label{subbackground}
A {\em topological $n$--dimensional manifold $M$ without boundary\/} is a Hausdorff topological space that is a countable union of open subsets $U_i$ labeled in a set $I$
($i \in I$), where every $U_i$ is identified with an open subset $V_i$ of $\RR^n$ by a homeomorphism $\phi_i: U_i \rightarrow V_i$, called a {\em chart.\/} Manifolds are considered up to homeomorphism so that homeomorphic manifolds are considered identical.

For $r=0,\dots,\infty$, the topological manifold $M$ {\em  has a $C^r$--structure\/} or {\em is a $C^r$--manifold\/}, if, for each
pair $\{i,j\} \subset I$, the map $\phi_j \circ \phi_i^{-1}$ defined on $\phi_i(U_i \cap U_j)$ is a $C^r$--diffeomorphism to its image. The notion
of $C^s$--maps, $s \leq r$, from such a manifold to another one can be naturally deduced from the known case where the manifolds are open subsets of some  $\RR^n$, thanks to the local identifications provided by the charts. $C^r$--manifolds are considered up to $C^r$--diffeomorphisms.

An {\em orientation\/} of a real vector space $V$ of positive dimension is a basis of $V$ up to a change of basis with positive determinant. When $V=\{0\}$, an orientation of $V$ is an element of $\{-1,1\}$.
For $n>0$, an orientation of $\RR^n$ identifies $H_n(\RR^n,\RR^n\setminus\{x\};\RR)$ with $\RR$.
(In these notes, we freely use basic algebraic topology, see \cite{greenberg} for example.)
A homeomorphism $h$ from an open subset $U$ of $\RR^n$ to another such $V$ is {\em orientation-preserving\/} at a point $x$, if 
$h_{\ast} \colon H_n(U,U\setminus\{x\}) \rightarrow H_n(V,V\setminus\{h(x)\})$ is orientation-preserving.
If $h$ is a diffeomorphism, $h$ is orientation-preserving at $x$ if and only if the determinant of the Jacobian $\tang_xh$ is positive.
If $\RR^n$ is oriented and if the transition maps $\phi_j \circ \phi_i^{-1}$
are orientation-preserving (at every point) for $\{i,j\} \subset I$, the manifold $M$ is {\em oriented}.

For $n=0$, $1$, $2$ or $3$, any topological $n$-manifold may be equipped with a unique smooth structure (up to diffeomorphism) (See Theorem~\ref{thmstructhree}, below). Unless otherwise mentioned, our manifolds are {\em smooth \/} (i. e. $C^{\infty}$), oriented and compact, and considered up oriented diffeomorphisms.
Products are oriented by the order of the factors. More generally, unless otherwise mentioned, the order of appearance of coordinates or parameters orients manifolds.

A point $y$ is {\em a regular value \/} of a smooth map $p \colon M \rightarrow N$ between two smooth manifolds $M$ and $N$, if for any $x \in p^{-1}(y)$ the tangent map $\tang_xp$ at $x$ is surjective. According to the Morse-Sard theorem \cite[p. 69]{hirsch}, the set of regular values of such a map is dense.
If $M$ is compact, it is furthermore open. 

When $M$ is oriented and compact, and when the dimension of $M$ coincides with the dimension of $N$, the {\em differential degree\/} of $p$ at a regular value $y$ of $N$ is the (finite) sum running over the $x \in p^{-1}(y)$ of the signs of the determinants of $\tang_xp$.
In our case where $M$ has no boundary, this differential degree is locally constant on the set of regular values, and it is the {\em degree\/} of $p$, if $N$ is connected. See \cite[Chapter 5]{Mil}.

Finally, recall a homological definition of the degree.
Let $[M]$ denote the class of an oriented {\em closed\/} (i.e. compact, connected, without boundary) $n$--manifold in $H_n(M;\ZZ)$. $H_n(M;\ZZ)=\ZZ[M]$. If $M$ and $N$ are two closed oriented $n$--manifolds and if $f \colon M \rightarrow N$ is a (continuous) map, then $H_n(f)([M])=\mbox{deg}(f)[N]$. In particular, for the Gauss map $p_{JK}$ of Subsection~\ref{subdefgauss}, $$H_2(p_{JK})([S^1\times S^1])=lk(J,K)[S^2].$$

\subsection{The Gauss linking number as a degree}

Since the differential degree of the Gauss map $p_{JK}$ is locally constant, $lk_G(J,K)=\int_{S^1 \times S^1}p_{JK}^{\ast}(\omega)$ for any $2$-form $\omega$ on $S^2$ such that $\int_{S^2}\omega=1$.

Let us compute $lk_G(J,K)$ as the differential degree of $p_{JK}$ at the vector $Y$ that points towards us. The set $p_{JK}^{-1}(Y)$ is made of the pairs of points $(w,z)$ where the projections of $J(w)$ and $K(z)$ coincide, and $J(w)$ is under $K(z)$. They correspond to the  {\em crossings\/} \pcbg and \ncbg of the diagram.

In a diagram, a crossing is {\em positive\/} if we turn counterclockwise from the arrow at the end of the upper strand to the arrow of the end of the lower strand like \pcortrig.
Otherwise, it is {\em negative\/} like \nc.

For the positive crossing \pcbg, moving $J(w)$ along $J$ following the orientation of $J$, moves $p_{JK}(w,z)$ towards the South-East direction, while moving $K(z)$ along $K$ following the orientation of $K$, moves $p_{JK}(w,z)$ towards the North-East direction, so that the local orientation induced by the image of $p_{JK}$ around $Y \in S^2$ is \begin{tikzpicture}
\useasboundingbox (-.1,.1) rectangle (1.5,.55);
\draw [->] (0,.2) -- (.2,0) node[right]{\tiny $\tang p\frac{\partial}{\partial w}$};
\draw [->,]  (0,.2) -- (.2,.4) node[right]{\tiny $\tang p\frac{\partial}{\partial z}$};
\end{tikzpicture}, which is \begin{tikzpicture}
\useasboundingbox (-.1,.1) rectangle (.5,.55);
\draw [->] (0,.2) -- (.2,0) node[right]{\tiny $1$};
\draw [->,]  (0,.2) -- (.2,.4) node[right]{\tiny $2$};
\end{tikzpicture}.
Therefore, the contribution of a positive crossing to the degree is $1$. Similarly, the contribution of a negative crossing is $(-1)$.

We have just proved the following formula $$\mbox{deg}_Y(p_{JK})= \sharp \pcbg - \sharp \ncbg$$
where $\sharp$ stands for the cardinality --here $\sharp \pcbg$ is the number of occurences of $\pcbg$ in the diagram--
so that
$$lk_G(J,K)= \sharp \pcbg - \sharp \ncbg .$$
Similarly, $\mbox{deg}_{-Y}(p_{JK})= \sharp \pcbgkj - \sharp \ncbgkj$
so that
$$lk_G(J,K)= \sharp \pcbgkj - \sharp \ncbgkj = \frac12\left(\sharp \pcbg +\sharp \pcbgkj - \sharp \ncbg - \sharp \ncbgkj\right)$$
and $lk_G(J,K)=lk_G(K,J)$.

In our first example, $lk_G(J,K)=2$. Let us draw some further examples.

For the {\em positive Hopf link \/} \begin{tikzpicture}[scale=.5] \useasboundingbox (-1,-.5) rectangle (1,.5);
\draw[->] (0,0) -- (0,.2) .. controls (0,.4) and (-.15,.5).. (-.3,.5) .. controls(-.45,.5) and (-.6,.2) .. (-.6,0) node[left]{\tiny $J$};
\draw[->] (-.6,0) .. controls (-.6,-.2) and (-.45,-.5) .. (-.3,-.5) .. controls (-.15,-.5) and (0,-.35) .. (0,-.2);
\draw[->,draw=gray] (-.1,.2) .. controls (-.2,.2) and (-.25,.15) .. (-.25,.05) .. controls (-.25,-.05) and (-.2,-.1) .. (0,-.1) .. controls (.2,-.1) and (.25,-.05) .. (.25,.05) .. controls (.25,.15) and (.2,.2) .. (.1,.2) node[right]{\tiny $K$};
\end{tikzpicture}, $lk_G(J,K)=1$.

For the {\em negative Hopf link \/} \begin{tikzpicture}[x=-30,scale=.5] \useasboundingbox (-1,-.5) rectangle (1,.5);
\draw[->] (0,0) -- (0,.2) .. controls (0,.4) and (-.15,.5).. (-.3,.5) .. controls(-.45,.5) and (-.6,.2) .. (-.6,0);
\draw[->] (-.6,0) .. controls (-.6,-.2) and (-.45,-.5) .. (-.3,-.5) .. controls (-.15,-.5) and (0,-.35) .. (0,-.2);
\draw[->,draw=gray] (-.1,.2) .. controls (-.2,.2) and (-.25,.15) .. (-.25,.05) .. controls (-.25,-.05) and (-.2,-.1) .. (0,-.1) .. controls (.2,-.1) and (.25,-.05) .. (.25,.05) .. controls (.25,.15) and (.2,.2) .. (.1,.2);
\end{tikzpicture}, $lk_G(J,K)=-1$.

For the {\em Whitehead link \/}
\begin{tikzpicture}
\useasboundingbox (-.5,-1.1) rectangle (1.5,1);
\draw (.4,0) -- (.4,.2) .. controls (.4,.45) and (.6,.6).. (.8,.6) .. controls (1,.6) and (1.4,.2) .. (1.4,-.1);
\draw (0,0) -- (0,.2) .. controls (0,.6) and (.5,.9) .. (.8,.9) .. controls (1.1,.9) and (1.7,.5) .. (1.7,-.1);
\draw [->,draw=gray] (.3,.2) -- (.1,.2) (-.1,.2) .. controls (-.2,.2) and (-.35,.15) .. (-.35,.05) .. controls (-.35,-.05) and (0,-.1) .. (.2,-.1) .. controls (.4,-.1) and (.75,-.05) .. (.75,.05) .. controls (.75,.15) and (.6,.2) .. (.5,.2);
\draw[->] (.2,-.3) .. controls (-.1,-.3) and (.2,-1.1) .. (.8,-1.1) .. controls (1.1,-1.1) and (1.7,-.7) .. (1.7,-.1);
\draw [draw=white,double=black,very thick] (0,-.2) .. controls (0,-.3) and (.1,-.5) .. (.2,-.5) .. controls (.3,-.5) and (.4,-.3) ..  (.4,-.2);
\draw [draw=white,double=black,very thick]  (.2,-.3) .. controls (.3,-.3) and (.5,-.8) .. (.8,-.8) .. controls (1.1,-.8) and (1.4,-.4) .. (1.4,-.1);
\end{tikzpicture}
, $lk_G(J,K)=0$.

\subsection{Some background material on manifolds with boundary and algebraic intersections}

A {\em topological $n$--dimensional manifold $M$ with possible boundary\/} is a Hausdorff topological space that is a union of open subsets $U_i$ with subscripts in a set $I$, 
($i \in I$), where every $U_i$ is identified with an open subset $V_i$ of $]-\infty,0]^k \times \RR^{n-k}$ by a chart $\phi_i: U_i \rightarrow V_i$. The {\em boundary\/} of $]-\infty,0]^k \times \RR^{n-k}$ is made of the points $(x_1,\dots,x_n)$ of $]-\infty,0]^k \times \RR^{n-k}$ such that there exists $i \leq k$ such that $x_i=0$.
 The {\em boundary\/} of $M$ is made of the points that are mapped to the boundary of $]-\infty,0]^k \times \RR^{n-k}$.

For $r=1,\dots,\infty$, the topological manifold $M$ {\em is a $C^r$--manifold with ridges (or with corners) (resp. with boundary)\/}, if, for each
pair $\{i,j\} \subset I$, the map $\phi_j \circ \phi_i^{-1}$ defined on $\phi_i(U_i \cap U_j)$ is a $C^r$--diffeomorphism to its image (resp. and if furthermore $k\leq 1$, for any $i$).
Then the {\em ridges\/} of $M$ are made of the points that are mapped to points $(x_1,\dots,x_n)$ of $]-\infty,0]^k \times \RR^{n-k}$ so that there are at least two $i \leq k$ such that $x_i=0$.

The tangent bundle to an oriented submanifold $A$ in a manifold $M$ at a point $x$ is denoted by $T_xA$. The {\em normal bundle\/} $T_xM/T_xA$ to $A$ in $M$ at $x$ is denoted by ${\normbun}_xA$. It is oriented so that (a lift of an oriented basis of) ${\normbun}_xA$ followed by (an oriented basis of) $T_xA$ induce the orientation of $T_xM$.
The boundary $\partial M$ of an oriented manifold $M$ is oriented by the {\em outward normal first\/} convention. If $x \in \partial M$ is not in a ridge, the outward normal to $M$ at $x$ followed by an oriented basis of $T_x \partial M$ induce the orientation of $M$. 
For example, the standard orientation of the disk in the plane induces the standard orientation of the circle, counterclockwise, as the following picture shows.
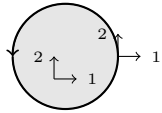
\begin{center}
\begin{tikzpicture}
\useasboundingbox (-.8,-.8) rectangle (.8,.8);
\draw [fill=gray!20] (0,0) circle (.7);
\draw [->] (-.15,-.3) -- (.15,-.3) node[right]{\tiny $1$};
\draw [->] (-.15,-.3) -- (-.15,0) node[left]{\tiny $2$};
\draw [->] (.7,0) -- (1,0) node[right]{\tiny $1$};
\draw [->] (.7,0) -- (.7,.3) node[left]{\tiny $2$};
\draw [thick,->] (-.7,0) arc (-180:180:.7);
\end{tikzpicture}
\end{center}
As another example, the sphere $S^2$ is oriented as the boundary of the ball $B^3$, which has the standard orientation
induced by (Thumb, index finger (2), middle finger (3)) of the right hand.
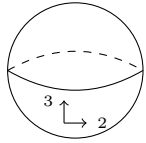
\begin{center}
\begin{tikzpicture} \useasboundingbox (-1,-.7) rectangle (1,1);
\draw (0,0) circle (.9);
\draw [out=-30,in=-150] (-.9,0) to (.9,0); 
\draw [out=29,in=151,dashed] (-.9,0) to (.9,0); 
\draw [->] (-.15,-.7) -- (.15,-.7) node[right]{\tiny $2$};
\draw [->] (-.15,-.7) -- (-.15,-.4) node[left]{\tiny $3$};
\end{tikzpicture}
\end{center}

Two submanifolds $A$ and $B$ in a manifold $M$ are transverse if at each intersection point $x$, $T_xM=T_xA+T_xB$.
The transverse intersection of two submanifolds $A$ and $B$ in a manifold $M$ is oriented so that the normal bundle to $A\cap B$ is $({\normbun}(A) \oplus {\normbun}(B))$, fiberwise. If the two manifolds are of complementary dimensions, then the sign of an intersection point is $+1$ if the orientation of its normal bundle coincides with the orientation of the ambient space, that is if $T_xM={\normbun}_xA \oplus {\normbun}_xB$ (as oriented vector spaces), this is equivalent to $T_xM=T_xA \oplus T_xB$ (as oriented vector spaces again, exercise). Otherwise, the sign is $-1$. If $A$ and $B$ are compact and if $A$ and $B$ are of complementary dimensions in $M$, their {\em algebraic intersection\/} is the sum of the signs of the intersection points, it is denoted by
$\langle A, B \rangle_M$.

When $M$ is an oriented manifold, $(-M)$ denotes the same manifold, equipped with the opposite orientation.
In a manifold $M$, a $k$-dimensional {\em chain (resp. rational chain)\/} is a finite combination with coefficients in $\ZZ$ (resp. in $\QQ$) of smooth $k$-dimensional oriented submanifolds $C$ of $M$ with boundary and ridges, up to the identification of $(-1)C$ with $(-C)$.

Again, unless otherwise mentioned, manifold are oriented.
The boundary $\partial$ of chains is a linear map that maps a smooth submanifold to its oriented boundary.  The canonical orientation of a point is the sign $+1$ so that $\partial [0,1]= \{1\}-\{0\}$.

\begin{lemma}
\label{lemsignint}
 Let $A$ and $B$ be two transverse submanifolds of a $d$--dimensional manifold $M$, of respective dimensions $\alpha$ and $\beta$, with disjoint boundaries. 
Then $$\partial (A\cap B)=(-1)^{d-\beta}\partial A \cap B + A \cap \partial B.$$
\end{lemma}
\bp Note that $\partial (A\cap B) \subset \partial A \cup \partial B$. At a point $a \in \partial A$, $T_aM$ is oriented by $({\normbun}_aA,o,T_a\partial A)$, where $o$ is the outward normal to $A$. If $a \in \partial A \cap B$, then $o$ is also an outward normal to $A\cap B$, and $\partial (A \cap B)$ is cooriented by $({\normbun}_aA,{\normbun}_aB,o)$ while $\partial A \cap B$ is cooriented by $({\normbun}_aA,o,{\normbun}_aB)$.
At a point $b \in A \cap \partial B$, $\partial (A \cap B)$ is cooriented by $({\normbun}_aA,{\normbun}_aB,o)$ like $A \cap \partial  B$.
\eop

\subsection{A general definition of the linking number}

\begin{lemma}
Let $J$ and $K$ be two rationally null-homologous disjoint cycles of respective dimensions $j$ and $k$ in a $d$-manifold $M$, where $d=j+k+1$.
There exists a rational $(j+1)$--chain $\Sigma_J$ bounded by $J$ transverse to $K$, and a rational $(k+1)$--chain $\Sigma_K$ bounded by $K$ transverse to $J$ and for any two such rational chains $\Sigma_J$ and $\Sigma_K$, $\langle J,\Sigma_K\rangle_M=(-1)^{j+1}\langle \Sigma_J,K\rangle_M $.
In particular, $\langle J,\Sigma_K\rangle_M$ is a topological invariant of $(J,K)$, which is denoted by $lk(J,K)$ and called the {\em linking number\/} of $J$ and $K$. 
$$lk(J,K)=(-1)^{(j+1)(k+1)}lk(K,J).$$
\end{lemma}
\bp
Since $K$ is rationally null-homologous, $K$ bounds a rational $(k+1)$--chain $\Sigma_K$. Without loss, $\Sigma_K$ is assumed to be transverse to $\Sigma_J$ so that $\Sigma_J \cap \Sigma_K$ is a rational $1$--chain (which is a rational combination of circles and intervals).
(As explained in \cite[Chapter 3]{hirsch}, generically, manifolds are transverse.) According to Lemma~\ref{lemsignint},
$$\partial (\Sigma_J \cap \Sigma_K) = (-1)^{d+k+1}J \cap \Sigma_K +\Sigma_J \cap K.$$
Furthermore, the sum of the coefficients of the points in the left-hand side must be zero, since this sum vanishes for the boundary of an interval. This shows that $\langle J,\Sigma_K\rangle_M=(-1)^{d+k} \langle \Sigma_J,K\rangle_M$, and therefore that this rational number is independent of the chosen $\Sigma_J$ and $\Sigma_K$. Since $(-1)^{d+k}\langle \Sigma_J,K\rangle_M=(-1)^{j+1}(-1)^{k(j+1)}\langle K, \Sigma_J\rangle_M $, $lk(J,K)=(-1)^{(j+1)(k+1)}lk(K,J).$
\eop

In particular, the {\em linking number\/} of two rationally null-homologous disjoint links $J$ and $K$ in a $3$-manifold $M$ is the algebraic intersection of a rational chain bounded by one of the knots and the other one.

For $\KK=\ZZ$ or $\QQ$, a {\em $\KK$-sphere or (integer or rational) homology $3$-sphere\/} (resp. a {\em $\KK$-ball\/}) is a smooth, compact, oriented $3$-manifold, without ridges, with the same $\KK$-homology as the sphere $S^3$ (resp. as a point).
In such a manifold, any knot is rationally null-homologous so that the linking number of two disjoint knots always makes sense.

A {\em meridian \/} of a knot $K$ is the (oriented) boundary of a disk that intersects $K$ once with a positive sign.
Since a chain $\Sigma_J$ bounded by a knot $J$ disjoint from $K$ in a $3$--manifold $M$ provides a rational cobordism between $J$ and a combination of meridians of $K$,
 $[J]=lk(J,K)[m_K]$ in $H_1(M\setminus K;\QQ)$ where $m_K$ is a meridian of $K$.

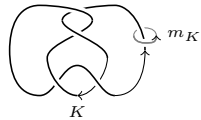
\begin{center}
\begin{tikzpicture}
\useasboundingbox (-1,-.3) rectangle (1,1.3);
\draw[->] (.5,0) .. controls (.8,0) and (.9,.3) .. (.9,.6);
\draw [draw=gray] (1.05,.8)  .. controls (1.05,.85) and (1,.9) .. (.9,.9) .. controls (.8,.9) and (.75,.85) .. (.75,.8); 
\draw [draw=white,double=black,very thick] (.9,.6) .. controls (.9,.9) and (.7,1.2) .. (.5,1.2) .. controls (.3,1.2) and (-.2,1.15) .. (-.2,1);
\draw [draw=white,double=gray,very thick] (.75,.8)  .. controls (.75,.75) and (.8,.7) .. (.9,.7) .. controls (1,.7) and (1.05,.7) .. (1.05,.8);
\draw [->,draw=gray] (.75,.8)  .. controls (.75,.75) and (.8,.7) .. (.9,.7) .. controls (1,.7) and (1.05,.7) .. (1.05,.8) node[right]{\tiny $m_K$} ; 
\draw[->] (.4,.5) .. controls (.4,.3) and (.2,0) .. (0,0) node[below]{\tiny $K$} ;
\draw [draw=white,double=black,very thick] (.2,1) .. controls (.2,1.15) and (-.3,1.2) .. (-.5,1.2) .. controls (-.8,1.2) and (-.9,.9) .. (-.9,.6) .. controls (-.9,.3) and (-.8,0) .. (-.5,0) .. controls (-.3,0) and (-.2,.4) .. (0,.4) .. controls (.2,.4) and (.3,0) .. (.5,0);
\draw[draw=white,double=black,very thick] (0,0) .. controls (-.2,0) and (-.4,.3) .. (-.4,.5) .. controls (-.4,.7) and (.2,.85) .. (.2,1);
\draw[draw=white,double=black,very thick] (-.2,1) .. controls (-.2,.85) and (.4,.7) .. (.4,.5);
\end{tikzpicture}
\end{center}

\begin{lemma}
\label{lemdeflkhom}
 When $K$ is a knot in a $\QQ$-sphere or a $\QQ$-ball $M$, $H_1(M\setminus K;\QQ)=\QQ[m_K]$, so that the equation $[J]=lk(J,K)[m_K]$ in $H_1(M\setminus K;\QQ)$ provides an alternative definition for the linking number.
\end{lemma}
\bp Exercise. \eop

The reader is also invited to check that $lk_G=lk$ as an exercise though it will be proved in the next subsection, see Proposition~\ref{propdeflkeq}.

\subsection{Generalizing the Gauss definition of the linking number and identifying the definitions}
\label{subgenGaussdef}

\begin{lemma}
\label{lempstt}
The map $$\begin{array}{llll}p_{S^2}\colon &((\RR^3)^2 \setminus \mbox{diag})&\rightarrow & S^2\\&(x,y) &\mapsto &\frac{1}{\parallel y-x \parallel}(y-x) \end{array}$$
is a homotopy equivalence.
In particular $$H_{i}(p_{S^2}) \colon H_i((\RR^3)^2 \setminus \mbox{diag};\ZZ) \rightarrow H_i(S^2;\ZZ)$$ is an isomorphism for all $i$,
$((\RR^3)^2 \setminus \mbox{diag})$ is a homology $S^2$, and $[S]=\left(H_{2}(p_{S^2})\right)^{-1}[S^2]$ is a canonical generator of $$H_2((\RR^3)^2 \setminus \mbox{diag};\ZZ)=\ZZ[S].$$
\end{lemma}
\bp $((\RR^3)^2 \setminus \mbox{diag})$ is homeomorphic to $\RR^3 \times ]0,\infty[ \times S^2$ via the map $$(x,y) \mapsto (x,\parallel y-x \parallel, p_{S^2}(x,y)).$$
\eop

As in Subsection~\ref{subdefgauss}, consider a two-component link
$J \sqcup  K: S^1 \sqcup S^1 \hookrightarrow \RR^3$. 
This embedding induces an embedding
$$\begin{array}{llll}J \times  K\colon & S^1 \times S^1 &\hookrightarrow &((\RR^3)^2 \setminus \mbox{diag})\\
&(z_1,z_2) &\mapsto & (J(z_1),K(z_2))\end{array}$$
the map $p_{JK}$ of Subsection~\ref{subdefgauss} reads $p_{S^2} \circ (J \times  K)$,
and since $H_2(p_{JK})[S^1 \times S^1]=\mbox{deg}(p_{JK})[S^2]=lk_G(J,K)[S^2]$ in $H_2(S^2;\ZZ)=\ZZ[S^2]$, 
$$[J\times K]=H_2(J \times  K)[S^1 \times S^1]=lk_G(J,K)[S]$$
in $H_2((\RR^3)^2 \setminus \mbox{diag};\ZZ)=\ZZ[S].$
We will see that this definition of $lk_G$ generalizes to links in rational homology $3$--spheres
and then prove that our generalized definition coincides with the general definition of linking numbers in this case.

For a $3$-manifold $M$, the normal bundle to the diagonal of $M^2$ in $M^2$ is identified with the tangent bundle to $M$, fiberwise, by the map
$$(u,v) \in \frac{(T_xM)^2}{\mbox{diag}((T_xM)^2)} \mapsto (v-u) \in T_xM.$$

A {\em parallelization\/} $\tau$ of an oriented $3$-manifold $M$ is a bundle isomorphism $\tau \colon M \times \RR^3 \longrightarrow TM$ that restricts to $ x\times \RR^3$ as an orientation-preserving linear isomorphism from $ x\times \RR^3$ to $T_xM$, for any $x \in M$.
It has long been known that any oriented $3$-manifold is parallelizable (i.e. admits a parallelization). (It is proved in Subsection~\ref{subproofpar}.)
Therefore, a tubular neighborhood of the diagonal in $M^2$ is diffeomorphic to $M \times \RR^3$.

\begin{lemma}
\label{lemhomstwo}
Let $M$ be a rational homology $3$--sphere, let $\infty$ be a point of $M$. Let $\check{M}=(M \setminus \{\infty\})$. Then $\check{M}^2 \setminus \mbox{diag}$ has the same rational homology as $S^2$. Let $B$ be a ball in $\check{M}$ and let $x$ be a point inside $B$, then
the class $[S]$ of $x \times \partial B$ is a canonical generator of $H_2(\check{M}^2 \setminus \mbox{diag};\QQ)=\QQ[S]$.
\end{lemma}
\bp
In this proof, the homology coefficients are in $\QQ$. 
Since $\check{M}$ has the homology of a point, the K\"unneth Formula implies that
 $\check{M}^2$ has the homology of a point. 
Now, by excision,
$$H_{\ast}(\check{M}^2,\check{M}^2 \setminus \mbox{diag}) \cong H_{\ast}(\check{M} \times \RR^3,\check{M}\times (\RR^3 \setminus 0))$$
$$ \cong H_{\ast}( \RR^3, S^2) \cong \left\{\begin{array}{ll} \QQ \;\;&\;\mbox{if} \;\ast =3,\\
0\;&\;\mbox{otherwise.} \end{array} \right.$$
Using the long exact sequence of the pair $(\check{M}^2,\check{M}^2 \setminus \mbox{diag})$, we get that
$H_{\ast}(\check{M}^2 \setminus \mbox{diag};\QQ)=H_{\ast}(S^2)$.
\eop

Define the {\em Gauss linking number\/} of two disjoint links $J$ and $K$ in $\check{M}$ so that $$[(J\times K)(S^1\times S^1)]=lk_G(J,K)[S]$$ in $H_2(\check{M}^2 \setminus \mbox{diag};\QQ)$. Note that the two definitions of $lk_G$ coincide when $\check{M}=\RR^3$.
\begin{proposition}
\label{propdeflkeq}
 $$lk_G=lk$$
\end{proposition}
\bp
First note that both definitions make sense when $J$ and $K$ are disjoint links: $[J\times K]=lk_G(J,K)[S]$ and
$lk(J,K)$ is the algebraic intersection of $K$ and a rational chain $\Sigma_J$ bounded by $J$.

If $K$ is a knot, then the chain $\Sigma_J$ of $\check{M}$ provides a rational cobordism $C$ between $J$ and a combination of meridians of $K$ in $\check{M}\setminus K$, and a rational cobordism $C\times K$ in $\check{M}^2 \setminus \mbox{diag}$, which allow us to see that $lk_G(.,K)$ and $lk(.,K)$ linearly depend on $[J] \in H_1(\check{M}\setminus K)$. Thus we are left with the proof that $lk_G(m_K,K)=lk(m_K,K)=1$. Since $lk_G(m_K,.)$ linearly depends on $[K] \in H_1(\check{M}\setminus m_K)$, we are left with the proof $lk_G(m_K,K)=1$ when $K$ is a meridian of $m_K$. Now, there is no loss in assuming that our link is a Hopf link in $\RR^3$ so that the equality follows from the equality for the positive Hopf link in $\RR^3$.
\eop

For a $2$--component link $(J,K)$ in $\RR^3$, the definition of $lk(J,K)$ can be rewritten as
$$lk(J,K)=\int_{J \times K}p_{S^2}^{\ast}(\omega)=\langle J \times K, p_{S^2}^{-1}(Y)\rangle_{(\RR^3)^2 \setminus \mbox{diag}}$$ for any regular value $Y$ of $p_{JK}$, and for any $2$-form $\omega$ of $S^2$ such that $\int_{S^2}\omega=1$.
Thus, $lk(J,K)$ is the evaluation of a $2$-form $p_{S^2}^{\ast}(\omega)$ of $(\RR^3)^2 \setminus \mbox{diag}$ at the $2$-cycle $[J \times K]$, or it is the intersection of the $2$-cycle $[J \times K]$ with a $4$-manifold $p_{S^2}^{-1}(Y)$, which will later be seen as the interior of a prototypical propagator.
We will adapt these definitions to rational homology $3$--spheres in Subsection~\ref{subprop}. The definition of the linking number that we will generalize in order to produce more powerful invariants is contained in Lemma~\ref{lemlkprop}.

\newpage
\section{Propagators and the \texorpdfstring{$\Theta$}{Theta}-invariant}
\label{secTheta}
\setcounter{equation}{0}

Propagators will be the key ingredient to define powerful invariants from graph configurations
in Section~\ref{secconstcsi}. They are defined in Subsection~\ref{subprop} below after needed preliminaries.
They allow us to define the $\Theta$--invariant as an invariant of parallelized homology $3$--balls in Subsection~\ref{subThetapar}. The $\Theta$--invariant is next turned to an invariant of rational homology $3$--spheres in Subsection~\ref{subThetarat}.

\subsection{Blowing up in real differential topology}

Let $A$ be a submanifold of a smooth manifold $B$, and let $U{\normbun}(A)$ denote its unit normal bundle. The fiber $U{\normbun}_a(A)=({\normbun}_a(A) \setminus \{0\})/\RR^{+\ast}$ of $U{\normbun}(A)$ is oriented as the boundary of a unit ball of ${\normbun}_a(A)$.

Here, {\em blowing up\/} such a submanifold $A$ of codimension $c$ of $B$ means replacing $A$ by $U{\normbun}(A)$. For small open subspaces $U_A$ of $A$, $\left((\RR^c=\{0\} \cup ]0,\infty[ \times S^{c-1}) \times U_A\right)$ is replaced by $([0,\infty[ \times S^{c-1} \times U_A)$, so that the blown-up
manifold $\blowup{B}{A}$ is homeomorphic to the complement in $B$ of an open tubular neighborhood (thought of as infinitely small) of $A$. In particular, $\blowup{ B }{ A }$ is homotopy equivalent to $B \setminus A$.
Furthermore, the blow up is canonical, so that the created boundary is $\pm U{\normbun}(A)$ and there is a canonical smooth projection from $\blowup{B}{A}$ to $B$ such that the preimage of $a \in A$ is $U{\normbun}_a(A)$.
If $A$ and $B$ are compact, then $\blowup{B}{A}$ is compact, it is a smooth compactification of $B \setminus A$.

In the following figure, we see the result of blowing up $(0,0)$ in $\RR^2$, and the closures in $\blowup{\RR^2}{(0,0)}$ of $\{0\} \times \RR$, $\RR \times \{0\}$ and the diagonal of $\RR^2$, successively.

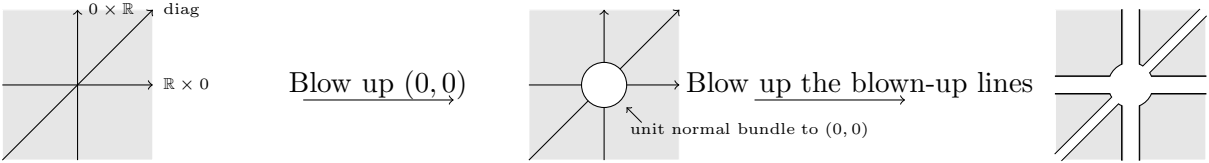
\begin{center}
\begin{tikzpicture}
\useasboundingbox (-2,-1) rectangle (14,1);
\begin{scope}[xshift=-1cm]
\draw [fill=gray!20, draw=white] (-1,-1) rectangle (1,1);
\draw [->] (-1,0) -- (1,0) node[right]{\tiny $\RR \times 0$};
\draw [->] (0,-1) -- (0,1) node[right]{\tiny $0 \times \RR$};
\draw [->] (-1,-1) -- (1,1) node[right]{\tiny diag};
\end{scope}
\draw (3,0) node{\small Blow up $(0,0)$} ;
\draw [->] (2,-.2) -- (4,-.2);
\begin{scope}[xshift=6cm]
\draw [fill=gray!20, draw=white] (-1,-1) rectangle (1,1);
\draw [->] (-1,0) -- (1,0);
\draw [->] (0,-1) -- (0,1);
\draw [->] (-1,-1) -- (1,1);
\draw [fill=white] (0,0) circle (.3);
\draw [->,very thin] (.5,-.5) -- (.3,-.3);
\draw (.2,-.6) node[right]{\tiny unit normal bundle to $(0,0)$};
\draw (3.4,0) node{\small Blow up the blown-up lines};
\draw [->] (2,-.2) -- (4,-.2);
\end{scope}
\begin{scope}[xshift=13cm]
\draw [fill=gray!20, draw=white] (-1,-1) rectangle (1,1);
\draw [very thick] (-.93,-1) -- (1,.93) (-1,-.93) -- (.93,1);
\draw [fill=black]  (-.12,-1) rectangle (.12,1) (-1,-.12) rectangle (1,.12);
\draw [fill=white] (0,0) circle (.3);
\draw [fill=white,draw=white] (-1,-1) -- (-.93,-1) -- (1,.93) -- (1,1) -- (.93,1) -- (-1,-.93) -- (-1,-1);
\draw [fill=white,draw=white]  (-.1,-1.1) rectangle (.1,1.1) (-1.1,-.1) rectangle (1.1,.1);
\end{scope}
\end{tikzpicture}
\end{center}

\subsection{The configuration space \texorpdfstring{$C_2(M)$}{C2(M)}}

See $S^3$ as $\RR^3 \cup \infty$ or as two copies of $\RR^3$ identified along $\RR^3 \setminus \{0\}$ by the (exceptionally orientation-reversing) diffeomorphism $x \mapsto x/\norm{x}^2$.
Then $\blowup{S^3}{\infty}=\RR^3 \cup S^2_{\infty}$ where the unit normal bundle $(-S^2_{\infty})$ to $\infty$ in $S^3$ is canonically diffeomorphic to $S^2$ via $p_{\infty} \colon  S^2_{\infty} \rightarrow S^2$, where $x \in S^2_{\infty}$ is the limit of a sequence of points of $\RR^3$  approaching $\infty$ along a line directed by $p_{\infty}(x) \in S^2$. $$\partial \blowup{S^3}{\infty}=S^2_{\infty}$$

Fix a rational homology $3$--sphere $M$, a point $\infty$ of $M$, and $\check{M}=M \setminus \{\infty\}$. 
Identify a neighborhood of $\infty$ in $M$ with the complement $\check{B}_{1,\infty}$ of the closed ball $B(1)$ of radius $1$ in $\RR^3$. Let $\check{B}_{2,\infty}$ be the complement of the closed ball ${B}(2)$ of radius $2$ in $\RR^3$, which is a smaller neighborhood of $\infty$ in $M$ via the understood identification.
Then $B_M=M \setminus \check{B}_{2,\infty}$ is a compact rational homology ball diffeomorphic to $\blowup{M}{\infty}$.

Define the {\em configuration space\/} $C_2(M)$ as the compact $6$--manifold with boundary and ridges obtained from $M^2$ by blowing up $(\infty,\infty)$, the closures in $\blowup{M^2}{(\infty,\infty)}$ of $\{\infty\} \times \check{M}$, $\check{M} \times \{\infty\}$ and the diagonal of $\check{M}^2$, successively.
Then $\partial C_2(M)$ contains the unit normal bundle $(\frac{T\check{M}^2}{\mbox{\scriptsize diag}} \setminus \{0\})/\RR^{+\ast}$ to the diagonal of $\check{M}^2$. This bundle is canonically isomorphic to the unit tangent bundle $U\check{M}$ to $\check{M}$ (again via the map $\left([(x,y)] \mapsto [y-x]\right)$).

\begin{lemma}
\label{lemstrucctwo}
Let $\check{C}_2(M)=\check{M}^2 \setminus \mbox{\rm diag}$.
The open manifold $C_2(M) \setminus \partial C_2(M)$ is $\check{C}_2(M)$ and the inclusion $\check{C}_2(M) \hookrightarrow C_2(M)$ is a homotopy equivalence. In particular, $C_2(M)$ is a compactification of $\check{C}_2(M)$ homotopy equivalent to $\check{C}_2(M)$.
The manifold $C_2(M)$ is a smooth compact $6$-dimensional manifold with boundary and ridges.
There is a canonical smooth projection $p_{M^2}\colon C_2(M) \rightarrow M^2$.
$$\partial C_2(M)=\pm p_{M^2}^{-1}(\infty,\infty) \cup  ( S^2_{\infty}\times \check{M}) \cup (- \check{M} \times S^2_{\infty}) \cup U\check{M}.$$
\end{lemma}
\bp Let $B_{1,\infty}$ be the complement of the open ball of radius one of $\RR^3$ in $S^3$. Blowing up $(\infty,\infty)$ in $B_{1,\infty}^2$ transforms a neighborhood of $(\infty,\infty)$
into the product $[0,1[ \times S^5$. Explicitly, there is a map
$$\begin{array}{llll}\psi \colon& [0,1[ \times S^5 & \rightarrow & \blowup{B_{1,\infty}^2}{(\infty,\infty)} \subset \blowup{M^2}{(\infty,\infty)}\\
  & (\lambda \in ]0,1[, (x \neq 0,y \neq 0) \in S^5 \subset (\RR^3)^2)  & \mapsto & (\frac1{\lambda\norm{x}^2}x,\frac1{\lambda\norm{y}^2}y)
  \end{array}$$
that is a diffeomorphism onto its image, which is a neighborhood of the preimage of $(\infty,\infty)$ under the blow-down map $\blowup{M^2}{(\infty,\infty)} \hfl{p_1} M^2$.
This neighborhood intersects $\infty \times \check{M}$, $\check{M} \times \infty$, and $\mbox{diag}(\check{M}^2)$  as $\psi(]0,1[ \times 0 \times S^2)$, $\psi(]0,1[ \times S^2\times 0)$ and $\psi(]0,1[ \times (S^5 \cap \mbox{diag}((\RR^3)^2)))$, respectively. 
In particular, the closures of $\infty \times \check{M}$, $\check{M} \times \infty$, and $\mbox{diag}(\check{M}^2)$ in $\blowup{M^2}{(\infty,\infty)}$ intersect the boundary $\psi(0 \times S^5)$ of $\blowup{M^2}{(\infty,\infty)}$ as three disjoint spheres in $S^5$, and they read
 $\infty \times \blowup{M}{\infty}$, 
$\blowup{M}{\infty} \times \infty$ and $\mbox{diag}(\blowup{M}{\infty}^2)$.
Thus, the next steps will be three blow-ups along these three disjoint smooth manifolds.

These blow-ups will preserve the product structure $\psi([0,1[ \times .)$.
Therefore, $C_2(M)$ is a smooth compact $6$-dimensional manifold with boundary, with three {\em ridges\/} $S^2 \times S^2$ in $p_{M^2}^{-1}(\infty,\infty)$. A neighborhood of these ridges in $C_2(M)$ is diffeomorphic to $[0,1[^2 \times S^2 \times S^2$.
\eop

\begin{lemma}
\label{lemextproj}
The map $p_{S^2}$ of Lemma~\ref{lempstt} smoothly extends to $C_2(S^3)$, and its extension $p_{S^2}$
satisfies:
$$p_{S^2}=\left\{\begin{array}{ll} -p_{\infty} \circ p_1 \;\;& \mbox{on} \;S^2_{\infty} \times \RR^3\\ p_{\infty} \circ p_2 \;\;& \mbox{on} \;  \RR^3 \times S^2_{\infty}\\ p_2 \;\;& \mbox{on} \;U\RR^3 
{=} \RR^3 \times S^2\end{array}\right.$$
where $p_1$ and $p_2$ denote the projections on the first and second factor with respect to the above expressions.
\end{lemma}
\bp 
Near the diagonal of $\RR^3$, we have a chart of $C_2(S^3)$
$$\psi_d: \RR^3 \times [0,\infty[  \times S^2 \longrightarrow C_2(S^3)$$
that maps $( x \in \RR^3, \lambda \in ]0,\infty[ ,y \in S^2)$ to $(x, x + \lambda y) \in (\RR^3)^2$. Here, $p_{S^2}$ extends as the projection onto the $S^2$ factor.\\
Consider the orientation-reversing embedding $\phi_{\infty}$
$$\begin{array}{llll}\phi_{\infty}: &\RR^3 &\longrightarrow &S^3\\
& \mu (x \in S^2) & \mapsto & \left\{\begin{array}{ll} \infty \;&\;\mbox{if}\; \mu=0\\
\frac{1}{\mu}x \;&\;\mbox{otherwise.} \end{array}\right.\end{array}$$
Note that this chart induces the already mentioned identification $p_{\infty}$ of the ill-oriented unit normal bundle $S^2_\infty$ to $\{\infty\}$ in $S^3$ with $S^2$. When $\mu \neq 0$,
$$p_{S^2}(\phi_{\infty}(\mu x), y \in \RR^3)= \frac{\mu y-x}{\norm{\mu y- x}}.$$
Then $p_{S^2}$ can be smoothly extended on $S^2_\infty \times \RR^3$ (where $\mu=0$) by
$$p_{S^2}(x \in S^2_\infty, y \in \RR^3) = -x.$$
Similarly, set
$p_{S^2}( x \in \RR^3, y \in S^2_\infty) = y.$
Now, with the map $\psi$ of the proof of Lemma~\ref{lemstrucctwo}, when $x$ and $y$ are  not equal to zero and when they are distinct,
$$p_{S^2}\circ \psi((\lambda ,(x,y)))=\frac{\frac{y}{\norm{y}^2}-\frac{x}{\norm{x}^2}}
{\norm{\frac{y}{\norm{y}^2}-\frac{x}{\norm{x}^2}}}
=\frac{\norm{x}^2y-\norm{y}^2x}
{\norm{\norm{x}^2y-\norm{y}^2x}}$$
when $\lambda \neq 0$. This map naturally extends to $\blowup{M^2}{(\infty,\infty)}$ outside the boundaries of $\infty \times \blowup{M}{\infty}$, 
$\blowup{M}{\infty} \times \infty$ and $\mbox{diag}(\blowup{M}{\infty})$ by keeping the same formula when $\lambda = 0$.

Let us check that $p_{S^2}$ smoothly extends over the boundary of the diagonal of $\blowup{M}{\infty}$. There is a chart of $C_2(M)$ near the preimage of this boundary
in  $C_2(M)$
$$\psi_2: [0,\infty[ \times [0,\infty[ \times S^2 \times S^2 \longrightarrow C_2(S^3)$$
that maps $(\lambda \in ]0,\infty[ , \mu \in  ]0,\infty[,  x \in S^2, y \in S^2)$
to $(\phi_{\infty}(\lambda x), \phi_{\infty}(\lambda(x + \mu y)))$ where $p_{S^2}$ reads
$$(\lambda,\mu,x,y) \mapsto \frac{y- 2\langle x,y \rangle x -\mu x}
{\norm{y- 2\langle x,y \rangle x -\mu x}},$$ and therefore smoothly extends when $\mu=0$. We similarly check that $p_{S^2}$ smoothly extends over the boundaries of $(\infty \times \blowup{M}{\infty})$ and 
$(\blowup{M}{\infty} \times \infty)$.
\eop

Let $\taust$ \index{N}{taus@$\taust$} denote the standard parallelization of $\RR^3$. Say that a parallelization $$\tau \colon \check{M} \times \RR^3 \rightarrow T \check{M}$$ of $\check{M}$ that coincides with $\taust$ 
on $\check{B}_{1,\infty}$ is {\em asymptotically standard.\/} According to Subsection~\ref{subproofpar}, such a parallelization exists.
Such a parallelization identifies $U\check{M}$ with $\check{M} \times S^2$.

\begin{proposition}
\label{propprojbord}
For any asymptotically standard parallelization $\tau$ of $\check{M}$,
there exists a smooth map $p_{\tau} \colon \partial C_2(M) \rightarrow S^2$ \index{N}{ptaus@$p_{\tau}$} such that
$$p_{\tau}= 
\left\{\begin{array}{ll} 
p_{S^2} \;\;& \mbox{on} \;  p_{M^2}^{-1}(\infty,\infty)\\
-p_{\infty} \circ p_1 \;\;& \mbox{on} \;S^2_{\infty} \times \check{M}\\
 p_{\infty} \circ p_2 \;\;& \mbox{on} \;  \check{M} \times S^2_{\infty}\\ 
p_2 \;\;& \mbox{on} \;U\check{M} \stackrel{\tau} {=} \check{M} \times S^2\end{array}\right.$$
where $p_1$ and $p_2$ denote the projections on the first and second factor with respect to the above expressions.
\end{proposition}
\bp This is a consequence of Lemma~\ref{lemextproj}.\eop

Since $C_2(M)$ is homotopy equivalent to $(\check{M}^2 \setminus \mbox{diag})$, according to Lemma~\ref{lemhomstwo}, $H_2(C_2(M);\QQ)=\QQ[S]$ where the canonical generator $[S]$
is the homology class of a fiber of $U\check{M} \subset \partial C_2(M)$.
For a $2$-component link $(J,K)$ of $\check{M}$, the homology class $[J \times K]$ of
$J \times K$ in $H_2(C_2(M);\QQ)$ reads $lk(J,K)[S]$, according to Subsection~\ref{subgenGaussdef} and to Proposition~\ref{propdeflkeq}.

Define an {\em asymptotic rational homology $\RR^3$\/} as a pair $(\check{M},\tau)$ where  $\check{M}$ is $3$-manifold that reads as the union over $]1,2] \times S^2$ of a rational homology ball $B_M$ and the complement $\check{B}_{1,\infty}$ of the unit ball of $\RR^3$, and $\tau$ is an asymptotically standard parallelization of $\check{M}$.
Since such a pair $(\check{M},\tau)$ canonically defines the rational homology $3$--sphere $M=\check{M} \cup \{\infty\}$, ``Let $(\check{M},\tau)$ be an asymptotic rational homology $\RR^3$'' is a shortcut for
``Let $M$ be a rational homology $3$--sphere equipped with an asymptotically standard parallelization $\tau$ of $\check{M}$''.

\subsection{On propagators}
\label{subprop}

\begin{definition}
Let $(\check{M},\tau)$ be an asymptotic rational homology $\RR^3$. A {\em propagating chain\/} of $(C_2(M),\tau)$ is a $4$--chain $\prop$ of $C_2(M)$
such that $\partial \prop =p_{\tau}^{-1}(a)$ for some $a \in S^2$.
A {\em propagating form\/} of $(C_2(M),\tau)$ is a closed $2$-form $\omega_p$ on $C_2(M)$ whose restriction to $\partial C_2(M)$ reads $p_{\tau}^{\ast}(\omega)$ for some $2$-form $\omega$ of $S^2$ such that $\int_{S^2}\omega=1$.
Propagating chains and propagating forms are simply called {\em propagators\/} when their nature is clear from the context.
\end{definition}

\begin{example}
Recall the map $p_{S^2} \colon C_2(S^3) \rightarrow S^2$ of Lemma~\ref{lemextproj}. For any $a \in S^2$, $p_{S^2}^{-1}(a)$ is a propagating chain of $(C_2(S^3),\taust)$, and for any $2$-form $\omega$ of $S^2$ such that $\int_{S^2}\omega=1$, $p_{S^2}^{\ast}(\omega)$ is a propagating form of $(C_2(S^3),\taust)$.
\end{example}

Propagating chains exist because the $3$-cycle $p_{\tau}^{-1}(a)$ of $\partial C_2(M)$ bounds in $C_2(M)$ since $H_3(C_2(M);\QQ)=0$.
Dually, propagating forms exist because the restriction induces a surjective map $H^2(C_2(M);\RR) \rightarrow H^2(\partial C_2(M);\RR)$ since $H^3(C_2(M),\partial C_2(M);\RR)=0$.
Explicit constructions of propagating chains associated to Morse functions or Heegaard diagrams can be found in \cite{lesHC}.

\begin{lemma}
\label{lemlkprop}
Let $(\check{M},\tau)$ be an asymptotic rational homology $\RR^3$. Let $C$ be a two-cycle of $C_2(M)$. For any propagating chain $\prop$ of $(C_2(M),\tau)$ transverse to $C$ and for any propagating form $\omega_p$ of $(C_2(M),\tau)$,
$[C]=\int_{C}\omega_p[S]=\langle C, \prop \rangle_{C_2(M)}[S]$ in $H_2(C_2(M);\QQ)=\QQ[S]$.
In particular, for any two-component link $(J,K)$ of $\check{M}$.
$$lk(J,K)=\int_{J \times K}\omega_p=\langle J \times K, \prop \rangle_{C_2(M)}.$$
\end{lemma}
\bp
Fix a propagating chain $\prop$, the algebraic intersection $\langle C, \prop \rangle_{C_2(M)}$ only depends on the homology class $[C]$ of $C$ in $C_2(M)$. Similarly, since $\omega_p$ is closed, $\int_{C}\omega_p$ only depends on $[C]$.
(Indeed, if $C$ and $C^{\prime}$ cobound a chain $D$, $C \cap \prop$ and $C^{\prime} \cap \prop$ cobound $\pm (D\cap \prop)$, and $\int_{\partial D=C^{\prime}-C}\omega_p=\int_Dd\omega_p$ according to the Stokes theorem.)
Furthermore, the dependance on $[C]$ is linear. Therefore it suffices to check the lemma for a cycle that represents the canonical generator $[S]$ of $H_2(C_2(M);\QQ)$. Any fiber of $U\check{M}$ is such a cycle.
\eop

\subsection{The \texorpdfstring{$\Theta$}{Theta}-invariant of \texorpdfstring{$(M,\tau)$}{(M,tau)}}
\label{subThetapar}

Note that the intersection of transverse (oriented) submanifolds is an associative operation, so that $A\cap B \cap C$ is well defined. Furthermore, for a connected manifold $N$, the class of a $0$-cycle in $H_0(M;\QQ)=\QQ[m]=\QQ$ is a well-defined number, so that the {\em algebraic intersection\/} of several transverse submanifolds whose codimension sum is the dimension of the ambient manifold is well defined as the homology class of their (oriented) intersection. This extends to rational chains, multilinearly.
Thus, for three such transverse submanifolds $A$, $B$, $C$ in a manifold $D$, their algebraic intersection $\langle A, B, C \rangle_D$ is the sum over the intersection points $a$ of the associated signs, where the sign of $a$ is positive if and only if the orientation of $D$ is induced by the orientation of ${\normbun}_aA \oplus {\normbun}_aB \oplus {\normbun}_aC$.

\begin{theorem}
Let $(\check{M},\tau)$ be an asymptotic rational homology $\RR^3$.
Let $\prop_a$, $\prop_b$ and $\prop_c$ be three pairwise transverse propagators of $(C_2(M),\tau)$ with respective boundaries $p_{\tau}^{-1}(a)$, $p_{\tau}^{-1}(b)$ and  $p_{\tau}^{-1}(c)$ for three distinct points $a$, $b$ and $c$ of $S^2$,
then \index{N}{ThetaMtau@$\Theta(M,\tau)$}
$$\Theta(M,\tau)=\langle \prop_a, \prop_b, \prop_c \rangle_{C_2(M)}$$
does not depend on the chosen propagators $\prop_a$, $\prop_b$ and $\prop_c$. It is a topological invariant of $(M,\tau)$. 
For any three propagating chains $\omega_a$, $\omega_b$ and $\omega_c$ of $(C_2(M),\tau)$,
$$\Theta(M,\tau)=\int_{C_2(M)} \omega_a \wedge \omega_b \wedge  \omega_c.$$
\end{theorem}
\bp
Since $H_4(C_2(M))=0$, if the propagator $\prop_a$ is changed to a propagator $\prop^{\prime}_a$ with the same boundary, $(\prop^{\prime}_a-\prop_a)$ bounds a $5$-dimensional chain $W$ transverse to $\prop_b \cap \prop_c$. The $1$-dimensional chain $W \cap \prop_b \cap \prop_c$ does not meet $\partial C_2(M)$ since $\prop_b \cap \prop_c$ does not meet $\partial C_2(M)$. Therefore, up to a well-determined sign, the boundary of $W \cap \prop_b \cap \prop_c$ is $\prop^{\prime}_a \cap \prop_b \cap \prop_c -\prop_a \cap \prop_b \cap \prop_c$. 
This shows that $\langle \prop_a, \prop_b, \prop_c \rangle_{C_2(M)}$ is independent of $\prop_a$ when $a$ is fixed. Similarly, it is independent of $\prop_b$ and $\prop_c$ when $b$ and $c$ are fixed. Thus, $\langle \prop_a, \prop_b, \prop_c \rangle_{C_2(M)}$ is a rational function on the connected set of triples $(a,b,c)$ of distinct point of $S^2$. It is easy to see that this function is continuous. Thus, it is constant.

Let us similarly prove that $\int_{C_2(M)} \omega_a \wedge \omega_b \wedge  \omega_c$ is independent of the propagating forms $\omega_a$, $\omega_b$ and $\omega_c$. Assume that the form $\omega_a$, which restricts to $\partial C_2(M)$ as $p_{\tau}^{\ast}(\omega_A)$, is changed to $\omega^{\prime}_a$, which restricts to $\partial C_2(M)$ as $p_{\tau}^{\ast}(\omega^{\prime}_A)$. 

\begin{lemma}
\label{lemetactwo}
There exists a one-form $\eta_A$ on $S^2$ such that $\omega^{\prime}_A=\omega_A +d\eta_A$.
For any such $\eta_A$, there exists a one-form $\eta$ on $C_2(M)$ such that $\omega^{\prime}_a-\omega_a=d\eta$, and the restriction of $\eta$ to $\partial C_2(M)$ is $p_{\tau}^{\ast}(\eta_A)$.
\end{lemma}
\noindent{\sc Proof of the lemma:}
Since $\omega_a$ and  $\omega^{\prime}_a$ are cohomologous, there exists a one-form $\eta$ on $C_2(M)$ such that $\omega^{\prime}_a=\omega_a +d\eta$. Similarly, since $\int_{S^2}\omega^{\prime}_A=\int_{S^2}\omega_A$, there exists a one-form $\eta_A$ on $S^2$ such that $\omega^{\prime}_A=\omega_A +d\eta_A$.
On $\partial C_2(M)$, $d(\eta-p_{\tau}^{\ast}(\eta_A))=0$.
Thanks to the exact sequence
$$ 0=H^1(C_2(M)) \longrightarrow H^1(\partial C_2(M)) \longrightarrow H^2(C_2(M), \partial C_2(M)) \cong H_4(C_2(M))=0, $$
$H^1(\partial C_2(M))=0$. 
Therefore, there exists a function $f$ from $\partial C_2(M)$ to $\RR$ such that $$df =\eta -p(\tau)^{\ast}(\eta_A)$$ on $ \partial C_2(M)$.
Extend $f$ to a $C^{\infty}$ map on $C_2(M)$ and change $\eta$ into $(\eta -df)$.
\eop

Then $$\begin{array}{ll}\int_{C_2(M)} \omega^{\prime}_a \wedge \omega_b \wedge  \omega_c - \int_{C_2(M)} \omega_a \wedge \omega_b \wedge  \omega_c & =\int_{C_2(M)}d(\eta \wedge \omega_b \wedge  \omega_c)=\int_{\partial C_2(M)}\eta \wedge \omega_b \wedge  \omega_c\\
        &=\int_{\partial C_2(M)}p(\tau)^{\ast}(\eta_A \wedge \omega_B \wedge  \omega_C)=0
       \end{array}$$
since any $5$-form on $S^2$ vanishes.
Thus, $\int_{C_2(M)} \omega_a \wedge \omega_b \wedge  \omega_c$ is independent of the propagating forms $\omega_a$, $\omega_b$ and $\omega_c$. Now, we can choose the propagating forms $\omega_a$, $\omega_b$ and $\omega_c$, Poincar\'e dual to $\prop_a$,
$\prop_b$ and $\prop_c$, and supported in very small neighborhoods of $\prop_a$,
$\prop_b$ and $\prop_c$, respectively, so that the intersection of the three supports is a very small neighborhood of $\prop_a \cap \prop_b \cap \prop_c$, where it can easily be seen that $\int_{C_2(M)} \omega_a \wedge \omega_b \wedge  \omega_c=\langle \prop_a, \prop_b, \prop_c \rangle_{C_2(M)}$.
\eop

In particular, $\Theta(M,\tau)$ reads $\int_{C_2(M)} \omega^3$ for any propagating chain $\omega$ of $(C_2(M),\tau)$. Since such a propagating chain represents the linking number, $\Theta(M,\tau)$ can be thought of as the {\em cube of the linking number with respect to $\tau$.\/}

When $\tau$ varies continuously, $\Theta(M,\tau)$ varies continuously in $\QQ$ so that $\Theta(M,\tau)$ is an invariant of the homotopy class of $\tau$.

\subsection{Parallelisations of \texorpdfstring{$3$}{3}--manifolds and Pontrjagin classes}

In this subsection, $M$ denotes a smooth, compact oriented $3$-manifold with possible boundary $\partial M$.
Recall that such a $3$-manifold is parallelizable.

Let $GL^+(\RR^3)$ denote the group of orientation-preserving linear isomorphisms of $\RR^3$.
Let $C^0\left((M,\partial M),(GL^+(\RR^3),1)\right)$ denote the set of maps $$g: (M, \partial M) \longrightarrow (GL^+(\RR^3),1)$$
from $M$ to $GL^+(\RR^3)$ that send $\partial M$ to the unit $1$ of $GL^+(\RR^3)$.
Let $[(M,\partial M),(GL^+(\RR^3),1)]$ denote the group of homotopy classes of such maps, with the group structure induced by the multiplication of maps, using the multiplication in $GL^+(\RR^3)$.
For a map $g$ in $C^0\left((M,\partial M),(GL^+(\RR^3),1)\right)$, set \index{N}{psiR@$\psi_{\RR}$}
$$\begin{array}{llll} 
\psi_{\RR}(g): &M \times \RR^3 &\longrightarrow  &M \times \RR^3\\
&(x,y) & \mapsto &(x,g(x)(y)).\end{array}$$
Let $\tau_M \colon M \times \RR^3 \rightarrow TM$ be a parallelization of $M$.
Then any parallelization $\tau$ of $M$ that coincides with $\tau_M$ on $\partial M$ reads
$$\tau = \tau_M \circ \psi_{\RR}(g)$$ for some $g \in C^0\left((M,\partial M),(GL^+(\RR^3),1)\right)$.

Thus, fixing $\tau_M$ identifies the set of homotopy classes of parallelizations of $M$ fixed on $\partial M$ with 
the group
$[(M,\partial M),(GL^+(\RR^3),1)]$.
Since $GL^+(\RR^3)$ deformation retracts onto the group $SO(3)$ of orientation-preserving linear isometries of $\RR^3$, $[(M,\partial M),(GL^+(\RR^3),1)]$ is isomorphic to $[(M,\partial M),(SO(3),1)]$.

See $S^3$ as $B^3/\partial B^3$ where $B^3$ is the standard ball of radius $2\pi$ of $\RR^3$ seen as $([0,2\pi] \times S^2)/(0 \sim \{0\} \times S^2)$.
Let $\rho \colon B^3 \rightarrow SO(3)$ map $(\theta \in[0,2\pi], v \in S^2)$ to the rotation $\rho(\theta,v)$ with axis directed by $v$ and with angle $\theta$. This map induces the double covering $\tilde{\rho} \colon S^3 \rightarrow SO(3)$, which orients $SO(3)$ and which allows one to deduce the first three homotopy groups of $SO(3)$ from the ones of $S^3$. They are $\pi_1(SO(3))=\ZZ/2\ZZ$, $\pi_2(SO(3))=0$ and $\pi_3(SO(3))=\ZZ[\tilde{\rho}]$. For $v\in S^2$, $\pi_1(SO(3))$ is generated by the class of the loop that maps $\exp(i \theta) \in S^1$ to the rotation $\rho(\theta,v)$.

Note that a map $g$ from $(M,\partial M)$ to $(SO(3),1)$ has a degree $\mbox{deg}(g)$, which may be defined as the differential degree at a regular value (different from $1$) of $g$. It can also be defined homologically, by $H_3(g)[M,\partial M] =\mbox{deg}(g)[SO(3),1]$.

The following theorem is proved in Section~\ref{secmorepar}.

\begin{theorem}
\label{thmpone}
For any smooth compact connected oriented $3$-manifold $M$, the group $$[(M,\partial M),(SO(3),1)]$$ is abelian, and
the degree $$\mbox{deg} \colon [(M,\partial M),(SO(3),1)] \longrightarrow \ZZ$$ is a group homomorphism,
which induces an isomorphism 
$$\mbox{deg} \colon [(M,\partial M),(SO(3),1)] \otimes_{\ZZ}\QQ \longrightarrow \QQ.$$
When $\partial M=\emptyset$, (resp. when $\partial M=S^2$), there exists a canonical map $p_1$ from the set of homotopy classes of parallelizations of $M$ (resp. that coincide with $\taust$ near $S^2$) such that for any map $g$ in $C^0\left((M,\partial M),(SO(3),1)\right)$, for any trivialization $\tau$ of $TM$
$$p_1(\tau \circ \psi_{\RR}(g))-p_1(\tau)=2\mbox{deg}(g).$$
\end{theorem}

The definition of the map $p_1$ is given in Subsection~\ref{subdefpont}, it involves relative Pontrjagin classes. When $\partial M=\emptyset$, the map $p_1$ coincides with the map $h$ that is studied by Kirby and Melvin in \cite{km} under the name {\em Hirzebruch defect.\/}  See also \cite[\S 3.1]{hirzebruchEM}.

Since $[(M,\partial M),(SO(3),1)]$ is abelian,
the set of parallelizations of $M$ that are fixed on $\partial M$ is an affine space with translation group $[(M,\partial M),(SO(3),1)]$.

Recall that $\rho \colon B^3 \rightarrow SO(3)$ maps $(\theta \in[0,2\pi], v \in S^2)$ to the rotation with axis directed by $v$ and with angle $\theta$. 
Let $M$ be an oriented connected $3$-manifold with possible boundary. For a ball $B^3$ embedded in $M$, let $\rho_M(B^3) \in C^0\left((M,\partial M),(SO(3),1)\right)$ be a (continuous) map that coincides with $\rho$ on $B^3$ and that maps the complement of $B^3$ to the unit of $SO(3)$. The homotopy class of $\rho_M(B^3)$ is well-defined.
\begin{lemma}
\label{lemdegrho}
 $\mbox{deg}(\rho_M(B^3))=2$
\end{lemma}
\bp Exercise.
\eop

\subsection{Defining a \texorpdfstring{$\QQ$}{Q}-sphere invariant from \texorpdfstring{$\Theta$}{Theta}}
\label{subThetarat}

Recall that an asymptotic rational homology $\RR^3$ is a pair $(\check{M},\tau)$ where  $\check{M}$ is $3$-manifold that reads as the union over $]1,2] \times S^2$ of a rational homology ball $B_M$ and the complement $\check{B}_{1,\infty}$ of the unit ball of $\RR^3$, and that is equipped with an asymptotically standard parallelization $\tau$.

In this subsection, we prove the following proposition.

\begin{proposition}
\label{propThetap}
Let $(\check{M},\tau)$ be an asymptotic rational homology $\RR^3$.
For any map $g$ in $C^0\left((B_M,B_M \cap \check{B}_{1,\infty}),(SO(3),1)\right)$ trivially extended to $\check{M}$,
$$\Theta(M, \tau \circ \psi_{\RR}(g))-\Theta(M, \tau)=\frac12\mbox{deg}(g).$$
\end{proposition}

Theorem~\ref{thmpone} allows us to derive the following corollary from Proposition~\ref{propThetap}.
\begin{corollary}
\label{corThetap} \index{N}{ThetaM@$\Theta(M)$}
 $\Theta(M)=\Theta(M,\tau) - \frac1{4}{p_1(\tau)}$ is an invariant of $\QQ$-spheres.
\end{corollary}
\eop

\begin{lemma}
\label{lemtheprimhom}
 $\Theta(M, \tau \circ \psi_{\RR}(g))-\Theta(M, \tau)$ is independent of $\tau$. Set $\Theta^{\prime}(g)=\Theta(M, \tau \circ \psi_{\RR}(g))-\Theta(M, \tau)$. Then $\Theta^{\prime}$ is a homomorphism from $[(B_M,B_M \cap \check{B}_{1,\infty}),(SO(3),1)]$ to $\QQ$.
\end{lemma}
\bp
For $d=a$, $b$ or $c$, the propagator $\prop_d$ of $(C_2(M),\tau)$ can be assumed to be a product $[-1,0] \times p_{\tau |UB_M}^{-1}(d)$ on a collar $[-1,0] \times UB_M$ of $UB_M$ in $C_2(M)$. Since $H_3([-1,0] \times UB_M;\QQ)=0$, $ \left(\partial ([-1,0] \times p_{\tau |UB_M}^{-1}(d)) \setminus (0 \times p_{\tau |UB_M}^{-1}(d)) \right) \cup (0 \times p_{\tau \circ \psi_{\RR}(g) |UB_M}^{-1}(d))$ bounds a chain $G_d$. 

The chains $G_a$, $G_b$ and $G_c$ can be assumed to be transverse. Construct the propagator $\prop_d(g)$ of $(C_2(M),\tau \circ \psi_{\RR}(g))$ from $\prop_d$ by replacing $[-1,0] \times p_{\tau |UB_M}^{-1}(d)$ by $G_d$ on $[-1,0] \times UB_M$. Then
$$\Theta(M, \tau \circ \psi_{\RR}(g))-\Theta(M, \tau)=\langle G_a,G_b,G_c\rangle_{[-1,0] \times UB_M}.$$
Using $\tau$ to identify $UB_M$ with $B_M \times S^2$ allows us to see that $\Theta(M, \tau \circ \psi_{\RR}(g))-\Theta(M, \tau)$ is independent of $\tau$. Then it is easy to observe that $\Theta^{\prime}$ is a homomorphism from $[(B_M,\partial B_M),(SO(3),1)]$ to $\QQ$.
\eop

According to Theorem~\ref{thmpone} and to Lemma~\ref{lemdegrho}, it suffices to prove that $\Theta^{\prime}(\rho_M(B^3))=1$ in order to prove Proposition~\ref{propThetap}. It is easy to see that $\Theta^{\prime}(\rho_M(B^3))=\Theta^{\prime}(\rho)$. Thus, we are left with the proof of the following lemma.

\begin{lemma}
\label{lemtheprimrho}
 $\Theta^{\prime}(\rho)=1$.
\end{lemma}

Again, see $B^3$ as $([0,2\pi] \times S^2)/(0 \sim \{0\} \times S^2)$.
We first prove the following lemma:

\begin{lemma}
\label{lemL}
Let $a$ be the North Pole. The point $(-a)$ is regular for the map
$$\begin{array}{llll}\rho_a \colon &B^3 &\rightarrow &S^2\\
   & m &\mapsto & \rho(m)(a)
  \end{array}
$$
and its preimage (cooriented by $S^2$ via $\rho_a$) is the knot $L_a=\{\pi\} \times E$, where $E$ is the equator that bounds the Southern Hemisphere.
\end{lemma}
\bp It is easy to see that $\rho_a^{-1}(-a)=\pm \{\pi\} \times E$.
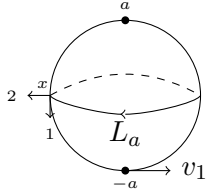
\begin{center}
\begin{tikzpicture} \useasboundingbox (-1,-1.2) rectangle (1,1.2);
\draw (0,0) circle (1) (-1.1,-.05) node[above]{\tiny $x$} (0,1) node[above]{\tiny $a$} (0,-.95) node[below]{\tiny $-a$} (0,-.2) node[below]{\small $L_a$};
\draw [-<] (-1,0) .. controls (-.95,-.1) and (-.3,-.25) .. (0,-.25);
\draw (0,-.25) .. controls (.3,-.25) and  (.95,-.1) .. (1,0);
\draw [out=29,in=151,dashed] (-1,0) to (1,0); 
\draw [->] (-1,0) -- (-1,-.3) node[below]{\tiny $1$};
\draw [->] (0,-1) -- (.6,-1) node[right]{\small $v_1$};
\draw [->] (-1,0) -- (-1.3,0) node[left]{\tiny $2$};
\fill (0,1) circle (1.5pt) (0,-1) circle (1.5pt);
\end{tikzpicture}
\end{center}
Let $x \in \{\pi\} \times E$.
When $m$ moves along the great circle that contains $a$ and $x$ from $x$ towards $(-a)$ in $\{\pi\} \times S^2$, $\rho(m)(a)$ moves from $(-a)$ in the same direction, which will be the direction of the tangent vector $v_1$ of $S^2$ at $(-a)$, counterclockwise in our picture, where $x$ is on the left. Then in our picture, $S^2$ is oriented at $(-a)$ by $v_1$ and by the tangent vector $v_2$ at $(-a)$ towards us. In order to move $\rho(\theta,v)(a)$ in the $v_2$ direction, one increases $\theta$
so that $L_a$ is cooriented and oriented like in the figure.
\eop

\noindent{\sc Proof of Lemma~\ref{lemtheprimrho}:}
We use the notation of the proof of Lemma~\ref{lemtheprimhom} and we construct an explicit $G_a$ in $[-1,0] \times UB^3 \stackrel{\taust} {=}[-1,0] \times B^3 \times S^2$.

When $\rho(m)(a) \neq - a$, there is a unique geodesic arc $[a,\rho(m)(a)]$ with length $(\ell \in [0, \pi[)$ from $a$ to $\rho(m)(a)=\rho_a(m)$.
For $t \in [0,1]$, let $X_t(m) \in [a,\rho_a(m)]$ be such that
the length of $[X_0(m)=a,X_t(m)]$ is $t\ell$.
This defines $X_t$ on $(M\setminus L_a)$, $X_1(m)=\rho_a(m)$. Let us show how the definition of $X_t$ smoothly extends on the manifold $\blowup{B^3}{L_a}$ obtained from $B^3$ by blowing up $L_a$.

The map $\rho_a$ maps the normal bundle to $L_a$ to a disk of $S^2$ around $(-a)$,
by an orientation-preserving diffeomorphism on every fiber (near the origin). In particular, $\rho_a$ induces a map $\tilde{\rho}_a$ from the unit normal bundle to $L_a$ to the unit normal bundle to $(-a)$ in $S^2$, which preserves the orientation of the fibers.
Then for an element $y$ of the unit normal bundle to $L_a$ in $M$, define $X_t(y)$ as before on the half great circle $[a,-a]_{\tilde{\rho}_a(-y)}$ from $a$ to $(-a)$ that is tangent to $\tilde{\rho}_a(-y)$ at $(-a)$ (so that $\tilde{\rho}_a(-y)$ is an outward normal to $[a,-a]_{\tilde{\rho}_a(-y)}$ at $(-a)$).
This extends the definition of $X_t$ , continuously.

The whole sphere is covered with degree $(-1)$ by the image of $([0,1] \times U{\normbun}_x(L_a))$,
where the fiber $U{\normbun}_x(L_a)$ of the unit normal bundle to $L_a$ is oriented as the boundary of a disk in the fiber of the normal bundle.
Let $G_h(a)$ be the closure of 
$\left(\cup_{t \in [0,1], m\in (B^3 \setminus L_a)} \left(m,{X_t}(m)\right)\right)$ in $UB^3$.
$$G_h(a)=\cup_{t \in [0,1], m\in \blowup{B^3}{L_a}}\left(p_{B^3}(m),{X_t}(m)\right).$$
Then
$$\partial G_h = -(B^3\times a) +\cup_{m\in B^3}(m,\rho_a(m)) + \cup_{t \in [0,1]}X_{t}(-\partial \blowup{S^3}{L_a})$$
where $(-\partial \blowup{S^3}{L_a})$ is oriented like  $\partial {\neightub}(L_a)$ so that the last summand reads $(-L_a \times S^2)$
because the sphere is covered with degree $(-1)$ by the image of $([0,1] \times U{\normbun}_x(L_a))$.

Let $D_a$ be a disk bounded by $L_a$ in $B^3$. 
Set $G(a)=G_h(a) +D_a \times S^2 $ so that $\partial G(a)=-(B^3\times a) +\cup_{m\in B^3}(m,\rho_a(m))$.
Now let $\iota$ be the endomorphism of $UB^3$ over $B^3$ that maps a unit vector to the opposite one.
Set $$\begin{array}{lllll}&G_a&=[-1,-2/3]\times B^3\times a &+ \{-2/3\}\times G(a) &+ [-2/3,0]\times \cup_{m\in B^3}(m,\rho_a(m))\\
 \mbox{and}&      G_{-a}&=[-1,-1/3]\times B^3 \times (-a) &+ \{-1/3\}\times \iota(G(a)) &+ [-1/3,0]\times \cup_{m\in B^3}(m,\rho(m)(-a)).
      \end{array}
$$

Then $$G_a \cap G_{-a}=[-2/3,-1/3] \times L_a \times (-a) + \{-2/3\} \times D_a \times (-a) -\{-1/3\} \times \cup_{m\in D_a}(m,\rho_a(m)).$$
Finally, $\Theta^{\prime}(\rho)$ is the algebraic intersection of $G_a \cap G_{-a}$ with $\prop_c(\rho)$ in $C_2(M)$. This intersection coincides with the algebraic intersection of $G_a \cap G_{-a}$ with any propagator of $(C_2(M),\tau)$ according to Lemma~\ref{lemlkprop}.
Therefore $$ \Theta^{\prime}(\rho)=\langle \prop_a,G_a \cap G_{-a} \rangle_{[-1,0]\times S^2 \times B^3}=-\mbox{deg}_a(\rho_a \colon D_a \rightarrow S^2).$$
The orientation of $L_a$ allows us to choose $(-D_a)$ as the Northern Hemisphere, the image of this hemisphere under $\rho_a$ covers the sphere with degree $1$ so that $\Theta^{\prime}(\rho)=1$.
\eop

\newpage
\section{An introduction to finite type invariants}
\setcounter{equation}{0}

This section contains the needed background from the theory of finite type invariants. 
It allows us to introduce the target space generated by Feynman-Jacobi diagrams, for the general invariants presented in Section~\ref{secconstcsi}, in a progressive way.

Theories of finite type invariants are useful to characterize invariants. Such a theory allowed Greg Kuperberg and Dylan Thurston to identify $\Theta/6$ with the Casson invariant $\lambda$ for integer homology $3$-spheres, in \cite{kt}. The invariant $\lambda$ was defined by Casson in 1984 as an algebraic count of conjugacy classes of irreducible representations from $\pi_1(M)$ to $SU(2)$. See \cite{akmc,gm,mar}.
The Kuperberg--Thurston result above was generalized to the case of rational homology $3$-spheres in \cite[Theorem~2.6 and Corollary~6.14]{lessumgen}. Thus, for any rational homology $3$--sphere $M$,
$$\Theta(M)=6\lambda(M),$$
where $\lambda$ is the Walker generalization of the Casson invariant to rational homology $3$--spheres, which is normalized like in \cite{akmc,gm,mar} for integer homology $3$--spheres, and like $\frac{1}{2}\lambda_W$ for rational homology $3$--spheres with respect to the Walker normalisation $\lambda_W$ of \cite{wal}.

For invariants of knots and links in $\RR^3$, the base of the theory of finite type invariants was mainly established by Bar-Natan in \cite{barnatan}. A more complete review of this theory has been written by Chmutov, Duzhin and Mostovoy in \cite{chmutovmostov}. For integer homology $3$--spheres, the theory was started by Ohtsuki in \cite{ohtkno} and further developed by Goussarov, Habiro, Le and others. See \cite{ggp,habiro,le}. Delphine Moussard developed a theory of finite type invariants for rational homology $3$--spheres in \cite{moussardAGT}. Her suitable theory is based on the Lagrangian-preserving surgeries defined below.

\subsection{Lagrangian-preserving surgeries}

\begin{definition}
\label{defhhh}
 An \emph{integer (resp. rational) homology handlebody} of genus $g$ is a compact oriented $3$-manifold $A$ that has the same integral (resp. rational) homology as the usual solid handlebody $\handlebody_g$ below.

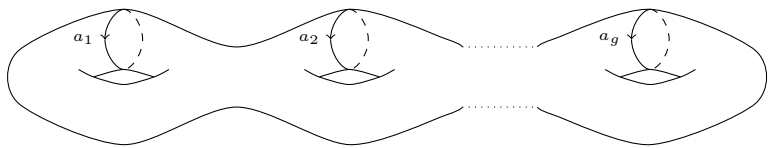
\begin{center}
 \begin{tikzpicture} \useasboundingbox (0,-1) rectangle (11,1); 
\draw [->] (2,.9) .. controls (1.9,.9) and (1.75,.7) .. (1.75,.5) node[left]{\tiny $a_1$};  
\draw (1.75,.5) .. controls (1.75,.3) and (1.9,.1) .. (2,.1);
\draw [dashed] (2,.9) .. controls (2.1,.9) and (2.25,.7) .. (2.25,.5) .. controls (2.25,.3) and (2.1,.1) .. (2,.1);
\draw plot[smooth] coordinates{(1.4,.1) (1.6,0) (2,-.1) (2.4,0) (2.6,.1)};
\draw plot[smooth] coordinates{(1.6,0) (2,.1) (2.4,0)};
\draw plot[smooth] coordinates{(6.5,.4) (6.3,.5) (5,.9) (3.5,.4) (2,.9) (.6,.35) (.6,-.35) (2,-.9) (3.5,-.4) (5,-.9) (6.3,-.5) (6.5,-.4)};
\draw [dotted] (6.5,-.4) -- (7.5,-.4)  (6.5,.4) -- (7.5,.4);
\draw plot[smooth] coordinates{(7.5,.4) (7.7,.5) (9,.9) (10.4,.35) (10.4,-.35) (9,-.9) (7.7,-.5) (7.5,-.4)};
\begin{scope}[xshift=3cm]
 \draw [->] (2,.9) .. controls (1.9,.9) and (1.75,.7) .. (1.75,.5) node[left]{\tiny $a_2$};  
\draw (1.75,.5) .. controls (1.75,.3) and (1.9,.1) .. (2,.1);
\draw [dashed] (2,.9) .. controls (2.1,.9) and (2.25,.7) .. (2.25,.5) .. controls (2.25,.3) and (2.1,.1) .. (2,.1);
\draw plot[smooth] coordinates{(1.4,.1) (1.6,0) (2,-.1) (2.4,0) (2.6,.1)};
\draw plot[smooth] coordinates{(1.6,0) (2,.1) (2.4,0)};
\end{scope}
\begin{scope}[xshift=7cm]
 \draw [->] (2,.9) .. controls (1.9,.9) and (1.75,.7) .. (1.75,.5) node[left]{\tiny $a_g$};  
\draw (1.75,.5) .. controls (1.75,.3) and (1.9,.1) .. (2,.1);
\draw [dashed] (2,.9) .. controls (2.1,.9) and (2.25,.7) .. (2.25,.5) .. controls (2.25,.3) and (2.1,.1) .. (2,.1);
\draw plot[smooth] coordinates{(1.4,.1) (1.6,0) (2,-.1) (2.4,0) (2.6,.1)};
\draw plot[smooth] coordinates{(1.6,0) (2,.1) (2.4,0)};
\end{scope}
\end{tikzpicture}
\end{center}

\begin{exo}
 Show that if $A$ is a rational homology handlebody of genus $g$, then $\partial A$ is a genus $g$ surface.
\end{exo}

The {\em Lagrangian\/} $\CL_A$ of a compact $3$-manifold $A$ is the kernel of the map induced by the inclusion from $H_1(\partial A;\QQ)$ to $H_1(\partial A;\QQ)$.
\end{definition}
In the figure, the Lagrangian of $\handlebody_g$ is freely generated by the classes of the curves $a_i$.

\begin{definition}
\label{deflagsur}
An {\em integral (resp. rational) Lagrangian-Preserving (or LP) surgery\/} $(A^{\prime}/A)$ is the replacement of an integral (resp. rational) homology handlebody $A$ embedded in the interior of a $3$-manifold $M$ by another such $A^{\prime}$ whose boundary is identified with $\partial A$
by an orientation-preserving diffeomorphism that sends $\CL_A$ to $\CL_{A^{\prime}}$.
The manifold $M(A^{\prime}/A)$ obtained by such an LP-surgery reads
$$M(A^{\prime}/A) = (M \setminus \mbox{Int}(A)) \cup_{\partial A} A^{\prime}.$$
(This only defines the topological structure of $M(A^{\prime}/A)$, but we equip $M(A^{\prime}/A)$ with its unique smooth structure.)
\end{definition}

\begin{lemma}
If $(A^{\prime}/A)$ is an integral (resp. rational) LP-surgery, then the homology of $M(A^{\prime}/A)$ with $\ZZ$-coefficients (resp. with $\QQ$-coefficients) is canonically isomorphic to
$H_{\ast}(M;\ZZ)$ (resp. to $H_{\ast}(M;\QQ)$). If $M$ is a $\QQ$-sphere, if $(A^{\prime}/A)$ is a rational LP-surgery, and if $(J,K)$ is a two-component link of $M \setminus A$, then the linking number of $J$ and $K$ in $M$ and the linking number of $J$ and $K$ in $M(A^{\prime}/A)$ coincide.
\end{lemma}
\bp Exercise. \eop

\subsection{Definition of finite type invariants}

Let $\KK=\QQ$ or $\RR$.

A $\KK$--valued {\em invariant\/} of oriented $3$-manifolds is a function from the set of $3$-manifolds, considered up to orientation-preserving diffeomorphisms to $\KK$. Let $\coprod_{i=1}^n S_i^1$ denote a disjoint union of $n$ circles, where each $S^1_i$ is a copy of $S^1$.
Here, an {\em $n$--component link\/} in a $3$-manifold $M$ is an equivalence class of smooth embeddings
$L \colon \coprod_{i=1}^n S_i^1 \hookrightarrow M$ under the equivalence relation that identifies two embeddings $L$ and $L^{\prime}$ if and only if there is an orientation-preserving diffeomorphism $h$ of $M$ such that $h(L)=L^{\prime}$.
A {\em knot\/} is a one-component link. 
A {\em link invariant\/} (resp. a {\em knot invariant\/}) is a function of links (resp. knots).
For example, $\Theta$ is an invariant of $\QQ$-spheres and the linking number is a rational invariant of two-component links in rational homology $3$--spheres

In order to study a function, it is usual to study its derivative, and the derivative of its derivative...
The derivative of a function is defined from its variations. 
For a function $f$ from $\ZZ^d = \oplus_{i=1}^d\ZZ e_i$ to $\KK$, one can define its first order derivatives $\frac {\partial f}{\partial e_i}\colon \ZZ^d \rightarrow \KK$ 
by $$\frac {\partial f}{\partial e_i}(z)=f(z+e_i)-f(z)$$
and check that all the first order derivatives of $f$ vanish if and only if $f$ is constant.
Inductively define an $n$-order derivative as a first order derivative of an $(n-1)$-order derivative for a positive integer $n$.
Then it can be checked that all the $(n+1)$-order derivatives of a function vanish if and only if $f$ is a polynomial of degree not greater than $n$.
In order to study topological invariants, we can similarly study their variations under {\em simple operations.\/}

Below, $X$ denotes one of the following sets
\begin{itemize}
\item $\ZZ^d$,
\item  the set $\CK$ of knots in $\RR^3$, the set $\CK_n$ of $n$-component links in $\RR^3$,
\item the set $\CM$ of $\ZZ$-spheres, the set $\CM_Q$ of $\QQ$-spheres.
\end{itemize}

and $\CO(X)$ denotes a set of {\em simple operations\/} acting on some elements of $X$.

For $X=\ZZ^d$, $\CO(X)$ will be made of the operations $(z \rightarrow z \pm e_i)$

For knots or links in $\RR^3$, the {\em simple operations\/} will be {\em crossing changes.\/}
A  {\em crossing change ball\/} of a link $L$ is a ball $B$ of the ambient space, where $L\cap B$ is a disjoint union of two arcs $\alpha_1$ and $\alpha_2$ properly embedded in $B$, and there exist two disjoint topological disks $D_1$ and $D_2$ embedded in $B$, such that, for $i\in \{1,2\}$, $\alpha_i \subset \partial D_i$ and $(\partial D_i \setminus \alpha_i) \subset \partial B$. After an isotopy, the projection of $(B,\alpha_1,\alpha_2)$ reads $\nccirc$ or $\pccirc$ (the corresponding pairs (ball,arcs) are isomorphic, but they are regarded in different ways), a {\em crossing change\/} is a change that does not change $L$ outside $B$ and that modifies $L$ inside $B$ by a local move $(\nccirc \rightarrow \pccirc)$ or $(\pccirc \rightarrow \nccirc)$. For the move $(\nccirc \rightarrow \pccirc)$, the crossing change is {\em positive\/}, it is {\em negative\/} for the move $(\pccirc \rightarrow \nccirc)$.

For integer (resp. rational) homology $3$--spheres, the simple operations will be integral (resp. rational) {\em $LP$-surgeries of genus $3$.\/}

Say that crossing changes are {\em disjoint\/} if they sit inside disjoint $3$-balls. Say that $LP$-surgeries $(A^{\prime}/A)$ and $(B^{\prime}/B)$ in a manifold $M$ are {\em disjoint\/} if $A$ and $B$ are disjoint in $M$.
Two operations on $\ZZ^d$ are always {\em disjoint\/} (even if they look identical).
In particular, disjoint operations commute, (their result does not depend on which one is performed first).
Let $ \underline{n}=\{1,2,\dots,n\}$. Consider the vector space $\CF_0(X)$ freely generated by $X$ over $\KK$.
For an element $x$ of $X$ and $n$ pairwise disjoint operations $o_1,\dots,o_n$ acting on $x$,
define $$[x;o_1,\dots,o_n]=\sum_{I \subseteq \underline{n}} (-1)^{\sharp I}x((o_i)_{i \in I}) \in \CF_0(X)$$
where $x((o_i)_{i \in I})$ denotes the element of $X$ obtained by performing the operations $o_i$ for $i \in I$ on $x$.
Then define $\CF_n(X)$ as the $\KK$-subspace of $\CF_0(X)$ generated by the $[x;o_1,\dots,o_n]$, for all $x \in X$ equipped with $n$ pairwise disjoint simple operations.
Since $$[x;o_1,\dots,o_n,o_{n+1}]= [x;o_1,\dots,o_n] - [x(o_{n+1});o_1,\dots,o_n],$$
$\CF_{n+1}(X) \subseteq \CF_{n}(X)$, for all $n \in \NN$.

\begin{definition}
A $\KK$--valued function $f$ on $X$, uniquely extends as a $\KK$--linear map of 
$$\CF_0(X)^{\ast}=\mbox{Hom}(\CF_0(X);\KK),$$ which is still denoted by $f$.
For an integer $n \in \NN$, the invariant (or function) $f$ is of {\em degree $\leq n$\/} if and only if $f(\CF_{n+1}(X))=0$. The {\em degree\/} of such an invariant is the smallest integer $n \in \NN$ such that $f(\CF_{n+1}(X))=0$. 
An invariant is of {\em finite type\/} if it is of degree $n$ for some $n \in \NN$. This definition depends on the chosen set of operations $\CO(X)$. We fixed our choices for our sets $X$, but other choices could lead to different notions. See \cite{ggp}.
\end{definition}

Let $\CI_n(X)=(\CF_0(X)/\CF_{n+1}(X))^{\ast}$ be the space of invariants of degree at most $n$.
Of course, for all $n \in \NN$, $\CI_n(X) \subseteq \CI_{n+1}(X)$.

\begin{example}
$\CI_n(\ZZ^d)$ is the space of polynomials of degree at most $n$ on $\ZZ^d$. (Exercise).
\end{example}

\begin{lemma}
 Any $n$-component link in $\RR^3$ can be transformed to the trivial $n$-component link below by a finite number of disjoint crossing changes.

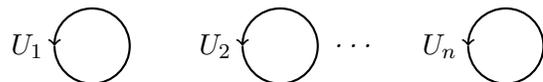
\begin{center}
 \begin{tikzpicture}
\useasboundingbox (-.5,-.6) rectangle (5.5,.6);
\draw [thick,->] (-.5,0) arc (-180:180:.5);
\draw  (-.5,0) node[left]{\small $U_1$};
\draw [thick,->] (2,0) arc (-180:180:.5);
\draw  (2,0) node[left]{\small $U_2$};
\draw  (3.5,0) node{\dots};
\draw [thick,->] (5,0) arc (-180:180:.5);
\draw  (5,0) node[left]{\small $U_n$};
\end{tikzpicture}\end{center}
\end{lemma}
\bp
Let $L$ be an $n$-component link in $\RR^3$.
Since $\RR^3$ is simply connected, there is a homotopy that carries $L$ to the trivial link. Such a homotopy $h \colon [0,1]\times \coprod_{i=1}^n S^1_i\rightarrow \RR^3$ can be chosen, so that $h(t,.)$ is an embedding except for finitely many times $t_i$, $0<t_1 < \dots <t_i < t_{i+1} < t_k <1$ where $h(t_i,.)$ is an immersion with one double point and no other multiple points, and the link $h(t,.)$ changes exactly by a crossing change when $t$ crosses a $t_i$. (For an alternative elementary proof of this lemma, see \cite[Subsection 7.1]{lescol}).
\eop

In particular, a degree $0$ invariant of $n$-component links of $\RR^3$ must be constant, since it is not allowed to vary under a crossing change.

\begin{exo}
1. Check that  $\CI_1(\CK)=\KK c_0$, where $c_0$ is the constant map that maps any knot to $1$.\\
2. Check that the linking number is a degree $1$ invariant of $2$--component links of $\RR^3$.\\
3. Check that $\CI_1(\CK_2)=\KK c_0 \oplus \KK lk$, where $c_0$ is the constant map that maps any two-component link to $1$.
\end{exo}

\subsection{Introduction to chord diagrams}

Let $f$ be a knot invariant of degree at most $n$. We want to evaluate $f([K;o_1,\dots,o_n])$ where the $o_i$ are disjoint negative crossing changes \pc $\rightarrow$ \nc to be performed on a knot $K$.
Such a $[K;o_1,\dots,o_n]$ is usually represented as a {\em singular knot with $n$ double points\/} that is an immersion of a circle with $n$ transverse double points \singtrefoil, where each double point \doublep  can be desingularized in two ways, the positive one \pc and the negative one \nc, and $K$ is obtained from the singular knot by desingularizing all the crossings in the positive way, which is \righttrefoil in our example. Note that the sign of the desingularization is defined
from the orientation of the ambient space.

Define the {\em chord diagram\/} $\Gamma([K;o_1,\dots,o_n])$ associated to $[K;o_1,\dots,o_n]$ as follows. Draw the preimage of the associated singular knot with $n$ double points as an oriented dashed circle equipped with the $2n$ preimages of the double points and join the pairs of preimages of a double point by a plain segment called a {\em chord\/}.
$$\Gamma(\singtrefoil)=\threexxchord$$
Formally, a {\em chord diagram\/} with $n$ chords is a cyclic order of the $2n$ ends of the $n$ chords, up to a permutation of the chords
and up to exchanging the two ends of a chord.

\begin{lemma}
\label{leminvtyf}
When $f$ is a knot invariant of degree at most $n$, $f([K;o_1,\dots,o_n])$ only depends on $\Gamma([K;o_1,\dots,o_n])$.
\end{lemma}
\bp
Since $f$ is of degree $n$, $f([K;o_1,\dots,o_n])$ is invariant under a crossing change outside the balls of the $o_i$, that is outside the double points of the associated singular knot. Therefore, $f([K;o_1,\dots,o_n])$ only depends on the cyclic order of the $2n$ arcs involved in the $o_i$ on $K$.
\eop

Let $\CD_n$ be the $\KK$-vector space freely generated by the $n$ chord diagrams on $S^1$.

\begin{center}
$\CD_0=\KK \zerochord$, $\CD_1=\KK \onechord$, $\CD_2=\KK \twoischord \oplus \KK \twoxchord$, $\CD_3=\KK \threeischord \oplus \KK \threeparchord \oplus \KK \threexchord \oplus \KK \threexparchord \oplus \KK \threexxchord$.
\end{center}

\begin{lemma}
The map $\phi_n$ from $\CD_n$ to $\frac{\CF_n(\CK)}{\CF_{n+1}(\CK)}$ that maps $\Gamma$ to some $[K;o_1,\dots,o_n]$ whose diagram is $\Gamma$ is well-defined and surjective.
\end{lemma}
\bp Use the arguments of the proof of Lemma~\ref{leminvtyf}.
\eop
$$\phi_3(\threexxchord)=[\singtrefoil].$$

The kernel of the composition of $\phi_n^{\ast}$ and the restriction below
 $$\CI_n(\CK)=\left(\frac{\CF_0(\CK)}{\CF_{n+1}(\CK)}\right)^{\ast}
\rightarrow \left(\frac{\CF_n(\CK)}{\CF_{n+1}(\CK)}\right)^{\ast} \hfl{\phi_n^{\ast}} \CD_n^{\ast}$$
is $\CI_{n-1}(\CK)$. Thus, $\frac{\CI_n(\CK)}{\CI_{n-1}(\CK)}$ injects into $\CD_n^{\ast}$ and $\CI_n(\CK)$ is finite dimensional for all $n$.
Furthermore, $$\frac{\CI_n(\CK)}{\CI_{n-1}(\CK)}=\mbox{Hom}(\frac{\CF_n(\CK)}{\CF_{n+1}(\CK)};\KK).$$

An {\em isolated chord\/} in a chord diagram is a chord between two points of $S^1$ that are consecutive on the circle.

\begin{lemma}
\label{lem14T}
Let $D$ be a diagram on $S^1$ that contains an isolated chord. Then
$\phi_n(D)=0.$
Let $D^1$, $D^2$, $D^3$, $D^4$ be four n-chord diagrams 
that are identical outside three portions of circles where they look like:
$$D^1=\dquatTun ,\;\;\;\; D^2= \dquatTdeux,\;\;\;\; D^3= \dquatTtrois\;\;\;\mbox{and}\;\;\; D^4=\dquatTqua,$$
then
$$\phi_n(-D^1+D^2+D^3-D^4)=0.$$
\end{lemma}
\bp For the first assertion, observe that
$\phi_n(\ischord)=[\pkink]-[\nkink]$.
For the second one, see \cite[Lemma 2.21]{lescol}, for example.
\eop

Let $\CA_n$ denote the quotient of $\CD_n$ by the {\em four--term relation,\/} which is the quotient of $\CD_n$
by the vector space generated by the $(-D^1+D^2+D^3-D^4)$ for all the $4$-tuples $(D^1,D^2,D^3,D^4)$ as above.
Call $(1T)$ the relation that identifies a diagram with an isolated chord with $0$ so that $\CA_n/(1T)$ is the quotient of $\CA_n$ by the vector space generated by diagrams with an isolated chord.

According to Lemma~\ref{lem14T} above, the map $\phi_n$ induces a map 
$$\overline{\phi}_n \colon\CA_n/(1T) \longrightarrow \frac{\CF_n(\CK)}{\CF_{n+1}(\CK)}$$

The fundamental theorem of {\em Vassiliev invariants\/} (which are finite type knot invariants) can now be stated.

\begin{theorem}
\label{thmbn}
There exists a family of linear maps $\left( Z^K_n \colon \CF_0(\CK) \rightarrow \CA_n \right)_{n \in \NN}$ such that
\begin{itemize}
\item $Z^K_n(\CF_{n+1}(\CK))=0$,
\item $Z^K_n$ induces the inverse of $\overline{\phi}_n$ from $\frac{\CF_n(\CK)}{\CF_{n+1}(\CK)}$ to $\CA_n/(1T)$.
\end{itemize}
In particular $\frac{\CF_n(\CK)}{\CF_{n+1}(\CK)} \cong \CA_n/(1T)$ and $\frac{\CI_n(\CK)}{\CI_{n-1}(\CK)}\cong (\CA_n/(1T))^{\ast}$.

\end{theorem}
This theorem has been proved by Kontsevich and Bar-Natan in \cite{barnatan} using the {\em Kontsevich integral\/} $Z^K=(Z^K_n)_{n\in \NN}$ described in \cite{chmutovduzhin} and in \cite[Chapter 8]{chmutovmostov}, for $\KK=\RR$. It is also true when $\KK=\QQ$.

\noindent {\bf Note}
The Kontsevich integral has been generalized to a functor from the category of framed tangles to a category of Jacobi diagrams by Le and Murakami in \cite{LeMur}. Le and Murakami showed how to derive the Reshetikhin-Turaev quantum invariants of framed links in $\RR^3$ defined in \cite{turaevop,reshturibbon} from their functor, in \cite[Theorem~10]{LeMur}.

\subsection{More spaces of diagrams}

\begin{definition}
\label{defunitrivgra}
A {\em uni-trivalent graph\/} $\Gamma$ is a $6$-tuple $(H(\Gamma),E(\Gamma),U(\Gamma),T(\Gamma),p_E,p_V)$ where
$H(\Gamma)$, $E(\Gamma)$, $U(\Gamma)$ and $T(\Gamma)$ are finite sets, which are called the set of half-edges of $\Gamma$, the set of edges of $\Gamma$, the set of univalent vertices of $\Gamma$ and the set of trivalent vertices of $\Gamma$, respectively, $p_E\colon H(\Gamma) \rightarrow E(\Gamma)$ is a two-to-one map (every element of $E(\Gamma)$ has two preimages under $p_E$) and $p_V\colon  H(\Gamma) \rightarrow U(\Gamma) \coprod T(\Gamma)$ is a map such that every element of $U(\Gamma)$ has one preimage under $p_V$ and every element of $T(\Gamma)$ has three preimages under $p_V$, up to isomorphism.
In other words, $\Gamma$ is a set $H(\Gamma)$ equipped with two partitions, a partition into pairs (induced by $p_E$), and a partition into singletons and triples (induced by $p_V$), up to the bijections that preserve the partitions. These bijections are the {\em automorphisms \/} of $\Gamma$.
\end{definition}

\begin{definition}
\label{defdia}
Let $C$ be an oriented one-manifold.
A {\em Jacobi diagram $\Gamma$ with support $C$,\/} also called Jacobi diagram on $C$, is a finite uni-trivalent graph $\Gamma$ 
equipped with an isotopy class of injections $i_{\Gamma}$ of the set $U(\Gamma)$ of univalent vertices of $\Gamma$ into the interior of $C$.
A {\em vertex-orientation\/} of a Jacobi diagram $\Gamma$ is an {\em orientation\/} of every trivalent vertex of $\Gamma$, which is a cyclic order
on the set of the three half-edges which meet at this vertex.
A Jacobi diagram is {\em oriented\/} if it is equipped with a vertex-orientation.
\end{definition}

Such an oriented Jacobi diagram $\Gamma$ is represented by a planar immersion of $\Gamma \cup C$ where the univalent vertices of $U(\Gamma)$ are located at their images under $i_{\Gamma}$, the one-manifold $C$ is represented by dashed lines, whereas
the diagram $\Gamma$ is plain. The vertices are represented by big points. The local orientation of a vertex is represented by the counterclockwise order of the three half-edges
that meet at it.

Here is an example of a picture of a Jacobi diagram $\Gamma$ on the disjoint 
union $M=S^1 \coprod S^1$ of two circles:

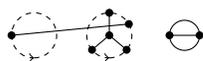
\begin{center}
 \exjactwosonebis
\end{center}

The {\em degree} of such a diagram is 
half the number of all the vertices of $\Gamma$. 

Of course, a chord diagram of $\CD_n$ is a degree $n$ Jacobi diagram on $S^1$ without trivalent vertices.

Let $\CD^t_n(C)$ denote the $\KK$-vector space generated by the degree $n$ oriented Jacobi diagrams
on $C$.
$$\CD^t_1(S^1)=\KK \onechord \oplus \KK \lolli \oplus \KK \tatasone \oplus \KK \tatasonetw \oplus \KK \haltsone$$
Let $\CA^t_n(C)$ \index{N}{AtnC@$\CA^t_n(C)$} denote the quotient of $\CD^t_n(C)$ by the following relations AS, Jacobi and STU:

\begin{center}
AS: \trivAS + \trivASop $=0$\\
Jacobi: \ihxone + \ihxtwo + \ihxthree $=0$\\
STU: \stuy $=$ \stui -\stux
\end{center}

As before, each of these relations relate oriented Jacobi diagrams which are identical outside the pictures
where they are like in the pictures. 

\begin{remark}
Lie algebras provide nontrivial linear maps, called {\em weight systems\/} from $\CA^t_n(C)$ to $\KK$, see \cite{barnatan} and \cite[Section 6]{lescol}. In the weight system constructions, the Jacobi relation for the Lie bracket ensures that the maps defined for oriented Jacobi diagrams factor through the Jacobi relation. In \cite{vogel}, Pierre Vogel proved that the maps associated to Lie (super)algebras are sufficient to detect nontrivial elements of $\CA^t_n(C)$ until degree $15$, and he exhibited a non trivial element of $\CA^t_{16}(\emptyset)$ that cannot be detected by such maps.
The Jacobi relation was originally called IHX by Bar-Natan in \cite{barnatan} because, up to AS, it can be written as $\ihxi=\ihxh - \ihxx.$
\end{remark}

Set $\CA_n(\emptyset)=\CA_n(\emptyset;\KK)=\CA^t_n(\emptyset)$. \index{N}{Anempty@$\CA_n(\emptyset)$}

When $C \neq \emptyset$, let $\CA_n(C)=\CA_n(C;\KK)$ \index{N}{AnC@$\CA_n(C)$} denote the quotient of $\CA^t_n(C)=\CA^t_n(C;\KK)$ by the vector space generated by the diagrams that have at least one connected component without univalent vertices.
Then $\CA_n(C)$ is generated by the oriented Jacobi diagrams whose (plain) connected components contain at least one univalent vertex.

\begin{proposition}
The natural map from $\CD_n$ to $\CA_n(S^1)$ induces an isomorphism from 
$\CA_n$ to $\CA_n(S^1)$.
\end{proposition}
\bsp
The natural map from $\CD_n$ to $\CA_n(S^1)$ factors though $4T$
since, according to $STU$,
$$\dquatTtrois - \dquatTun = \dquatSTU= \dquatTqua -\dquatTdeux$$
in $\CA^t_n(S^1)$.
Since STU allows us to inductively write any oriented Jacobi diagram whose connected components contain at least a univalent vertex as a combination of chord diagrams, the induced map from $\CA_n$ to $\CA_n(S^1)$ is surjective. In order to prove injectivity, one constructs an inverse map. See
\cite[Subsection 3.4]{lescol}.
\eop

 The Le fundamental theorem on {\em finite type invariants of $\ZZ$-spheres\/} is the following one.

\begin{theorem}
\label{thmle}
There exists a family $\left(Z^{LMO}_n \colon \CF_0(\CM) \rightarrow \CA_n(\emptyset) \right)_{n \in \NN}$ of linear maps such that
\begin{itemize}
\item $Z^{LMO}_n(\CF_{2n+1}(\CM))=0$,
\item $Z^{LMO}_n$ induces an isomorphism from $\frac{\CF_{2n}(\CM)}{\CF_{2n+1}(\CM)}$ to $\CA_n(\emptyset)$,
\item $\frac{\CF_{2n-1}(\CM)}{\CF_{2n}(\CM)}=\{0\}$.
\end{itemize}
In particular $\frac{\CF_{2n}(\CM)}{\CF_{2n+1}(\CM)} \cong \CA_n(\emptyset)$ and $\frac{\CI_{2n}(\CM)}{\CI_{2n-1}(\CM)}\cong \CA^{\ast}_n(\emptyset)$.
\end{theorem}
This theorem has been proved by Le \cite{le} using the Le-Murakami-Ohtsuki invariant $Z^{LMO}=(Z^{LMO}_n)_{n\in \NN}$ of \cite{lmo}. As explained in \cite{kuriyalo}, this LMO invariant contains the quantum Witten-Reshetikhin-invariants of rational homology $3$--spheres defined in \cite{reshtu}.

In \cite{moussardAGT}, Delphine Moussard obtained a similar fundamental theorem for {\em finite type invariants of $\QQ$-spheres\/} using
the configuration space integral $Z_{KKT}$ described in \cite{kt}, \cite{lesconst} and in Theorem~\ref{thmmain} below.

As in the knot case, the hardest part of these theorems is the construction of an invariant
$Z=(Z_n)_{n\in \NN}$ that has the required properties. We will define such an invariant by ``counting Jacobi diagram configurations'' in Subsection~\ref{subconstcsi} and explain why it satisfies the required so-called universality properties in Subsection~\ref{subuniv}.

\subsection{Multiplying diagrams}

Set $\CA^t(C)=\prod_{n \in \NN}\CA^t_n(C)$ and $\CA(C)=\prod_{n \in \NN}\CA_n(C)$.

Assume that a one-manifold $C$ is decomposed as a union of two one-manifolds $C = C_1 \cup C_2$ whose interiors in $C$ do not intersect. Define the
{\em product associated to this decomposition\/}:
$$\CA^t(C_1) \times \CA^t(C_2) \longrightarrow \CA^t(C)$$
as the continuous bilinear map which maps $([\Gamma_1],[\Gamma_2])$ to $[\Gamma_1 \coprod \Gamma_2]$, if $\Gamma_1$ is a diagram with support $C_1$ and if $\Gamma_2$ is a diagram with support $C_2$, where $\Gamma_1 \coprod \Gamma_2$ denotes their disjoint union.

In particular, the disjoint union of diagrams turns $\CA(\emptyset)$ into a commutative algebra graded by the degree, and it turns $\CA^t(C)$ into a $\CA(\emptyset)$-module, for any $1$-dimensional manifold $C$.

An orientation-preserving diffeomorphism from a manifold $C$ to another one $C^{\prime}$ induces
an isomorphism from $\CA_n(C)$ to $\CA_n(C^{\prime})$, for all $n$.

Let $I=[0,1]$ be the compact oriented interval. If $I=C$, and if we identify $I$ with $C_1=[0,1/2]$
and with $C_2=[1/2,1]$ with respect to the orientation, then
the above process turns $\CA(I)$ into an algebra where the elements with
non-zero degree zero part admit an inverse.

\begin{proposition}
\label{propdiagcom}
The algebra $\CA([0,1])$ is commutative.
 The projection from $[0,1]$ to $S^1=[0,1]/(0\sim 1)$ induces an isomorphism from 
$\CA_n([0,1])$ to $\CA_n(S^1)$ for all $n$, so that $\CA(S^1)$ inherits a commutative algebra structure from this isomorphism.
The choice of a connected component $C_j$ of $C$ equips $\CA(C)$ with an {\em $\CA([0,1])$-module structure $\sharp_j$  \/}, induced by the inclusion from $[0,1]$ to a little part of $C_j$ outside the vertices, and the insertion of diagrams with support $[0,1]$ there.
\end{proposition}

In order to prove this proposition, we present a useful trick in diagram spaces.

First adopt a convention.
So far, in a diagram picture, or in a chord diagram picture, the plain edge of a univalent vertex, has always been attached on the left-hand side of the oriented one-manifold. Now, if $k$ plain edges are attached on the other side on a diagram picture, then we agree that the corresponding represented element of $\CA^t_n(M)$ is $(-1)^k$ times the underlying
diagram. With this convention, we have the new antisymmetry relation in $\CA^t_n(M)$:

$$\asunor +\asunorc=0$$
and we can draw the STU relation like the Jacobi relation:
$$\stubone + \stubtwo + \stubthree=0.$$
\begin{lemma}
\label{lemcom}
Let $\Gamma_1$ be a Jacobi diagram with support $C$. Assume that $\Gamma_1 \cup C$ is immersed in the plane so that  $\Gamma_1 \cup C$ meets an open annulus $A$ embedded in the plane
exactly along $n+1$ embedded arcs $\alpha_1$, $\alpha_2$, \dots, $\alpha_n$ and $\beta$,
and one vertex $v$ so that:
\begin{enumerate}
\item The $\alpha_i$ may be dashed or plain, they run from a boundary component of $A$ to the other one,
\item $\beta$ is a plain arc which runs from the boundary of $A$ to
$v \in \alpha_1$,
\item The bounded component $D$ of the complement of $A$ does not contain a boundary point of $C$. 
\end{enumerate}
Let $\Gamma_i$ be the diagram obtained from $\Gamma_1$ by attaching the 
endpoint $v$ of $\beta$ to $\alpha_i$ instead of $\alpha_1$ on the same side,
where
the side of an arc is its side when going from the outside boundary component
of $A$ to the inside one $\partial D$.
Then $\sum_{i=1}^n\Gamma_i=0$ in $\CA^t(C)$.
\end{lemma}
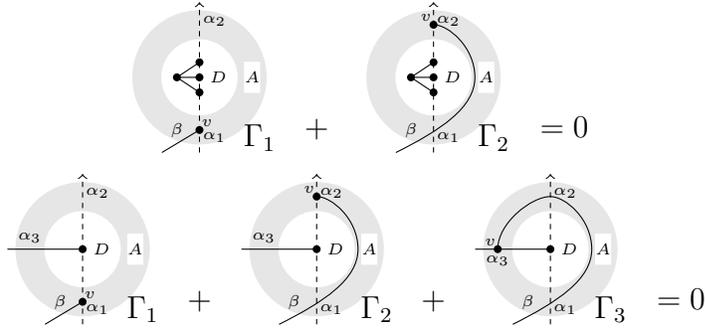
\begin{examples}
$$\annulone + \annultwo =0$$
$$\annulthree + \annulfour + \annulfive =0$$
\end{examples}

\bp
The second example shows that the STU relation is equivalent to this relation when the bounded component $D$ of $\RR^2 \setminus A$ intersects $\Gamma_1$ in the neighborhood of a univalent vertex on $C$. Similarly, the Jacobi relation is easily seen as given by this relation  when
 $D$ intersects $\Gamma_1$ in the neighborhood of a trivalent vertex. Also note that
AS  corresponds to the case when $D$  intersects $\Gamma_1$ along a
dashed or plain arc.
Now for the Bar-Natan \cite[Lemma 3.1]{barnatan} proof. See also \cite[Lemma 3.3]{vogel}. Assume without loss that $v$ is always attached on the right-hand-side of 
the $\alpha$'s. Add to the sum the trivial (by Jacobi and STU) contribution of the sum of the diagrams
obtained from $\Gamma_1$ by attaching $v$ to each of the three (dashed or plain) half-edges of each vertex $w$
of $\Gamma_1 \cup C$ in $D$ on the left-hand side when the half-edges are oriented towards $w$. Now, group the terms of the obtained sum by edges of $\Gamma_1 \cup C$
where $v$ is attached, and observe that the sum is zero edge by edge by AS.
\eop

\noindent{\sc Proof of Proposition~\ref{propdiagcom}:}
To each choice of a connected component $C_j$ of $C$, we associate an {\em $\CA(I)$-module structure $\sharp_j$ on $\CA(C)$\/}, which is given by the continuous bilinear map:
$$\CA(I) \times \CA(C) \longrightarrow \CA(C)$$
such that:
If  $\Gamma^{\prime}$ is a diagram with support $C$ and if $\Gamma$ is a diagram with support $I$, then $([\Gamma],[\Gamma^{\prime}])$ is mapped to the class of the diagram obtained by inserting $\Gamma$ along $C_j$ outside the vertices of $\Gamma$, according to the given orientation. 
For example, $$\dmer \cirt = \cirtd =\cirtdp $$
As shown in the first example that illustrates  Lemma \ref{lemcom}, the independence of the choice of the insertion locus is a consequence of Lemma~\ref{lemcom}, where $\Gamma_1$ is the disjoint union $\Gamma \coprod \Gamma^{\prime}$ and $\Gamma_1$ intersects $D$ along $\Gamma \cup I$.
This also proves that $\CA(I)$ is a commutative algebra.
Since the morphism from $\CA(I)$ to $\CA(S^1)$ induced by the identification of the two endpoints of $I$ amounts to quotient out $\CA(I)$ by the relation that identifies two diagrams that are obtained from one another by moving the nearest univalent vertex to an endpoint of I near the other endpoint, a similar application of Lemma \ref{lemcom} also proves that this morphism is an isomorphism from $\CA(I)$ to $\CA(S^1)$. (In this application,
$\beta$ comes from the inside boundary of the annulus.)
\eop

\newpage
\section{Configuration space construction of universal finite type invariants}
\label{secconstcsi}
\setcounter{equation}{0}

In this section, we finally describe the promised invariants, which generalize both the linking number and $\Theta$. These invariants count configurations of Jacobi diagrams with support some link, in an asymptotic rational homology $\RR^3$. In Subsection~\ref{subdefconfspace}, we introduce the relevant configuration spaces. In Subsection~\ref{subconfint}, we define integrals over these spaces from propagating forms. The wanted invariants are obtained by combining these integrals in Subsection~\ref{subconstcsi}. These integrals will be expressed in terms of algebraic intersections, which involve propagating chains, in Subsection~\ref{subrat}. Important universality properties of the constructed invariants are presented in Subsection~\ref{subuniv}.

\subsection{Configuration spaces of links in \texorpdfstring{$3$}{3}--manifolds}
\label{subdefconfspace}

Let $(\check{M},\tau)$ be an asymptotic rational homology $\RR^3$.

Let $C$ be a disjoint union of $k$ circles $S^1_i$, $i \in \underline{k}$ and let 
$$L : C \longrightarrow \check{M}$$ denote a $C^{\infty}$ embedding from $C$ to $\check{M}$. Let $\Gamma$ be a Jacobi diagram with support $C$. Let $U=U(\Gamma)$ denote the set of univalent vertices of $\Gamma$, and let $T=T(\Gamma)$ denote the set of trivalent vertices of $\Gamma$. A {\em configuration\/} of $\Gamma$ is an embedding 
$$c:U \cup T \hookrightarrow \check{M}$$
whose restriction $c_{|U}$ to $U$ may be written as $L \circ j$ for some injection 
$$j:U \hookrightarrow C$$
in the given isotopy class $[i_{\Gamma}]$ of embeddings of $U$ into the interior of $C$. Denote the set of these configurations by
$\check{C}(L;\Gamma)$, $$\check{C}(L;\Gamma)=\left\{c:U \cup T \hookrightarrow \check{M} \; ; 
\exists j \in [i_{\Gamma}], c_{|U}=L \circ j\right\}.$$ 
In $\check{C}(L;\Gamma)$, the univalent vertices move along $L(C)$ while the trivalent vertices move in the ambient space, and $\check{C}(L;\Gamma)$ is naturally an open submanifold of 
$C^U \times \check{M}^T$.

An {\em orientation\/} of a set of cardinality at least $2$ is a total order of its elements up to an even permutation.

Cut each edge of $\Gamma$ into two half-edges. When an edge is oriented, define its {\em first\/} half-edge and its {\em second\/} one, so that following the orientation of the edge, the first half-edge is met first.
Recall that $H(\Gamma)$ denotes the set of half-edges of $\Gamma$.

\begin{lemma}
When $\Gamma$ is equipped with a vertex-orientation, orientations of the manifold $\check{C}(L;\Gamma)$ are in canonical
one-to-one correspondence with orientations of the set $H(\Gamma)$.
\end{lemma}
\bp
Since $\check{C}(L;\Gamma)$ is naturally an open submanifold of 
$C^U \times \check{M}^T$, it inherits $\RR^{\sharp U + 3\sharp T}$-valued charts from $\RR$-valued orientation-preserving charts of $C$ and $\RR^3$-valued orientation-preserving charts of $\check{M}$. In order to define the orientation of $\RR^{\sharp U + 3\sharp T}$, one must identify its factors and order them (up to even permutation).
Each of the factors may be labeled by an element of $H(\Gamma)$: the $\RR$-valued local coordinate of an element of $C$ corresponding to the image under $j$ of an element of $U$ sits in the factor labeled by the half-edge of $U$;
the $3$ cyclically ordered (by the orientation of $\check{M}$) $\RR$-valued local coordinates of the image under a configuration $c$ of an element of $T$ live in the factors labeled by the three half-edges that are cyclically ordered by the vertex-orientation of $\Gamma$, so that the cyclic orders match.
\eop

The dimension of $\check{C}(L;\Gamma)$ is
$$\sharp U(\Gamma) + 3\sharp T(\Gamma)= 2 \sharp E(\Gamma)$$
where $E=E(\Gamma)$ denotes the set of edges of $\Gamma$.
Since $n=n(\Gamma)=\frac12(\sharp U(\Gamma) +\sharp T(\Gamma))$, $$\sharp E(\Gamma)=3n-\sharp U(\Gamma).$$

\subsection{Configuration space integrals}
\label{subconfint}

A {\em numbered\/} degree $n$ Jacobi diagram is a degree $n$ Jacobi diagram $\Gamma$ whose edges are oriented, equipped with an injection $j_E \colon E(\Gamma) \hookrightarrow \underline{3n}$.
Such an injection numbers the edges.
Note that this injection is a bijection when $U(\Gamma)$ is empty.
Let $\CD^e_n(C)$ \index{N}{De@$\CD^e_n(C)$} denote the set of numbered degree $n$ Jacobi diagrams with support $C$ without {\em looped edges\/} like \loopedge.

Let $\Gamma$ be a numbered degree $n$ Jacobi diagram. 
The orientations of the edges of $\Gamma$ induce the following orientation of the set $H(\Gamma)$ of half-edges of $\Gamma$: Order $E(\Gamma)$ arbitrarily, and order the half-edges as 
(First half-edge of the first edge, second half-edge of the first edge, \dots, second half-edge of the last edge). The induced orientation is called the \emph{edge-orientation} of $H(\Gamma)$. Note that it does not depend on the order of $E(\Gamma)$. Thus, as soon as $\Gamma$ is equipped with a vertex-orientation $o(\Gamma)$, the edge-orientation of $\Gamma$ orients $\check{C}(L;\Gamma)$.

An edge $e$ oriented from a vertex $v_1$ to a vertex $v_2$ of $\Gamma$ induces the following canonical map $$\begin{array}{llll}p_e \colon &\check{C}(L;\Gamma) &\rightarrow &C_2(M)\\
           & c & \mapsto & (c(v_1),c(v_2)).\end{array}$$

For any $i \in \underline{3n}$, let $\omega(i)$ be a propagating form of $(C_2(M),\tau)$.
Define $$I(\Gamma,o(\Gamma),(\omega(i))_{i \in \underline{3n}})=\int_{(\check{C}(L;\Gamma),o(\Gamma))} \bigwedge_{e \in E(\Gamma)}p_e^{\ast}(\omega(j_E(e)))$$
where $(\check{C}(L;\Gamma),o(\Gamma))$ denotes the manifold $\check{C}(L;\Gamma)$ equipped with the orientation induced by the vertex-orientation $o(\Gamma)$ and by the edge-orientation of $\Gamma$.

The convergence of this integral is a consequence of the following proposition, which will be proved in Subsection~\ref{subcompconf}.

\begin{proposition}
\label{propcompext}
 There exists a smooth compactification $C(L;\Gamma)$ of $\check{C}(L;\Gamma)$ where the maps $p_e$ smoothly extend.
\end{proposition}

According to this proposition, $\bigwedge_{e \in E(\Gamma)}p_e^{\ast}(\omega(j_E(e)))$ smoothly extends to $C(L;\Gamma)$, and
$$\int_{(\check{C}(L;\Gamma),o(\Gamma))} \bigwedge_{e \in E(\Gamma)}p_e^{\ast}(\omega(j_E(e)))=
\int_{(C(L;\Gamma),o(\Gamma))} \bigwedge_{e \in E(\Gamma)}p_e^{\ast}(\omega(j_E(e))).$$

\begin{examples}
For any three propagating forms $\omega(1)$, $\omega(2)$ and $\omega(3)$ of $(C_2(M),\tau)$,
 $$I(\haltere,(\omega(i))_{i \in \underline{3}})=lk(K_i,K_j)$$
and $$I(\tata,(\omega(i))_{i \in \underline{3}})=\Theta(M,\tau)$$
for any numbering of the (plain) diagrams (exercise).
\end{examples}

Let us now study the case of $I(\onechordsmalljnum,(\omega(i))_{i \in \underline{3}})$, which depends on the chosen propagating forms, and on the diagram numbering.

A {\em dilation\/} is a homothety with positive ratio.

Let $U^+K_j$ denote the fiber space over $K_j$ made of the tangent vectors to the knot $K_j$ of $\check{M}$ that orient $K_j$, up to dilation. The fiber of $U^+K_j$ is made of one point, so that the total space of this {\em unit positive tangent bundle to $K_j$\/} is $K_j$. Let $U^-K_j$ denote the fiber space over $K_j$ made of the opposite tangent vectors to $K_j$, up to dilation.

For a knot $K_j$ in $\check{M}$, $$\check{C}(K_j;\onechordsmallj)=\{(K_j(z),K_j(z\exp(i\theta))); (z,\theta) \in S^1 \times ]0,2\pi[\} .$$
Let $C_j={C}(K_j;\onechordsmallj)$ be the closure of $\check{C}(K_j;\onechordsmallj)$ in $C_2(M)$. This closure is
diffeomorphic to $S^1 \times [0,2\pi]$ where $S^1 \times 0$ is identified with $U^+K_j$, $S^1 \times \{2\pi\}$ is identified with $U^-K_j$ and $\partial C(K_j;\onechordsmallj) =U^+K_j - U^-K_j$.

\begin{lemma}

\label{lemvaritheta}
For any $i\in \underline{3}$, let\/ $\omega(i)$ and\/ $\omega^{\prime}(i)=\omega(i)+d\eta(i)$ be propagating forms of $(C_2(M),\tau)$, where $\eta(i)$ is a one-form on $C_2(M)$.

$$I(\onechordsmalljnum,(\omega^{\prime}(i))_{i \in \underline{3}})-I(\onechordsmalljnum,(\omega(i))_{i \in \underline{3}})=\int_{U^+K_j}\eta(k) -\int_{U^-K_j}\eta(k). $$
\end{lemma}
\bp Apply the Stokes theorem to $\int_{C_j}(\omega^{\prime}(k)-\omega(k))=\int_{C_j}d\eta(k)$. \eop

\begin{exo}
 Find a knot $K_j$ of $\RR^3$ and a form $\eta(k)$ of $C_2(\RR^3)$ such that the right-hand side of Lemma~\ref{lemvaritheta} does not vanish. (Use Lemma~\ref{lemetactwo}, hints can be found in Subsection~\ref{subsecstraight}.)
\end{exo}

Say that a propagating form $\omega$ of $(C_2(M),\tau)$ is {\em homogeneous\/} if its restriction to $\partial C_2(M)$ is $p_{\tau}^{\ast}(\omega_{S^2})$ for the homogeneous volume form $\omega_{S^2}$ of $S^2$ of total volume $1$.

\begin{lemma}
\label{lemdefItheta}
For any $i\in \underline{3}$, let $\omega(i)$ be a homogeneous propagating form of $(C_2(M),\tau)$.
Then $I(\onechordsmalljnum,(\omega(i))_{i \in \underline{3}})$ does not depend on the choices of the $\omega(i)$, it is denoted by $I_{\theta}(K_j,\tau)$.\index{N}{Itheta@$I_{\theta}(K_j,\tau)$}
\end{lemma}
\bp Apply Lemma~\ref{lemetactwo} with $\eta_A=0$, so that $\eta(k)=0$ in Lemma~\ref{lemvaritheta}.
\eop

\subsection{An invariant for links in \texorpdfstring{$\QQ$}{Q}--spheres from configuration spaces}
\label{subconstcsi}

Let $\KK=\RR$.
Let
$[\Gamma,o(\Gamma)]$ denote the class in $\CA^t_n(C)$ of a numbered Jacobi diagram $\Gamma$ of $\CD^e_n(C)$ equipped with a vertex-orientation $o(\Gamma)$,
then $I(\Gamma,o(\Gamma),(\omega(i))_{i \in \underline{3n}})[\Gamma,o(\Gamma)] \in \CA^t_n(C)$ is independent of the orientation of $o(\Gamma)$, it will be simply denoted by 
$I(\Gamma,(\omega(i))_{i \in \underline{3n}})[\Gamma]$.

\begin{theorem}
\label{thmmain}
Let $(\check{M},\tau)$ be an asymptotic rational homology $\RR^3$.
Let $L \colon \coprod_{j=1}^kS^1_j \hookrightarrow \check{M}$ be an embedding.
For any $i\in \underline{3n}$, let $\omega(i)$ be a homogeneous propagating form of $(C_2(M),\tau)$.

Set \index{N}{Z@$Z$}$$Z_n(L,\check{M},\tau)=\sum_{\Gamma \in \CD^e_n(C)}\frac{(3n-\sharp E(\Gamma))!}{(3n)!2^{\sharp E(\Gamma)}}I(\Gamma,(\omega(i))_{i \in \underline{3n}})[\Gamma] \in \CA_n^t(\coprod_{j=1}^kS^1_j).$$
Then $Z_n(L,\check{M},\tau)$ is independent of the chosen
$\omega(i)$, it only depends on the diffeomorphism class of $(M,L)$, on $p_1(\tau)$ and on the $I_{\theta}(K_j,\tau)$, for the components $K_j$ of $L$.

More precisely,
set $$Z(L,\check{M},\tau)=(Z_n(L,\check{M},\tau))_{n\in \NN} \in \CA^t(\coprod_{j=1}^kS^1_j).$$
There exist two constants $\alpha \in \CA(S^1;\QQ)$ and $\ansothree \in \CA(\emptyset;\QQ)$ such that $$\exp(-\frac14 p_1(\tau)\ansothree)\prod_{j=1}^k\left(\exp(-I_{\theta}(K_j,\tau)\alpha)\sharp_j\right) Z(L,\check{M},\tau)=Z(L,M)$$
only depends on the diffeomorphism class of $(M,L)$. Here $\exp(- I_{\theta}(K_j)\alpha)$ acts on $Z(L,\check{M},\tau)$, on the copy $S^1_j$ of $S^1$ as indicated by the subscript $j$.

$$Z(L,M) \in \CA^t(\coprod_{j=1}^kS^1_j;\QQ).$$

Furthermore, if $\check{M}=\RR^3$, then the projection $Z^u(L,S^3)$ \index{N}{Zu@$Z^u$} of $Z(L,S^3)$ on $\CA(\coprod_{j=1}^kS^1_j)$ is a {\em universal\/} finite type invariant of links in $\RR^3$, i.e. $Z_n^u$ satisfies the properties stated for $Z_n^K$ in Theorem~\ref{thmbn}. It is the configuration space invariant studied by Altsch\"uler, Freidel \cite{af}, Dylan Thurston \cite{thurstonconf}, Sylvain Poirier \cite{poirier} and others \footnote{after work of many people including Witten \cite{witten}, Guadagnini, Martellini and Mintchev~\cite{gmm}, Kontsevich~\cite{ko,Kon}, Bott and Taubes \cite{botttaubes}, Bar-Natan~\cite{barnatanper}, Axelrod and Singer~\cite{axelsingI,axelsingII}}.
If $k=0$, then $Z(\emptyset,M)$ is the Kontsevich configuration space invariant $Z_{KKT}(M)$, which is a universal invariant for $\ZZ$-spheres according to a theorem of Kuperberg and Thurston \cite{kt,lessumgen}, and which was completed to a universal finite type invariant for $\QQ$-spheres by Delphine Moussard \cite{moussardAGT}.
\end{theorem}

The proof of this theorem is sketched in Section~\ref{secproofinv}.

Under its assumptions, let $\omega_0$ be a homogeneous propagating form of $(C_2(M),\tau)$, let $\iota$ be the involution of $C_2(M)$ that permutes two elements in $\check{M}^2 \setminus \mbox{diagonal}$, set $\omega=\frac12(\omega_{0}-\iota_{\ast}(\omega_0))$, and set $\omega(i)=\omega$ for any $i$.

Let $\mbox{Aut}(\Gamma)$ be the set of automorphisms of $\Gamma$, which is the set of permutations of the half-edges that map a pair of half-edges of an edge to another such and a triple of half-edges that contain a vertex to another such, and that map half-edges of univalent vertices on a component $K_j$ to half-edges of univalent vertices on $K_j$ so that the cyclic order among such vertices is preserved.
Set \index{N}{betaGamma@$\beta_{\Gamma}$}
$$\beta_{\Gamma}=\frac{(3n-\sharp E(\Gamma))!}{(3n)!2^{\sharp E(\Gamma)}}.$$ 

Then $$\sum_{\Gamma \in \CD^e_n(C)}\beta_{\Gamma}I(\Gamma,(\omega(i))_{i \in \underline{3n}})[\Gamma]=\sum_{\Gamma \,\mbox{\tiny unnumbered, unoriented}} \frac1{\sharp \mbox{Aut}(\Gamma)}I(\Gamma,(\omega)_{i \in \underline{3n}})[\Gamma]$$
where the sum of the right-hand side runs over the degree $n$ Jacobi diagrams on $C$ without looped edges.

Indeed, for a numbered graph $\Gamma$, there are $\frac{1}{\beta_{\Gamma}}$ ways of renumbering it, and $\sharp \mbox{Aut}(\Gamma)$ of them will produce the same numbered graph.

\subsection{On the universality proofs}
\label{subuniv}

\begin{theorem}
\label{thmunivone} Let $y $, $z \in \NN$. Recall $\underline{y}=\{1,2,\dots,y\}$. Set $(\underline{z}+y)=\{y+1,y+2,\dots,y+z\}$. Let $\check{M}$ be an asymptotically standard $\QQ$-homology $\RR^3$.
 Let $L$ be a link in $\check{M}$. Let $(B_b)_{\in \underline{y}}$ be a collection of pairwise disjoint balls in $\check{M}$ such that every $B_b$ intersects $L$ as a ball of a crossing change that contains a positive crossing $c_b$, and let $L((B_b)_{b \in \underline{y}})$ be the link obtained by changing the positive crossings $c_b$ to negative crossings. Let $(A_a)_{a\in(\underline{z}+y)}$ be a collection of pairwise disjoint rational homology handlebodies in 
$\check{M} \setminus (L \cup_{b=1}^{y} B_b)$. Let $(A^{\prime}_a/A_a)$ be rational LP surgeries in $\check{M}$. Set
$X=[M,L;(A^{\prime}_a/A_a)_{a \in(\underline{z}+y)},(B_b,c_b)_{b \in \underline{y}}]$ and $$Z_{n}(X)=  \sum_{I\subset \underline{y+z}}(-1)^{\sharp I}Z_n\left(L((B_b)_{b \in I \cap \underline{y}}),M((A^{\prime}_a/A_a)_{a \in I \cap  (\underline{z}+y)})\right).$$
If $2n <2 y +  z $, then $Z_{n}(X)$ vanishes .
\end{theorem}
\bsp
As in \cite{lessumgen}, one can use (generalized) propagators for the $M((A^{\prime}_a/A_a)_{a \in I \cap  (\underline{z}+y)})$ that coincide for different $I$ wherever it makes sense (for example, for configurations that do not involve points in surgered pieces $A_a$). See also \cite{kt}.
Then contributions to the alternate sum of the integrals over parts that do not involve at least one point in an $A_a$ or in an $A^{\prime}_a$, for all $a$ cancel.
Assume that every crossing change is performed by moving only one strand.
Again, contributions to the alternate sum of the integrals that do not involve at least one point on a moving strand cancel.
Furthermore, if the moving strand of $c_b$ is moved very slightly, and if no other vertex is constrained to lie on the other strand in the ball of the crossing change, then the alternate sum is close to zero. Thus in order to produce a contribution to the alternate sum, a graph must have at least $(2 y +  z)$ vertices. See \cite{af} or \cite[Section~5.4]{lescol}, and \cite[Section~3]{lessumgen} for more details.
\eop

This implies that $Z_n^u$ is of degree at most $n$ for links in $\RR^3$, and that $Z_n$ is of degree at most $2n$ for $\ZZ$-spheres or $\QQ$-spheres.

Now, under the hypotheses of Theorem~\ref{thmunivone},
assume that $A_a$ is the standard genus $3$ handlebody with three handles with meridians $m^{(a)}_j$ and longitudes $l^{(a)}_j$ such that $\langle m^{(a)}_i, \ell^{(a)}_j\rangle_{\partial A_a} =\delta_{ij}$. See  $A_a$ as a thickening of the trivalent graph below.

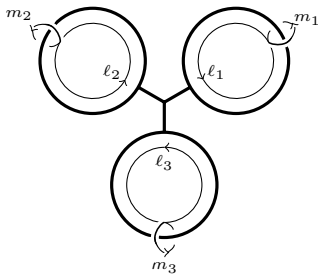
\begin{center}
 \begin{tikzpicture} \useasboundingbox (-2.5,-2.5) rectangle (2.5,2.2);
\draw [very thick] (30:.4) -- (0,0) -- (0,-.4) (150:.4) -- (0,0) (30:1.1) circle (.7) (150:1.1) circle (.7)  (-90:1.1) circle (.7);
\draw [->] (30:.6) arc (-150:210:.5);
\draw [->] (150:.6) arc (-30:330:.5);
\draw [->] (-90:.6) arc (-270:90:.5);
\draw (30:.8) node{\tiny $\ell_1$} (150:.8) node{\tiny $\ell_2$} (-90:.8) node{\tiny $\ell_3$} (30:2.2) node{\tiny $m_1$} (149:2.25) node{\tiny $m_2$} (-90:2.2) node{\tiny $m_3$};
\draw [->,draw=white,double=black,very thick] (30:1.6) .. controls (29:1.6) and (26:1.7) .. (26:1.8) .. controls (26:1.9) and (29:2) .. (30:2);
\draw (26:1.8) .. controls (26:1.9) and (29:2) .. (30:2);
\draw [>-] (30:2) .. controls (31:2) and (34:1.95) .. (34:1.9);
\draw (30:1.6) .. controls (31:1.6) and (34:1.65) .. (34:1.7);
\draw [->,draw=white,double=black,very thick] (150:1.6) .. controls (149:1.6) and (146:1.7) .. (146:1.8) .. controls (146:1.9) and (149:2) .. (150:2);
\draw (146:1.8) .. controls (146:1.9) and (149:2) .. (150:2);
\draw [>-] (150:2) .. controls (151:2) and (154:1.95) .. (154:1.9);
\draw (150:1.6) .. controls (151:1.6) and (154:1.65) .. (154:1.7);
\draw [->,draw=white,double=black,very thick] (-90:1.6) .. controls (-91:1.6) and (-94:1.7) .. (-94:1.8) .. controls (-94:1.9) and (-91:2) .. (-90:2);
\draw (-94:1.8) .. controls (-94:1.9) and (-91:2) .. (-90:2);
\draw [>-] (-90:2) .. controls (-89:2) and (-86:1.95) .. (-86:1.9);
\draw (-90:1.6) .. controls (-89:1.6) and (-86:1.65) .. (-86:1.7);
\end{tikzpicture}
\end{center}

Also assume that $A^{\prime}_a$ is an integer homology handlebody.
In $A_a \cup_{\partial A_a}(- A^{\prime}_a)$, there is a surface $S_j$ such that $\partial (S_j \cap A_a)=m^{(a)}_j$. Assume that $\langle S_1,S_2,S_3\rangle_{A_a \cup_{\partial A_a}(- A^{\prime}_a)} =1$. ( For example, choose $A^{\prime}_a$ such that $A_a \cup_{\partial A_a}(- A^{\prime}_a)=(S^1)^3$, like in the case of the Matveev Borromean surgery of \cite{matveev}. )
Assume that the $l^{(a)}_j$ bound surfaces $D^{(a)}_j$ in $\check{M}$.

Assume that the collection of surfaces $\{D^{(a)}_j\}_{a \in(\underline{z}+y),j \in\underline{3}}$ reads 
$\{D_{p,1}\}_{p \in \underline{P}} \sqcup \{D_{p,2}\}_{p \in \underline{P}}$ so that \\
for any $q \in \underline{P}$, for $\delta \in \underline{2}$, if $D_{q,\delta}=D^{(a(q,\delta))}_{j(q,\delta)}$, the interior of $D_{q,\delta}$ intersects 
$$L \cup \bigcup_{a \in(\underline{z}+y)} \left(A_a \cup \cup_{j \in\underline{3}, D^{(a)}_j \neq D_{q,\delta}} D^{(a)}_j \right) \cup \cup_{b \in \underline{y}} (B_b)$$
 only in $A_{a(q,3-\delta)} \cup D^{(a(q,3-\delta))}_{j(q,3-\delta)}$.
 
Note that $\langle D_{q,\delta},\ell^{(a(q,3-\delta))}_{j(q,3-\delta)} \rangle_{M}=lk(\partial D_{q,1} ,\partial D_{q,2})$.

\begin{example}
 Note that these assumptions are realised in the following case.
Start with an embedding of a Jacobi diagram $\Gamma$ whose univalent vertices belong to chords (plain edges between two univalent vertices) on $\cup_{i=1}^k S^1_i$ in $\check{M}$. Assume that the trivalent vertices of $\Gamma$ are labeled in $(\underline{z}+y)$, and assume that its chords are labeled in $\underline{y}$. Apply the following operations\\
replace edges \edget without univalent vertices by \edgehopf ,\\
replace a chord \edgeuu labeled by $b$ by a crossing change $c_b$ $\noedgeuu \rightarrow \xchangeuu $ in a ball $B_b$ that is a neighborhood of the plain edge.\\
Thicken the trivalent graph \smally associated to the trivalent vertex labeled by $a$, and call it $A_a$. Then the surfaces $D^{(a)}_j$ are the disks bounded by the small loops of \smally.
\end{example}

Conversely, under the assumptions before the example, define the following vertex-oriented Jacobi diagram $$\Gamma([M,L;(A^{\prime}_a/A_a)_{a \in (\underline{z}+y)},(B_b,c_b)_{b \in \underline{y}}])$$ on $\cup_{i=1}^k S^1_i$, with 
 \begin{itemize}
  \item two univalent vertices joined by a chord for each crossing change ball $B_b$ at the corresponding places on $\cup_{i=1}^k S^1_i$ (in $L^{-1}(B_b)$), 
\item one trivalent vertex for each $A_a$, where the three adjacent half-edges of the vertex correspond to the three $D^{(a)}_j$, with the fixed cyclic order,
 \end{itemize}
such that any pair of half-edges corresponding to some $D_{p,1}$ and its friend $D_{p,2}$ forms an edge between two trivalent vertices.

\begin{theorem}
\label{thmkeyuniv}
Under the assumptions above, 
let $X=[M,L;(A^{\prime}_a/A_a)_{a \in(\underline{z}+y)},(B_b,c_b)_{b \in \underline{y}}]$.
When $2n=2y+z$,
$$Z_{n}(X)=\left(\prod_{p \in \underline{P}}lk(\partial D_{p,1} ,\partial D_{p,2})\right) [\Gamma(X)] \;\mbox{mod 1T}\;\left(\mbox{or in}\; \frac{\CA_n^t(\coprod_{j=1}^kS^1_j)}{(1T)}\right).$$
\end{theorem}
\bsp
When $z=0$, the proof of Theorem~\ref{thmunivone} can be pushed further in order to prove the result like in \cite{af} or \cite[Section 5.4]{lescol}. 
In general, when $y=0$, it is a consequence of the main theorem in \cite{lessumgen} (Theorem 2.4).
The general result can be obtained by mixing the arguments of \cite[Section 3]{lessumgen} with the arguments of the link case.
\eop

This theorem is the key to proving the universality of $Z^u$ among Vassiliev invariants for links in $\RR^3$ and to proving the universality of $Z$ among finite type invariants of $\ZZ$-spheres. This universality implies that all finite type invariants factor through $Z$.

\begin{remark}
Theorem~\ref{thmkeyuniv} with $Z^{LMO}$ instead of $Z$ is proved in \cite{le}, when $y=0$, when the $(A^{\prime}_a/A_a)$ are Matveev's Borromean surgeries and when the $D^{(a)}_j$ are disks such that $lk(\partial D_{p,1} ,\partial D_{p,2})=1$. Then the main theorem of \cite{aucles} implies Theorem~\ref{thmkeyuniv} with $Z^{LMO}$ instead of $Z$, when $y=0$ and when the $A_a$ and the $A^{\prime}_a$ are integral homology handlebodies. 
\end{remark}

\newpage
\section{Compactifications, anomalies, proofs and questions}
\label{secproofinv}
\setcounter{equation}{0}

In this section, we state Theorem~\ref{thmmainstraight}. This is another version of Theorem~\ref{thmmain}, which leads
to a definition of $Z$ involving algebraic intersections rather than integrals in Subsection~\ref{subrat}. It is based on the concept of \emph{straight links} introduced in Subsection~\ref{subsecstraight}.

This section also contains sketches of proofs of Theorems~\ref{thmmain} and \ref{thmmainstraight}.
We begin with the introduction of appropriate compactifications of configuration spaces to justify the convergence of our integrals stated in Proposition~\ref{propcompext}.

\subsection{Compactifications of configuration spaces}
\label{subcompconf}

Let ${\nfinset}$ be a finite set. See the elements of $M^{\nfinset}$ as maps $m \colon {\nfinset} \rightarrow M$.

For a non-empty $I \subseteq {\nfinset}$, let $E_I$ be the set of maps that map $I$ to $\infty$.
For $I \subseteq {\nfinset}$ such that $\sharp I \geq 2$, let $\Delta_I$ be the set of maps that map $I$ to a single element of $M$.
When $I$ is a finite set, and when $V$ is a vector space of positive dimension, $\check{S}_I({V})$ denotes the space of injective maps from $I$ to $V$ up to translation and dilation. When $\sharp I \geq 2$, $\check{S}_I({V})$ embeds in the compact space ${S}_I({V})$ of non-constant maps from $I$ to $V$ up to translation and dilation.

\begin{lemma}
 The fiber of the unit normal bundle to $\Delta_I$
 in $M^{\nfinset}$ over a configuration $m$ is ${S}_I(T_{m(I)}M)$.
\end{lemma}
\bp Exercise. \eop

Let $\check{C}_{\nfinset}({M})$ denote the space of injective maps from ${\nfinset}$ to $\check{M}$.

Define a compactification ${C}_{\nfinset}({M})$ of $\check{C}_{\nfinset}({M})$ by generalizing the previous construction of $C_2(M)=C_{\underline{2}}(M)$ as follows.

Start with $M^{\nfinset}$.
Blow up $E_{{\nfinset}}$, which is the point $m=\infty^{\nfinset}$ such that $m^{-1}(\infty)={\nfinset}$.\\
Then for $k=\sharp {\nfinset},\sharp {\nfinset}-1, \dots, 3, 2$, in this decreasing order, successively blow up 
the (closures of the preimages under the composition of the previous blow-down maps of the) $\Delta_{I}$ such that $\sharp I=k$ (choosing an arbitrary order among them) and, next, the (closures of the preimages under the composition of the previous blow-down maps of the) $E_J$ such that $\sharp J=k-1$ (again, choosing an arbitrary order among them).

\begin{lemma}
\label{lemcompconf}
The successive manifolds that are blown-up in the above process are smooth and transverse to the boundaries.
The manifold ${C}_{\nfinset}({M})$ is a smooth compact $(3\sharp {\nfinset})$-manifold independent of the possible order choices in the process. For $i,j \in {\nfinset}$, $i\neq j$,
the map $$\begin{array}{llll}p_{i,j} \colon & \check{C}_{\nfinset}({M}) &\rightarrow &C_2(M) \\ & m & \mapsto &(m(i),m(j))\end{array}$$  smoothly extends to ${C}_N({M})$.
\end{lemma}
\bsp
A configuration $m_0$ of $M^{\nfinset}$ induces the following partition $\CP(m_0)$ of $${\nfinset} =m_0^{-1}(\infty) \coprod \coprod_{x \in \check{M} \cap m_0({\nfinset})}m_0^{-1}(x).$$

Pick disjoint neighborhoods $V_x$ in $M$ of the points $x$ of $m_0({\nfinset})$ that are furthermore in $\check{M}$ for $x$ in $\check{M}$ and that are identified with balls of $\RR^3$ by $C^{\infty}$-charts. Consider the neighborhood $\prod_{x \in m_0({\nfinset})}V_x^{m_0^{-1}(x)}$ of $m_0$ in $M^{\nfinset}$.
The first blow-ups that transformed this neighborhood are 
\begin{itemize}
 \item the blow-up of $E_{m_0^{-1}(\infty)}$ if $m_0^{-1}(\infty) \neq \emptyset$, which changed (a smaller neighborhood of $\infty^{m_0^{-1}(\infty)}$ in) $V_{\infty}^{m_0^{-1}(\infty)}$ to $[0,\varepsilon_{\infty}[ \times S^{3 \sharp m_0^{-1}(\infty) -1}$,
\item and the blow-ups of the $\Delta_{m_0^{-1}(x)}$, for the $x \in \check{M}$ such that $\sharp m_0^{-1}(x)\geq 2$, which changed (a smaller neighborhood of $x^{m_0^{-1}(x)}$ in) $V_{x}^{m_0^{-1}(x)}$ to $[0,\varepsilon_{x}[ \times F(U_{x}^{m_0^{-1}(x)})$, where
$U_x \subset V_x$ and $F(U_{x}^{m_0^{-1}(x)})$ fibers over $U_x$, and the fiber over $y \in U_x$ is ${S}_{m_0^{-1}(x)}(T_{y}M)$.
\end{itemize}
When considering how the next blow-ups affect the preimage of a neighborhood of $m_0$, we can
restrict to our new factors.

First consider a factor $[0,\varepsilon_{x}[ \times F(U_{x}^{m_0^{-1}(x)})$.
Picking $i \in m_0^{-1}(x)$ and fixing a Riemannian structure on $TU_x$ identifies ${S}_{m_0^{-1}(x)}(T_{y}M)$ with the space of maps $c \colon m_0^{-1}(x)\rightarrow T_yM$ such that
$c(i)=0$ and $\sum_{j \in m_0^{-1}(x)}\norm{c(j)}^2=1$. Then $(\lambda,c)$ is identified with $y + \lambda c$ in $V_x^{m_0^{-1}(x)}$ (where $V_x$ is identified with an open subset of $\RR^3$), for $\lambda \neq 0$.
Now, $[0,\varepsilon_{x}[ \times F(U_{x}^{m_0^{-1}(x)})$ must be blown-up along its intersections
with the preimage closures of the $\Delta_I$ such that $\sharp I \geq 2$, $I \subset m_0^{-1}(x)$ and $I$ is maximal. These intersections respect the product structure by $[0,\varepsilon_{x}[$ and the fibration over $U_x$ so that we only need to understand the blow-ups of the intersections of the $\Delta_I$ with a fiber of $F(U_{x}^{m_0^{-1}(x)})$. These are nothing but configurations in a ball of $\RR^3$, and we can iterate our process.

Now consider the possible factor $[0,\varepsilon_{\infty}[ \times S^{3 \sharp m_0^{-1}(\infty) -1}$ and blow up its intersections with the   preimage closures of the $E_J$ for $J \subset m_0^{-1}(\infty)$ maximal and with the preimage closures of the  $\Delta_I$ with $I \subset m_0^{-1}(\infty)$ in an order compatible with the algorithm. Here, $S^{3 \sharp m_0^{-1}(\infty) -1}$ is the unit sphere of $(\RR^3_{\infty})^{m_0^{-1}(\infty)}$.
A point $d\in (\RR^3_{\infty})^{m_0^{-1}(\infty)}$ is in the preimage closure of $E_J$ under the previous blow-up if $d(J)=0$. In particular, the $E_J$ and the $\Delta_I$ again read as products by $[0,\varepsilon_{\infty}[$, and we study what happens near a given $d$ of $S^{3 \sharp m_0^{-1}(\infty) -1}$. For such a $d$, we proceed as before if $d^{-1}(0) = \emptyset$. Otherwise the factor of $d^{-1}(0)$ must be treated differently, namely by blowing up $0^{d^{-1}(0)}$ in $S^{3 \sharp m_0^{-1}(\infty) -1}$. Then iterate.

This produces a compact manifold $C_{\nfinset}(M)$ with boundary and ridges, which is finally independent of the order of the blow-ups (when this order is compatible with the algorithm), since it is locally independent. The interior of $C_{\nfinset}(M)$ is $\check{C}_{\nfinset}(\check{M})$. Since the blow-ups separate all the pairs of points at some scale, $p_e$ naturally extends there. The introduced local coordinates show that the extension is smooth. See \cite[Section 3]{lesconst} for more details.
\eop

\begin{lemma}
\label{lemcompconfL}
The closure of $\check{C}(L;\Gamma)$ in $C_{V(\Gamma)}(M)$ is a smooth compact submanifold of $C_{V(\Gamma)}(M)$, which is denoted by ${C}(L;\Gamma)$.
\end{lemma}
\bp Exercise. \eop

Proposition~\ref{propcompext} is a consequence of Lemmas~\ref{lemcompconf} and \ref{lemcompconfL}. \eop

\subsection{Straight links}
\label{subsecstraight}

A one-chain $c$ of $S^2$ is {\em algebraically trivial\/} if for any two points $x$ and $y$ outside its support, the algebraic intersection of an arc from $x$ to $y$ transverse to $c$ with $c$ is zero, or equivalently if the integral of any one form of $S^2$ along $c$ is zero.

Let $(\check{M},\tau)$ be an asymptotic rational homology $\RR^3$.
Say that $K_j$ is {\em straight\/} with respect to $\tau$ if the curve $p_{\tau}(U^+K_j)$ of $S^2$ is algebraically trivial (recall the notation from Proposition~\ref{propprojbord} and Subsection~\ref{subconfint}). A link is {\em straight\/} with respect to $\tau$ if all its components are. If $K_j$ is straight, then $p_{\tau}(\partial C(K_j;\onechordsmallj))$ is algebraically trivial.

\begin{lemma}
Recall $C_j=C(K_j;\onechordsmallj)$, $C_j \subset C_2(M)$.

If $p_{\tau}(\partial C_j)$ is algebraically trivial, then for any propagating chain $\prop$ of $(C_2(M),\tau)$ transverse to $C_j$ and for any propagating form $\omega_p$ of $(C_2(M),\tau)$,
$$\int_{C_j}\omega_p=\langle C_j, \prop \rangle_{C_2(M)}=I_{\theta}(K_j,\tau)$$
where $I_{\theta}(K_j,\tau)$ is defined in Lemma~\ref{lemdefItheta}.
In particular, $I_{\theta}(K_j,\tau) \in \QQ$ and $I_{\theta}(K_j,\tau) \in \ZZ$ when $M$ is an integer homology $3$--sphere.
\end{lemma}
\bp Exercise. Recall Lemmas~\ref{lemetactwo} and \ref{lemvaritheta}.
\eop

\begin{proposition}
Let $\check{M}$ be an asymptotically standard $\QQ$-homology $\RR^3$.
For any parallel $K_{\parallel}$ of a knot $K$ in $\check{M}$, there exists an asymptotically standard parallelization $\tilde{\tau}$ homotopic to $\tau$, such that $K$ is straight with respect to $\tilde{\tau}$, and
$I_{\theta}(K_j,\tilde{\tau})=lk(K,K_{\parallel})$ or $I_{\theta}(K_j,\tilde{\tau})=lk(K,K_{\parallel})+1$.

For any embedding $K\colon S^1 \rightarrow \check{M}$ that is straight with respect to $\tau$, $I_{\theta}(K,\tau)$ is the linking number of $K$ and a parallel of $K$.

\end{proposition}
\bsp
For any knot embedding $K$, there is an asymptotically standard parallelization $\tilde{\tau}$ homotopic to $\tau$ such that $p_{\tilde{\tau}}(U^+K)$ is one point. Thus $K$ is straight with respect to $(M,\tilde{\tau})$.
Then $\tilde{\tau}$ induces a parallelization of $K$, and $I_{\theta}(K,\tilde{\tau})$ is the linking number of $K$ with the parallel induced by $\tilde{\tau}$. (Exercise).

In general, for two homotopic asymptotically standard parallelizations $\tau$ and $\tilde{\tau}$ such that $K$ is straight with respect to $\tau$ and $\tilde{\tau}$, $I_{\theta}(K,\tau)-I_{\theta}(K,\tilde{\tau})$ is an even integer (exercise) so that $I_{\theta}(K,\tau)$ is always the linking number of $K$ with a parallel of $K$.

In $\RR^3$ equipped with $\taust$, any link is represented by an embedding $L$ that sits in a horizontal plane except when it crosses under, so that the non-horizontal arcs crossing under are in vertical planes. Then the non-horizontal arcs have an algebraically trivial contribution to $p_{\tau}(U^+K_j)$, while the horizontal contribution  can be changed by adding kinks $\pkink$ or $\mkink$ so that $L$ is straight with respect to $\taust$.
In this case $I_{\theta}(K_j,\taust)$ is the {\em writhe\/} of $K_j$, which is the number of positive self-crossings of $K_j$ minus the number of negative self-crossings of $K_j$. In particular, up to isotopy of $L$, $I_{\theta}(K_j,\taust)$ can be assumed to be $\pm 1$ (Exercise).

Similarly, for any number $\iota$ that is congruent mod $2\ZZ$ to $I_{\theta}(K,\tau)$ there exists an embedding $K^{\prime}$ isotopic to $K$ and straight such that $I_{\theta}(K^{\prime},\tau)=\iota$ (Exercise).
\eop

\subsection{Rationality of \texorpdfstring{$Z$}{Z}}
\label{subrat}

Let us state another version of Theorem~\ref{thmmain} using straight links instead of homogeneous propagating forms. Recall $\beta_{\Gamma}=\frac{(3n-\sharp E(\Gamma))!}{(3n)!2^{\sharp E(\Gamma)}}.$

\begin{theorem}
\label{thmmainstraight}
Let $(\check{M},\tau)$ be an asymptotic rational homology $\RR^3$.
Let $L \colon \coprod_{j=1}^kS^1_j \hookrightarrow \check{M}$ be a straight embedding with respect to $\tau$.
For any $i\in \underline{3n}$, let $\omega(i)$ be a propagating form of $(C_2(M),\tau)$.
Set $$Z^s_n(L,\check{M},\tau)=\sum_{\Gamma \in \CD^e_n(C)}\beta_{\Gamma}I(\Gamma,(\omega(i))_{i \in \underline{3n}})[\Gamma] \in \CA_n^t(\coprod_{j=1}^kS^1_j).$$
Then $Z^s_n(L,\check{M},\tau)$ is independent of the chosen
$\omega(i)$. In particular, with the notation of Theorem~\ref{thmmain}, $$Z^s_n(L,\check{M},\tau)=Z_n(L,\check{M},\tau).$$
\end{theorem}

This version of Theorem~\ref{thmmain} allows us to replace the configuration space integrals by algebraic intersections in configuration spaces, and thus to prove the rationality of $Z$ for straight links as follows.

For any $i \in \underline{3n}$, let $\prop(i)$ be a propagating chain of $(C_2(M),\tau)$.
Say that a family $(\prop(i))_{i \in \underline{3n}}$ is {\em in general $3n$ position\/} with respect to $L$ if for any $\Gamma \in \CD^e_n(C)$, the $p_e^{-1}(\prop(j_E(e)))$ are pairwise transverse chains in $C(L;\Gamma)$.
In this case, define $I(\Gamma,o(\Gamma),(\prop(i))_{i \in \underline{3n}})$ as the algebraic intersection in $(C(L;\Gamma),o(\Gamma))$ of the codimension $2$ rational chains $p_e^{-1}(\prop(j_E(e)))$. If the $\omega(i)$ are propagating forms of $(C_2(M),\tau)$ Poincar\'e dual to the $\prop(i)$ and supported in sufficiently small neighborhoods of the $\prop(i)$, then 
$$I(\Gamma,o(\Gamma),(\prop(i))_{i \in \underline{3n}})=I(\Gamma,o(\Gamma),(\omega(i))_{i \in \underline{3n}})$$ for any $\Gamma \in \CD^e_n(C)$,
and $I(\Gamma,o(\Gamma),(\omega(i))_{i \in \underline{3n}})$ is rational, in this case.

\subsection{On the anomalies}
\label{subanom}

The constants $\alpha=(\alpha_n)_{n\in \NN}$ and $\ansothree=(\ansothree_n)_{n\in \NN}$ of Theorem~\ref{thmmain} are called {\em anomalies.\/}
The anomaly $\ansothree$ is the opposite of the constant $\xi$ defined in \cite[Section 1.6]{lesconst}, $\ansothree_{2n}=0$ for any integer $n$, and 
$\ansothree_1=\frac1{12}[\tata]$ according to \cite[Proposition 2.45]{lesconst}. The computation of $\ansothree_1$ can also be deduced from Corollary~\ref{corThetap}.

We define $\alpha$ below.
Let $v \in S^2$. Let $\linearmapanv$ denote the linear map
$$\begin{array}{llll}\linearmapanv: &\RR &\longrightarrow& \RR^3\\
& 1&\mapsto&v.\end{array}$$
Let $\Gamma$ be a numbered Jacobi diagram on $\RR$.
Define $\check{C}(\linearmapanv;\Gamma)$ like in Subsection~\ref{subdefconfspace} where the line $\linearmapanv$ of $\RR^3$ replaces the link $L$ of $\check{M}$.
Let $\check{Q}(v;\Gamma)$ be the quotient of $\check{C}(\linearmapanv;\Gamma)$ by the translations parallel to $\linearmapanv$ and by the dilations. Then the map $p_{e,S^2}$ associated to an edge $e$ of $\Gamma$ maps a configuration to the direction of the vector from its origin to its end in $S^2$. It factors through
$\check{Q}(v;\Gamma)$, which has two dimensions less. Now, define $\check{Q}(\Gamma)$ as the total space
of the fibration over $S^2$ whose fiber over $v$ is $\check{Q}(v;\Gamma)$. The configuration space $\check{Q}(\Gamma)$ carries a natural smooth structure, it can be compactified as before,  and it can be oriented as follows, when a vertex-orientation $o(\Gamma)$ is given. Orient $\check{C}(\linearmapanv;\Gamma)$ as before, orient $\check{Q}(v;\Gamma)$ so that $\check{C}(\linearmapanv;\Gamma)$ is locally homeomorphic to the oriented product
(translation vector $z$ in $\RR v$, ratio of homothety $\lambda \in ]0,\infty[$) $\times \check{Q}(v;\Gamma)$ and orient $\check{Q}(\Gamma)$ with the $(\mbox{base} (=S^2) \oplus \mbox{fiber})$ convention. (This can be summarized by saying that the $S^2$-coordinates replace $(z,\lambda)$.)

\begin{proposition}
\label{propanom}
For $i \in \underline{3n}$, let $\omega(i,S^2)$ be a two-form of $S^2$ such that $\int_{S^2}\omega(i,S^2)=1$. Define 
$$I(\Gamma,o(\Gamma),\omega(i,S^2))=\int_{\check{Q}(\Gamma)}\bigwedge_{e \in E(\Gamma)}p_{e,S^2}^{\ast}(\omega(j_E(e),S^2)).$$
Let $\CD^c_n(\RR)$ \index{N}{Dcn@$\CD^c_n(\RR)$} denote the set of connected numbered diagrams on $\RR$ with at least one univalent vertex, without looped edges.
Set
$$2\alpha_n=\sum_{\Gamma \in \CD^c_n(\RR)}\frac{(3n-\sharp E(\Gamma))!}{(3n)!2^{\sharp E(\Gamma)}}I(\Gamma,o(\Gamma),\omega(i,S^2))[\Gamma,o(\Gamma)] \;\; \in {\cal A}(\RR).$$
Then $\alpha_n$ does not depend on the chosen $\omega(i,S^2)$, $\alpha_1= \frac12\left[ \onechordR \right]$ and $\alpha_{2k}=0$ for all $k \in \NN$.
The series $\alpha=\sum_{n \in \NN}\alpha_n$ is called the {\em Bott and Taubes anomaly\/}.
\end{proposition}
\bp The independence of the choices of the $\omega(i,S^2)$ will be a consequence of Lemma~\ref{leminvone} below.
Let us prove that $\alpha_{2k}=0$ for all $k \in \NN$.
Let $\Gamma$ be a numbered graph and let $\overline{\Gamma}$ be obtained from $\Gamma$ by reversing the orientations of the $(\sharp E)$ edges of $\Gamma$.
Consider the map $r$ from $\check{Q}(\overline{\Gamma})$ to $\check{Q}({\Gamma})$ that composes a configuration by the multiplication by $(-1)$ in $\RR^3$. It sends a configuration over $v \in S^2$ to a configuration over $(-v)$, and it is therefore a fibered map over the orientation-reversing antipode of $S^2$. Equip $\Gamma$ and $\overline{\Gamma}$ with the same vertex-orientation. Then our map $r$ is orientation-preserving if and only if $\sharp T(\Gamma) +1 +\sharp E(\Gamma)$ is even. Furthermore for all the edges $e$ of $\overline{\Gamma}$, $p_{e,S^2} \circ r=p_{e,S^2}$, then since $\sharp E =n+\sharp T$, $$I(\overline{\Gamma},o(\Gamma),\omega(i,S^2))=(-1)^{n+1}I({\Gamma},o(\Gamma),\omega(i,S^2)).$$
\eop

It is known that $\alpha_3=0$ and $\alpha_5=0$ \cite{poirier}.
Furthermore, according to \cite{lesunikon}, $\alpha_{2n+1}$ is a combination of diagrams with two univalent vertices, and $Z^u(S^3,L)$ is obtained from the Kontsevich integral by inserting  $d$
times the plain part of $2\alpha$ on each degree $d$ connected component of a diagram.

\subsection{The dependence on the forms in the invariance proofs}

The variation of $I(\Gamma,o(\Gamma),(\omega(j))_{j \in \underline{3n}})$ when some
$\omega(i=j_E(f \in E(\Gamma)))$ is changed to $\omega(i) +d\eta$ for a one-form $\eta$ on $C_2(M)$ reads
$$\int_{(C(L;\Gamma),o(\Gamma))} \left( p_f^{\ast}(d \eta) \wedge \bigwedge_{e \in (E(\Gamma)\setminus \{f\})}p_e^{\ast}(\omega(j_E(e))) \right).$$
According to the Stokes theorem, it reads $\int_{\partial (C(L;\Gamma),o(\Gamma))}\left( p_f^{\ast}( \eta) \wedge \bigwedge_{e \in (E(\Gamma)\setminus \{f\})}p_e^{\ast}(\omega(j_E(e))) \right)$ 
where the integral along $\partial (C(L;\Gamma),o(\Gamma))$ is actually the integral along the codimension one faces of $C(L;\Gamma)$, which are considered as open. Such a codimension one face only involves one blow-up. 

For any non-empty subset $B$ of $V(\Gamma)$, the codimension one face associated to the blow-up of $E_B$ in $M^{V(\Gamma)}$ is denoted by
$F(\Gamma,\infty,B)$, it lies in the preimage of $\infty^B \times \check{M}^{V(\Gamma) \setminus B}$, in $C(L;\Gamma)$.

The other codimension one faces are associated to the blow-ups of the $\Delta_B$ in $M^{V(\Gamma)}$, for subsets $B$ of $V(\Gamma)$ of cardinality at least $2$. The face of $C(L;\Gamma)$ associated to $\Delta_B$ is denoted by $F(\Gamma,B)$.
Let $b \in B$. Assume that $b \in U(\Gamma)$ if $U(\Gamma) \cap B \neq \emptyset$. The image of $F(\Gamma,B)$ in $M^{V(\Gamma)}$
is in the set of maps $m$ of $\Delta_B$ that define an injection from $(V(\Gamma) \setminus B ) \cup\{b \in B\}$ to $\check{M}$, which factors through an injection isotopic to the restriction of $i_{\Gamma}$ on $U(\Gamma) \cap \left((V(\Gamma) \setminus B ) \cup\{b\} \right)$. This set of maps
$\check{C}_{(V(\Gamma) \setminus B) \cup\{b \}}(\check{M},i_{\Gamma})$ is a submanifold of $\check{C}_{(V(\Gamma) \setminus B) \cup\{b \}}(\check{M})$.
Thus, $F(\Gamma,B)$ is a bundle over $\check{C}_{(V(\Gamma) \setminus B ) \cup\{b \}}(\check{M},i_{\Gamma})$.

When $B$ has no univalent vertices, the fiber over a map $m$ is the space
$\check{S}_B(T_{m(b)})$ of injective maps from $B$ to $T_{m(b)}$
up to translations and dilations.

When $B$ contains univalent vertices of a component $K_j$, the fiber over $m$ is the submanifold $\check{S}_B(T_{m(b)}M,\Gamma)$ of $\check{S}_B(T_{m(b)}M)$, made of the configurations that map the univalent vertices of $B$ to a line of $T_{m(b)}M$ directed by $U^+K_j$ at $m(b)$, in an order prescribed by $\Gamma$.
If $B$ does not contain  all the univalent vertices of $\Gamma$ on $S^1_j$, this order is unique.
Otherwise, $F(\Gamma,B)$ has $\sharp (B \cap U(\Gamma))$ connected components corresponding to the total orders that induce the cyclic order of $B \cap U(\Gamma)$.

When $B$ is a subset of the set of vertices $V(\Gamma)$ of a numbered graph $\Gamma$,
$E(\Gamma_B)$ denotes the set of edges of $\Gamma$ between two elements of $B$ (edges of $\Gamma$ are plain), and $\Gamma_B$ is the subgraph of $\Gamma$ made of the vertices of $B$ and the edges of $E(\Gamma_B)$.

\begin{lemma}
\label{leminvone}
Let $(\check{M},\tau)$ be an asymptotic rational homology $\RR^3$. Let $C=\coprod_{j=1}^kS^1_j$.

For $i \in \underline{3n}$, let $\omega(i)$ be a closed $2$-form on $[0,1] \times C_2(M)$ whose restriction to $\{t\} \times C_2(M)$ is denoted by $\omega(i,t)$, for any $t \in [0,1]$.
Assume that for $t \in [0,1]$, $\omega(i,t)$ 
restricts to $(\partial C_2(M) \setminus UB_M)$ as $p_{\tau}^{\ast}(\omega(i,t)(S^2))$, for some two-form $\omega(i,t)(S^2)$ of $S^2$ such that $\int_{S^2}\omega(i,t)(S^2)=1$.
Set $$Z_n(t)=\sum_{\Gamma \in \CD^e_n(C)}\beta_{\Gamma}I(\Gamma,(\omega(i,t))_{i \in \underline{3n}})[\Gamma] \in \CA_n^t(\coprod_{j=1}^kS^1_j).$$
Then $$Z_n(1)-Z_n(0)=\sum_{\left\{\begin{array}{c}(\Gamma,B);\Gamma \in \CD^e_n(C), B \subset V(\Gamma), \sharp B \geq 2;\\\Gamma_B \;\mbox{\small is a connected component of }\;\Gamma\end{array}\right\}}I(\Gamma,B)$$
where
$$I(\Gamma,B)= \beta_{\Gamma}\int_{[0,1]\times F(\Gamma,B)}\bigwedge_{e \in E(\Gamma)}p_e^{\ast}(\omega(j_E(e)))[\Gamma].$$
Under the assumptions of Theorem~\ref{thmmain} (where the $\omega(i)$ are homogeneous) or Theorem~\ref{thmmainstraight} (where $L$ is straight with respect to $\tau$), when $(M,L,\tau)$ is fixed,
$Z_n(L,\check{M},\tau)$ is independent of the chosen $\omega(i)$.

In particular, when $k=0$, $Z(\check{M},\tau)$ coincides with the Kontsevich configuration space integral invariant described in \cite{lesconst}.

Furthermore, the $\alpha_n$ of Proposition~\ref{propanom} are also independent of the forms $\omega(i,S^2)$.

\end{lemma}
\bsp According to the Stokes theorem, for any $\Gamma \in \CD^e_n(C)$,
$$I(\Gamma,(\omega(i,1))_{i \in \underline{3n}})-I(\Gamma,(\omega(i,0))_{i \in \underline{3n}})=\sum_F\int_{[0,1]\times F} \bigwedge_{e \in E(\Gamma)}p_e^{\ast}(\omega(j_E(e)))$$
where the sum runs over the codimension one faces $F$ of $C(L;\Gamma)$.
Below, we sketch the proof that the only contributing faces are the faces $F(\Gamma,B)$ such that $\sharp B \geq 2$ and $\Gamma_B$ is a connected component of $\Gamma$, or equivalently, that the other faces do not contribute.

Like in \cite[Lemma 2.17]{lesconst} faces $F(\Gamma,\infty,B)$ do not contribute.
When the product of all the $p_e$ factors through a quotient of $[0,1]\times F(\Gamma,B)$ of smaller dimension, the face $F(\Gamma,B)$ does not contribute. This allows us to get rid of 
\begin{itemize}
 \item the faces $F(\Gamma,B)$ such that $B$ is not a pair of univalent vertices of $\Gamma$, and $\Gamma_B$ is not connected  (see \cite[Lemma 2.18]{lesconst}),
\item the faces $F(\Gamma,B)$ such that $\sharp B \geq 3$ where $\Gamma_B$ has a univalent vertex that was trivalent in $\Gamma$ (see \cite[Lemma 2.19]{lesconst}).
\end{itemize}
We also have faces that cancel each other, for graphs that are identical outside their $\Gamma_B$ part.
\begin{itemize}
 \item The faces $F(\Gamma,B)$ (that are not already listed) such that $\Gamma_B$ has at least a bivalent vertex cancel (mostly by pairs) by the parallelogram identification (see \cite[Lemma 2.20]{lesconst}).
\item The faces $F(\Gamma,B)$ where $\Gamma_B$ is an edge between two trivalent vertices cancel by triples, thanks to the Jacobi (or IHX) relation (see \cite[Lemma 2.21]{lesconst}).
\item Similarly, two faces where $B$ is made of two (necessarily consecutive in $C$) univalent vertices of $\Gamma$ cancel $(3n-\sharp E(\Gamma))$ faces $F(\Gamma^{\prime},B^{\prime})$ where $\Gamma^{\prime}_{B^{\prime}}$ is an edge between a univalent vertex of $\Gamma$ and a trivalent vertex of $\Gamma$, thanks to the STU relation.
\end{itemize}

Thus, we are left with the faces $F(\Gamma,B)$ such that $\Gamma_B$ is a (plain) connected component of $\Gamma$, and we get the wanted formula for $(Z_n(1)-Z_n(0))$.

In the anomaly case, the same analysis of faces leaves no contributing faces, so that the $\alpha_n$ are independent of the forms $\omega(i,S^2)$ in Proposition~\ref{propanom}.

Back to the behaviour of $Z(L,\check{M},\tau)$ under the assumptions of Theorem~\ref{thmmain} or Theorem~\ref{thmmainstraight}, assume that $(M,L,\tau)$ is fixed and apply the formula of the lemma to compute the variation
of $Z_n(L,\check{M},\tau)$ when some propagating chain $\omega(i,0)$ of $(C_2(M),\tau)$ is changed to some other propagating chain $\omega(i,1)=\omega(i,0)+d\eta$. According to Lemma~\ref{lemetactwo}, under our assumptions, $\eta$ can be chosen so that $\eta=p_{\tau}^{\ast}(\eta_{S^2})$ on $\partial C_2(M)$ and $\eta_{S^2}=0$ if $\omega(i,0)$ and $\omega(i,1)$ are homogeneous.
Define $\omega(i)=\omega(i,0)+d(t \eta)$ on $[0,1]\times C_2(M)$ ($t\in [0,1]$), and extend the other $\omega(j)$ trivially. 

Then $(Z_n(1)-Z_n(0))$ vanishes if $\omega(i,0)$ and $\omega(i,1)$ are homogeneous, as all the involved $I(\Gamma,B)$ do, so that $Z_n(L,\check{M},\tau)$ is independent from the chosen homogeneous propagating forms $\omega(i)$ of $C_2(M,\tau)$ in Theorem~\ref{thmmain}. Now, assume that $L$ is straight.

When $i \notin j_E(E(\Gamma))$, the integrand of $I(\Gamma,B)$ factors through the natural projection of $[0,1]\times F(\Gamma,B)$ onto $F(\Gamma,B)$, and $I(\Gamma,B)=0$, consequently. 
Assume $i=j_E(e_i \in E(\Gamma))$, then
$$I(\Gamma,B)=\beta_{\Gamma}\int_{[0,1]\times F(\Gamma,B)}p_{e_i}^{\ast}(d(t \eta )) \wedge \bigwedge_{e \in E(\Gamma)\setminus e_i}p_e^{\ast}(\omega(j_E(e)).$$
The form $\bigwedge_{e \in E(\Gamma_B)}p_e^{\ast}(\omega(j_E(e))$ pulls back through $[0,1]\times F(\Gamma_B,B)$, and through $F(\Gamma_B,B)$ when $e_i \notin E(\Gamma_B)$, so that, for dimension reasons, $I(\Gamma,B)$ vanishes unless $e_i \in E(\Gamma_B)$. Therefore, we assume $e_i \in E(\Gamma_B)$.

When $B$ contains no univalent vertices, $I(\Gamma,B)$ 
factors through
$$\int_{[0,1] \times \cup_{m(b)\in \check{M}}\check{S}_B(T_{m(b)}M)} p_{e_i}^{\ast}(d(t \eta ))  \wedge \bigwedge_{e \in E(\Gamma_B)\setminus e_i}p_e^{\ast}(\omega(j_E(e))).$$
Here the parallelization $\tau$ identifies the bundle $\cup_{m(b)\in \check{M}}\check{S}_B(T_{m(b)}M)$ with $\check{M} \times \check{S}_B(\RR^3)$, and the integrand factors through the projection of $[0,1]\times \check{M} \times \check{S}_B(\RR^3)$ onto $ [0,1]\times \check{S}_B(\RR^3)$ whose dimension is smaller (by $3$). In particular, $I(\Gamma,B)=0$ in this case,  the independence of the choice of the $\omega(i)$ is proved when $k=0$ (when the link is empty), and $Z(\check{M},\tau)$ coincides with the Kontsevich configuration space integral invariant described in \cite{lesconst}.

Let us now study the sum of the $I(\Gamma,B)$, where $(\Gamma \setminus \Gamma_B)$ is a fixed labeled graph and $\Gamma_B$ is a fixed numbered connected diagram with at least one univalent vertex
on $S^1_j$.

This sum factors through $\int_{[0,1] \times \cup_{m(b)\in K_j}\check{S}_B(T_{m(b)}M,\Gamma)}p_{e_i}^{\ast}(d(t \eta )) \wedge \bigwedge_{e \in E(\Gamma_B) \setminus e_i}p_e^{\ast}(\omega(j_E(e)))$.

At a collapse, the univalent vertices of $\Gamma_B$ are equipped with a linear order, which makes $\Gamma_B$ a numbered graph $\tilde{\Gamma}_B$ on $\RR$.
The corresponding connected component of $[0,1] \times \cup_{m(b)\in K_j}\check{S}_B(T_{m(b)}M,\Gamma)$ reads $[0,1] \times \cup_{x \in U^+K_j} \check{Q}(p_{\tau}(x);\tilde{\Gamma}_B)$ ($\check{Q}(v;\tilde{\Gamma}_B)$ was defined in Subsection~\ref{subanom}). This allows us to see the contribution of such a connected component as the integral of a one-form (defined by partial integrations) over $p_{\tau}(U^+K_j)$. Such an integral is zero when $K_j$ is straight.

\eop

Now, Theorem~\ref{thmmainstraight} is a corollary of Theorem~\ref{thmmain} (which is not yet completely proved).

\subsection{The dependence on the parallelizations in the invariance proofs}

Recall that $\CA^t_n(C)$ splits according to the number of connected components without univalent vertices of the graphs. Then it is easy to observe that
$$Z(L,\check{M},\tau)=\sum_{n \in \NN}Z_n(L,\check{M},\tau)=Z^u(L,\check{M},\tau)Z(M;\tau)$$
where $Z^u$ is obtained from $Z$ by sending the graphs with components that have no univalent vertices to $0$, and $Z(M;\tau)=Z(\emptyset,\check{M},\tau)$.
According to \cite[Theorem 1.9]{lesconst}, $Z(M)=Z(M;\tau)\exp(-\frac14 p_1(\tau) \ansothree)$ is a topological invariant of $M$. Here, we will now focus on $Z^u(L,\check{M},\tau)$, and define it with a given homogeneous propagating form, $\omega=\omega(i)$ for all i, so that $Z^u(L,\check{M},\tau)$ is an invariant of the diffeomorphism class of $(L,\check{M},\tau)$. We study its variation under a continuous deformation of $\tau$ and we prove the following lemma.

\begin{lemma}
\label{lemvartauanom}
Let $(\tau(t))_{t \in [0,1]}$ define a smooth homotopy of asymptotically standard parallelizations of $\check{M}$.
$$\frac{\partial}{\partial t}Z^u(L,\check{M},\tau(t))=(\sum_{j=1}^k\frac{\partial}{\partial t}I_{\theta}(K_j,\tau(t))\alpha \sharp_j)Z^u(L,\check{M},\tau(t)).$$
\end{lemma}
\bp 
Set $Z_n(t)=Z_n^u(L,\check{M},\tau(t))$, observe that $Z_n$ (which is valued in a finite-dimensional vector space) is differentiable thanks to the expression of $Z_n(t)-Z_n(0)$ in Lemma~\ref{leminvone} (any function $\int_{[0,t]\times C}\omega$ for a smooth compact manifold $C$ and a smooth form $\omega$ on $[0,1]\times C$ is differentiable with respect to $t$).
Now, the forms associated to edges of $\Gamma_B$ do not depend on the configuration of $ (V(\Gamma) \setminus B)$. They will be integrated along $[0,1] \times (\cup_{m(b)\in K_j}\check{S}_B(T_{m(b)}M,\Gamma_B))$, while the other ones will be integrated along $\check{C}(L;\Gamma \setminus \Gamma_B)$ at $u \in [0,1]$.

Therefore, the global variation $(Z(t)-Z(0))$ reads
$$\sum_{j=1}^k\int_0^t \left(\sum_{\Gamma_B \in \CD^c(\RR)}\beta_{\Gamma_B}\int_{c \in \cup_{m(b)\in K_j}\check{S}_B(T_{m(b)}M,\Gamma_B)} \left(\bigwedge_{e \in E(\Gamma_B)}p_e^{\ast}(\omega_{S^2})\right)(u,c)[\Gamma_B]\sharp_j\right) Z(u) du$$
where $\CD^c(\RR)=\cup_{n \in \NN} \CD^c_n(\RR).$ \index{N}{Dc@$\CD^c(\RR)$}
Define $$I(\Gamma_B,K_j)(t)=\int_{(u,c); u \in [0,t],c \in \cup_{m(b)\in K_j}\check{S}_B(T_{m(b)}M,\Gamma_B)} \bigwedge_{e \in E(\Gamma_B)}p_e^{\ast}(\omega_{S^2})(u,c)$$
so that $\frac{\partial}{\partial u}I(\Gamma_B,K_j)(u)du$ is the integral of $ \bigwedge_{e \in E(\Gamma_B)}p_e^{\ast}(\omega)$
along $u \times (\cup_{m(b)\in K_j}\check{S}_B(T_{m(b)}M,\Gamma_B))$
and 
$$Z(t)-Z(0)=\int_0^t\left(\sum_{j=1}^k\left(\sum_{\Gamma_B \in \CD^c(\RR)}\beta_{\Gamma_B} \frac{\partial}{\partial u}I(\Gamma_B,K_j)(u)[\Gamma_B]\sharp_j\right) Z(u)\right) du.$$
Therefore, $$\frac{\partial}{\partial t}Z(t)=\sum_{j=1}^k\left(\sum_{\Gamma_B \in \CD^c(\RR)}\beta_{\Gamma_B} \frac{\partial}{\partial t}I(\Gamma_B,K_j)(t)[\Gamma_B]\sharp_j\right) Z(t)$$
and we are left with the computation of $\frac{\partial}{\partial t}I(\Gamma_B,K_j)(t)$.

The restriction of $p_{\tau(.)}$ from $[0,1]\times U^+K_j$ to $S^2$ induces a map 
$$p_{a,\tau,\Gamma_B}\colon [0,1]  \times \cup_{m(b)\in K_j}\check{S}_B(T_{m(b)}M,\Gamma_B) \rightarrow \check{Q}({\Gamma}_B)$$ for any $\Gamma_B$.
$$I(\Gamma_B,K_j)(t)=\int_{\mbox{Im}(p_{a,\tau,\Gamma_B})}\bigwedge_{e \in E(\Gamma_B)}p_e^{\ast}(\omega_{S^2}).$$
Integrating $\bigwedge_{e \in E(\Gamma_B)}p_e^{\ast}(\omega_{S^2})[\Gamma_B]$ along the fiber in $\check{Q}(\Gamma_B)$ yields a two--form
on $S^2$, which is homogeneous, because everything is. Thus this form reads $2\alpha(\Gamma_B)\omega_{S^2}[\Gamma_B]$ where $\alpha(\Gamma_B) \in \RR$, and where $\sum_{\Gamma_B \in \CD^c(\RR)}\beta_{\Gamma_B}\alpha(\Gamma_B)[\Gamma_B]=\alpha$.
Therefore
$$I(\Gamma_B,K_j)(t)=2\beta_{\Gamma_B}\alpha(\Gamma_B)\int_{[0,t]\times U^+K_j}p_{\tau(.)}^{\ast}(\omega_{S^2}).$$
Since $\frac{\partial}{\partial t}\int_{[0,t]\times U^+K_j}p_{\tau(.)}^{\ast}(\omega_{S^2})=\frac12 \frac{\partial}{\partial t}I_{\theta}(K_j,\tau(t))$, we conclude easily.
\eop

Then the derivative of $$\prod_{j=1}^k\exp(-I_{\theta}(K_j,\tau(t))\alpha)\sharp_j Z^u(L,\check{M},\tau(t))$$ vanishes so that this expression does not change when $\tau$ smoothly varies.

\subsection{End of the proof of Theorem~\ref{thmmain}}

Thanks to \cite[Theorem 1.9]{lesconst}, in order to conclude the (sketch of) proof of Theorem~\ref{thmmain}, we are left with the proof that $$\prod_{j=1}^k\left(\exp(-I_{\theta}(K_j,\tau)\alpha)\sharp_j\right) Z^u(L,\check{M},\tau)$$
does not depend on the homotopy class of $\tau$.

When $\tau$ changes in a ball that does not meet the link, the forms can be changed only in the neighborhoods of the unit tangent bundle to this ball. Using Lemma~\ref{leminvone} again, the variation will be seen on faces $F(\Gamma,B)$, where $\Gamma_B$ has at least one univalent vertex, and where the forms associated to the edges of $\Gamma_B$ do not depend on the parameter in $[0,1]$ so that their product vanishes. In particular, 
$$\prod_{j=1}^k\left(\exp(-I_{\theta}(K_j,\tau)\alpha)\sharp_j\right) Z^u(L,\check{M},\tau)$$
is invariant under the natural action of $\pi_3(SO(3))$ on the homotopy classes of parallelizations.

We now examine the effect of the twist of the parallelization by a map $g \colon (B_M,1) \rightarrow (SO(3),1)$.
Without loss, assume that $p_{\tau}(U^+K_j)=v$ for some $v$ of $S^2$ and that $g$ maps $K_j$ to rotations with axis $v$.
We want to compute $Z^u(L,\check{M},\tau \circ \psi_{\RR}(g))- Z^u(L,\check{M},\tau)$.
Identify $UB_M$ with $B_M \times S^2$ via $\tau$. There exists a form $\omega$ on
$[0,1]\times B_M \times S^2$ that reads $p_{\tau}^{\ast}(\omega_{S^2})$ on $\partial ([0,1]\times B_M \times S^2) \setminus (1 \times B_M \times S^2)$) and that reads $p_{\tau \circ \psi_{\RR}(g)}^{\ast}(\omega_{S^2})$ on $ 1 \times B_M \times S^2$.
Extend this form to a form $\Omega$ on $[0,1]\times C_2(M)$, that restricts to $0 \times \partial C_2(M)$ as $p_{\tau}^{\ast}(\omega_{S^2})$, and to $1 \times \partial C_2(M)$ as $p_{\tau \circ \psi_{\RR}(g)}^{\ast}(\omega_{S^2})$, where $p_{\tau \circ \psi_{\RR}(g)} =p_{\tau} \circ \psi_{\RR}(g^{-1})$ on $B_M \times S^2$ so that $p_{\tau \circ \psi_{\RR}(g)}^{\ast}(\omega_{S^2})=\psi_{\RR}(g^{-1})^{\ast}\left(p_{\tau
}^{\ast}(\omega_{S^2})\right) $, there.
Let $\CD^{e,u}_n(C)$ denote the set of diagrams of $\CD^e_n(C)$ without components without univalent vertices.
Define $$Z_n(t)=\sum_{\Gamma \in \CD^{e,u}_n(C)}\beta_{\Gamma}I(\Gamma,(\Omega_{|t \times C_2(M)})_{i \in \underline{3n}})[\Gamma] \in \CA_n(\coprod_{j=1}^kS^1_j).$$
For $\Gamma_B \in \CD^c(\RR)$,
set $$I(\Gamma_B,K_j,\Omega)(t)=\int_{(u,c); u \in [0,t],c \in \cup_{m(b)\in K_j}\check{S}_B(T_{m(b)}M,\Gamma_B)} \bigwedge_{e \in E(\Gamma_B)}p_e^{\ast}(\Omega)[\Gamma_B].$$
Set $\beta_j(t)=\sum_{\Gamma_B \in \CD^c(\RR)}\beta_{\Gamma_B}I(\Gamma_B,K_j,\Omega)(t)$ and $\gamma_j(t)=\frac{\partial }{\partial t}\beta_j(t)$.
Thanks to Lemma~\ref{leminvone}, like in the proof of Lemma~\ref{lemvartauanom}, $Z(t)$ is differentiable, and
$Z^{\prime}(t) = (\sum_{j=1}^k \gamma_j(t) \sharp_j) Z(t)$.

By induction on the degree, it is easy to see that this equation determines $Z(t)$ as a function of the $\beta_j(t)$ and
$Z(0)$ whose degree $0$ part is $1$, and that $Z(t)=\prod_{j=1}^k\exp(\beta_j(t)) \sharp_j Z(0)$.

Extend $\Omega$ over $[0,2]\times C_2(M)$ so that its restriction to $[1,2]\times B_M \times S^2$ is obtained by applying $\left(\psi_{\RR}(g^{-1})\right)^{\ast}$ to the $\Omega$ translated, and extend all the introduced maps, then $\gamma_j(t+1)=\gamma_j(t)$ because everything is carried by $\left(\psi_{\RR}(g^{-1})\right)^{\ast}$.
In particular $\beta_j(2)=2\beta_j(1)$.

Now, $Z(2)=Z^u(L,M,\tau \circ \psi_{\RR}(g)^2)=\prod_{j=1}^k\exp((I_{\theta}(K_j,\tau \circ \psi_{\RR}(g)^2)-I_{\theta}(K_j,\tau))\alpha)\sharp_j Z^u(L,\check{M},\tau)$, since $g^2$ is homotopic to the trivial map outside a ball (see Lemma~\ref{lempreptrivun},~2).
By induction on the degree of diagrams, this shows
$\beta_j(2)=(I_{\theta}(K_j,\tau \circ \psi_{\RR}(g)^2)-I_{\theta}(K_j,\tau))\alpha$.
Conclude by observing that under our assumptions, where $I_{\theta}(K_j,\tau \circ \psi_{\RR}(g)^i)$ is the linking number of $K_j$ and its parallel induced by $\tau \circ \psi_{\RR}(g)^i$,
$$I_{\theta}(K_j,\tau \circ \psi_{\RR}(g)^2)-I_{\theta}(K_j,\tau)=2(I_{\theta}(K_j,\tau \circ \psi_{\RR}(g))-I_{\theta}(K_j,\tau)).$$
This finishes the (sketch of) proof of Theorem~\ref{thmmain} in general.

\subsection{Some open questions}
\begin{enumerate}
\item A Vassiliev invariant is {\em odd\/} if it distinguishes some knot from the same knot with the opposite orientation. Are there odd Vassiliev invariants ?
\item More generally, do Vassiliev invariants distinguish knots in $S^3$ ?
\item According to a theorem of Bar-Natan and Lawrence \cite{barnatanlawr}, the LMO invariant fails to distinguish rational homology $3$--spheres with isomorphic $H_1$, so that, according to a Moussard theorem \cite{moussardAGT}, rational finite type invariants fail to distinguish $\QQ$-spheres. Do finite type invariants distinguish $\ZZ$-spheres ?
\item Find relationships between $Z$ or other finite type invariants and Heegaard Floer homologies. See \cite{lesHC} to get propagators associated to Heegaard diagrams. Also see related work by Shimizu and Watanabe \cite{shimizu,watanabeMorse}.
\item Compare $Z$ with the LMO invariant $Z_{LMO}$.
\item Compute the anomalies $\alpha$ and $\ansothree$.
\item Find surgery formulae for $Z$.
\item 
Kricker defined a lift $\tilde{Z}^K$ of the Kontsevich integral $Z^K$ (or the LMO invariant) for null-homologous knots in $\QQ$-spheres \cite{kricker,garkri}.
The Kricker lift is valued in a space $\tilde{A}$ that is mapped to $\CA_n(S^1)$ by a map $H$, which allows one to recover $Z^K$ from $\tilde{Z}^K$.
The space $\tilde{A}$ is a space of trivalent diagrams whose edges are decorated by rational functions whose denominators divide the Alexander polynomial.
Compare the Kricker lift $\tilde{Z}^K$ with the equivariant configuration space invariant $\tilde{Z}^c$ of \cite{lesbonn}  valued in the same diagram space $\tilde{A}$. See \cite{lesuniveq} for alternative definitions and further properties of $\tilde{Z}^c$.
\item Is $Z$ obtained from $\tilde{Z}^c$ in the same way as $Z^K$ is obtained from $\tilde{Z}^K$ ?

\end{enumerate}

\newpage
\section{More on parallelizations of \texorpdfstring{$3$}{3}--manifolds and Pontrjagin classes}
\setcounter{equation}{0}
\label{secmorepar}

In order to make the definition of $\Theta$ complete, we give a detailed self-contained presentation of $p_1(\tau)$.
In this section, $M$ is a smooth oriented connected $3$-manifold with possible boundary.

\subsection{\texorpdfstring{$[(M,\partial M),(SO(3),1)]$}{[(M, boundary of M),(SO(3),1)]} is an abelian group.}

Again, see $S^3$ as $B^3/\partial B^3$ and see $B^3$ as $([0,2\pi] \times S^2)/(0 \sim \{0\} \times S^2)$.
Recall that $\rho \colon B^3 \rightarrow SO(3)$ maps $(\theta \in[0,2\pi], v \in S^2)$ to the rotation $\rho(\theta,v)$ with axis directed by $v$ and with angle $\theta$. 

Also recall that the group structure of $[(M,\partial M),(SO(3),1)]$ is induced by the multiplication of maps, using the multiplication of $SO(3)$.

Any $g \in C^0\left((M,\partial M),(SO(3),1)\right)$
 induces a map $$H_1(g;\ZZ) \colon H_1(M,\partial M;\ZZ)\longrightarrow (H_1(SO(3),1)=\ZZ/2\ZZ).$$
Since $$H_1(M,\partial M;\ZZ/2\ZZ)=H_1(M,\partial M;\ZZ)/2H_1(M,\partial M;\ZZ)=H_1(M,\partial M;\ZZ)\otimes_{\ZZ}\ZZ/2\ZZ,$$
$\mbox{Hom}(H_1(M,\partial M;\ZZ),\ZZ/2\ZZ)=\mbox{Hom}(H_1(M,\partial M;\ZZ/2\ZZ),\ZZ/2\ZZ)=H^1(M,\partial M;\ZZ/2\ZZ)$, and the image of $H_1(g;\ZZ)$ under the above isomorphisms is denoted by
$H^1(g;\ZZ/2\ZZ)$. (Formally, this $H^1(g;\ZZ/2\ZZ)$ denotes the image of the generator of $H^1(SO(3),1;\ZZ/2\ZZ)=\ZZ/2\ZZ$ under $H^1(g;\ZZ/2\ZZ)$ in $H^1(M,\partial M;\ZZ/2\ZZ)$.)

\begin{lemma}
\label{lempreptrivun}
Let $M$ be an oriented connected $3$-manifold with possible boundary. Recall that $\rho_M(B^3) \in C^0\left((M,\partial M),(SO(3),1)\right)$ is a map that coincides with $\rho$ on a ball $B^3$ embedded in $M$ and that maps the complement of $B^3$ to the unit of $SO(3)$.
\begin{enumerate}
\item Any homotopy class of a map $g$ from $(M,\partial M)$ to $(SO(3),1)$, such that $H^1(g;\ZZ/2\ZZ)$
is trivial, belongs to the subgroup $<[\rho_M(B^3)]>$ of 
$[(M,\partial M),(SO(3),1)]$
generated by $[\rho_M(B^3)]$. 
\item For any $[g] \in [(M,\partial M),(SO(3),1)]$, 
$[g]^2 \in <[\rho_M(B^3)]>.$
\item The group $[(M,\partial M),(SO(3),1)]$ is abelian.
\end{enumerate}
\end{lemma}
\bp Let $g \in C^0\left((M,\partial M),(SO(3),1)\right)$. Assume that $H^1(g;\ZZ/2\ZZ)$ is trivial.
Choose a cell decomposition of $M$ with respect to its boundary,
with only one three-cell, no zero-cell if $\partial M \neq \emptyset$, one zero-cell if $\partial M = \emptyset$,  one-cells, and two-cells. Then after a homotopy relative to $\partial M$, we may assume that $g$ maps the one-skeleton of $M$ to $1$. 
Next, since $\pi_2(SO(3)) = 0$, we may assume that $g$ maps the two-skeleton of $M$ to $1$, and therefore that $g$ maps the exterior of some $3$-ball to $1$. Now $g$ becomes a map from $B^3/\partial B^3=S^3$ to $SO(3)$, and its homotopy class is $k[\tilde{\rho}]$ in $\pi_3(SO(3))=\ZZ[\tilde{\rho}]$.
Therefore $g$ is homotopic to $\rho_M(B^3)^k$. This proves the first assertion.

Since $H^1(g^2;\ZZ/2\ZZ)=2H^1(g;\ZZ/2\ZZ)$ is trivial, the second assertion follows.

For the third assertion, first note that $[\rho_M(B^3)]$ belongs to the center of
$[(M,\partial M),(SO(3),1)]$ because it can be supported in a small ball disjoint 
from the support (preimage of $SO(3) \setminus \{1\}$) of a representative of any other element. Therefore, according to the second assertion any square will be in the center. Furthermore, since any commutator induces the trivial map on $\pi_1(M)$, any commutator is in $<[\rho_M(B^3)]>$.
In particular, if $f$ and $g$ are elements of $[(M,\partial M),(SO(3),1)]$, 
$$(gf)^2=(fg)^2=(f^{-1}f^2g^2f)(f^{-1}g^{-1}fg)$$ 
where the first factor equals $f^2g^2=g^2f^2$. Exchanging $f$ and $g$ yields
$f^{-1}g^{-1}fg=g^{-1}f^{-1}gf$. Then the commutator, which is a power of $[\rho_M(B^3)]$, has a vanishing square, and thus a vanishing degree. Then it must be trivial.
\eop

\subsection{Any oriented \texorpdfstring{$3$}{3}--manifold is parallelizable.}
\label{subproofpar}

In this subsection, we prove the following standard theorem. The spirit of our proof is the same as the Kirby proof in \cite[p.46]{Kir}. But instead of assuming familiarity with the obstruction theory described by Steenrod in \cite[Part III]{St}, we use this proof as an introduction to this theory.

\begin{theorem}[Stiefel]
\label{thmpar}
Any oriented $3$-manifold is parallelizable.
\end{theorem}

\begin{lemma}
 The restriction of the tangent bundle $TM$ to an oriented $3$-manifold $M$ to any closed (non-necessarily orientable) surface $S$ immersed in $M$ is trivializable.
\end{lemma}
\bp Let us first prove that this bundle is independent of the immersion. It is the direct sum of the tangent bundle to the surface and of its normal
one-dimensional bundle. This normal bundle is trivial when $S$ is orientable, and its unit bundle is the $2$-fold orientation cover of the surface, otherwise. (The orientation cover of $S$ is its $2$-fold orientable cover, which is trivial over annuli embedded in the surface). Then since any surface $S$ can be immersed in $\RR^3$, the restriction $TM_{|S}$ is the pull-back of the trivial bundle of $\RR^3$ by such an immersion, and it is trivial.
\eop

Then using Stiefel-Whitney classes, the proof of Theorem~\ref{thmpar} quickly goes as follows.
Let $M$ be an orientable smooth $3$-manifold, equipped with a smooth triangulation. (A theorem of Whitehead proved in the Munkres book \cite{Mu} ensures the existence of such a triangulation.)
By definition, the {\em first Stiefel-Whitney class\/} $w_1(TM) \in H^1(M;\ZZ/2\ZZ=\pi_0(GL(\RR^3)))$ seen as a map from
$\pi_1(M)$ to $\ZZ/2\ZZ$ maps the class of a loop $c$ embedded in $M$ to $0$ if $TM_{|c}$ is orientable and to $1$ otherwise.
It is the obstruction to the existence of a trivialization of $TM$ over the one-skeleton of $M$.
Since $M$ is orientable, the first Stiefel-Whitney class $w_1(TM)$ vanishes and $TM$ can be trivialized over the one-skeleton of $M$.
The {\em second Stiefel-Whitney class\/} $w_2(TM) \in H^2(M;\ZZ/2\ZZ=\pi_1(GL^+(\RR^3)))$ seen as a map from $H_2(M;\ZZ/2\ZZ)$ to $\ZZ/2\ZZ$ maps the class of a connected closed surface $S$ to $0$ if $TM_{|S}$ is trivializable and to $1$ otherwise.
The second Stiefel-Whitney class $w_2(TM)$ is the obstruction to the existence of a trivialization of $TM$ over the two-skeleton of $M$, when $w_1(TM)=0$.
According to the above lemma, $w_2(TM)=0$, and $TM$ can be trivialized over the two-skeleton of $M$.
Then since $\pi_2(GL^+(\RR^3))=0$, any parallelization over the two-skeleton of $M$ can be extended as a parallelization of $M$.
\eop

We detail the involved arguments below without mentioning Stiefel-Whitney classes, (actually by almost defining $w_2(TM)$). The elementary proof below can be thought of as an introduction to the obstruction theory used above.

\noindent{\sc Elementary proof of Theorem~\ref{thmpar}:}
Let $M$ be an oriented $3$-manifold. Choose a triangulation of $M$. For any cell $c$ of the triangulation, define an arbitrary trivialization $\tau_c \colon c \times \RR^3 \rightarrow TM_{|c}$ such that $\tau_c$ induces the orientation of $M$.
This defines a trivialization $\tau^{(0)}\colon M^{(0)} \times \RR^3 \rightarrow TM_{|M^{(0)}}$ of $M$ over the $0$-skeleton $M^{(0)}$ of $M$. Let $C_k(M)$ be the set of $k$--cells of the triangulation. Every cell is equipped with an arbitrary orientation.
For an edge $e \in C_1(M)$ of the triangulation, on $\partial e$, $\tau^{(0)}$ reads $\tau^{(0)}= \tau_e \circ \psi_{\RR}(g_e)$ for a map $g_e \colon \partial e \rightarrow GL^+(\RR^3)$. Since $GL^+(\RR^3)$ is connected, $g_e$ extends to $e$, and $\tau^{(1)}= \tau_e  \circ \psi_{\RR}(g_e)$ extends $\tau^{(0)}$ to $e$. Doing so for all the edges extends $\tau^{(0)}$ to a trivialization $\tau^{(1)}$ of the one-skeleton $M^{(1)}$ of $M$.

For an oriented triangle $t$ of the triangulation, on $\partial t$, $\tau^{(1)}$ reads $\tau^{(1)}=  \tau_t \circ \psi_{\RR}(g_t)$ for a map $g_t \colon \partial t \rightarrow GL^+(\RR^3)$. Let $E(t,\tau^{(1)})$ be the homotopy class of $g_t$ in $(\pi_1(GL^+(\RR^3))=\pi_1(SO(3))=\ZZ/2\ZZ)$, $E(t,\tau^{(1)})$ is independent of $\tau_t$.
Then $E(.,\tau^{(1)}) \colon C_2(M) \rightarrow \ZZ/2\ZZ$ is a cochain. When $E(.,\tau^{(1)})=0$, $\tau^{(1)}$ may be extended to
a trivialization $\tau^{(2)}$ over the two-skeleton of $M$, as before.

Since $\pi_2(GL^+(\RR^3))=0$, $\tau^{(2)}$ can next be extended over the three-skeleton of $M$, that is over $M$.

Let us now study the obstruction cochain $E(.,\tau^{(1)})$ whose vanishing guarantees the existence of a parallelization of $M$.

If the map $g_e$ associated to $e$ is changed to $d(e)g_e$ for some $d(e) \colon (e,\partial e) \rightarrow (GL^+(\RR^3),1)$ for every edge $e$, define the associated trivialization $\tau^{(1)\prime}$, and the cochain 
$D(\tau^{(1)},\tau^{(1)\prime}) \colon C_1(M) \rightarrow \ZZ/2\ZZ$ that maps $e$ to the homotopy class of $d(e)$.
Then $(E(.,\tau^{(1)\prime})-E(.,\tau^{(1)}))$ is the coboundary of $D(\tau^{(1)},\tau^{(1)\prime})$.

Let us show that $E(.,\tau^{(1)})$ is a cocycle.
Consider a $3$-simplex $T$, then $\tau^{(0)}$ extends to $T$. Without loss of generality, assume that $\tau_T$
coincides with this extension, that for any face $t$ of $T$, $\tau_t$ is the restriction of $\tau_T$ to $t$, and that the above $\tau^{(1)\prime}$ coincides with $\tau_T$ on the edges of $\partial T$. Then $E(.,\tau^{(1)\prime})(\partial T)=0$. Since a coboundary
also maps $\partial T$ to $0$, $E(.,\tau^{(1)})(\partial T)=0$.

Now, it suffices to prove that the cohomology class of $E(.,\tau^{(1)})$ (which is actually $w_2(TM)$) vanishes in order to prove that there is an extension $\tau^{(1)\prime}$ of $\tau^{(0)}$ on $M^{(1)}$ that extends on $M$.

Since $H^2(M;\ZZ/2\ZZ)=\mbox{Hom}(H_2(M;\ZZ/2\ZZ));\ZZ/2\ZZ)$, it suffices to prove that $E(.,\tau^{(1)})$ maps any $2$--dimensional $\ZZ/2\ZZ$-cycle $C$ to $0$. 

We represent the class of such a cycle $C$ by a non-necessarily orientable closed surface $S$ as follows.
Let ${\neigh}(M^{(0)})$ and ${\neigh}(M^{(1)})$ be small regular neighborhoods of $M^{(0)}$ and $M^{(1)}$ in $M$, respectively, such that 
${\neigh}(M^{(1)}) \cap (M \setminus {\neigh}(M^{(0)}))$ is a disjoint union, running over the edges $e$, of solid cylinders $B_e$ identified with $]0,1[ \times D^2$. The core $]0,1[ \times \{0\}$ of $B_e= ]0,1[ \times D^2$ is a connected part of the interior of the edge $e$. (${\neigh}(M^{(1)})$ is thinner than ${\neigh}(M^{(0)})$.)

Construct $S$ in the complement of ${\neigh}(M^{(0)}) \cup {\neigh}(M^{(1)})$ as the intersection of the support of $C$ with this complement.
Then the closure of $S$ meets the part $[0,1] \times S^1$ of every $\overline{B_e}$ as an even number of parallel intervals from $\{0\} \times S^1$ to $\{1\} \times S^1$. Complete $S$ in $M  \setminus {\neigh}(M^{(0)})$ by connecting the intervals pairwise in $\overline{B_e}$ by disjoint bands. After this operation, the boundary of the closure of $S$ is a disjoint union of circles in the boundary of ${\neigh}(M^{(0)})$, where ${\neigh}(M^{(0)})$ is a disjoint union of balls around the vertices. Glue disjoint disks of ${\neigh}(M^{(0)})$ along these circles to finish the construction of $S$.

Extend $\tau^{(0)}$ to ${\neigh}(M^{(0)})$, assume that $\tau^{(1)}$ coincides with this extension over $M^{(1)}\cap {\neigh}(M^{(0)})$, and extend $\tau^{(1)}$ to ${\neigh}(M^{(1)})$.
Then $TM_{|S}$ is trivial, and we may choose a trivialization $\tau_S$ of $TM$ over $S$ that coincides with our extension of $\tau^{(0)}$ over ${\neigh}(M^{(0)})$, over $S \cap {\neigh}(M^{(0)})$. We have a cell decomposition of $(S,S \cap {\neigh}(M^{(0)}))$ with only 
$1$-cells and $2$-cells, where the $2$-cells of $S$ are in one-to-one canonical correspondence with the $2$-cells of $C$, and one-cells bijectively correspond to bands connecting two-cells in the cylinders $B_e$. These one-cells are equipped with the trivialization of $TM$ induced by $\tau^{(1)}$.
Then we can define $2$--dimensional cochains $E_S(.,\tau^{(1)})$ and $E_S(.,\tau_S)$ from $C_2(S)$ to $\ZZ/2\ZZ$ as before, with respect to this cellular decomposition of $S$, where $(E_S(.,\tau^{(1)})-E_S(.,\tau_S))$
is again a coboundary and $E_S(.,\tau_S)=0$ so that $E_S(C,\tau^{(1)})=0$, and since $E(C,\tau^{(1)})=E_S(C,\tau^{(1)})$, $E(C,\tau^{(1)})=0$ and we are done.
\eop

\subsection{The homomorphism induced by the degree on \texorpdfstring{$[(M,\partial M),(SO(3),1)]$}{[(M, boundary of M),(SO(3),1)]}}

Let $S$ be a non-necessarily orientable closed surface embedded in the interior of $M$, and let $\tau$ be a parallelization of $M$.
We define a twist $g(S,\tau) \in C^0\left((M,\partial M),(SO(3),1)\right)$ below.

The surface $S$ has a tubular neighborhood $\neightub(S)$, which is a $[-1,1]$--bundle over $S$ that admits (orientation-preserving) bundle charts with domains $[-1,1]\times D$ for disks $D$ of $S$ so that the changes of coordinates restrict to the fibers as $\pm \mbox{Identity}$.
Then $$g(S,\tau)\colon (M, \partial M) \longrightarrow (GL^+(\RR^3),1)$$ is the continuous map that maps $M\setminus {\neightub}(S)$ to $1$ such that
$g(S,\tau)((t,s) \in [-1,1]\times D)$ is the rotation with angle $\pi(t+1)$ and with axis
$p_2(\tau^{-1}({\normnu}_s)=(s,p_2(\tau^{-1}({\normnu}_s))))$ where
${\normnu}_s=T_{(0,s)}([-1,1] \times s)$ is the tangent vector to the fiber $[-1,1] \times s$ at $(0,s)$.
Since this rotation coincides with the rotation with opposite axis and with opposite angle $\pi(1-t)$, our map $g(S,\tau)$ is a well-defined continuous map.

Clearly, the homotopy class of $g(S,\tau)$ only depends on the homotopy class of $\tau$ and on the isotopy class of $S$. When $M=B^3$, when $\tau$ is the standard parallelization of $\RR^3$, and when $\frac12 S^2$ denotes the sphere $\frac12 \partial B^3$ inside $B^3$, the homotopy class of $g(\frac12 S^2,\tau)$ coincides with the homotopy class of $\rho$.

\begin{lemma}
\label{lemgStau}
 $H^1(g(S,\tau);\ZZ/2\ZZ)$ is the mod $2$ intersection with $S$.\\
$H^1(.;\ZZ/2\ZZ)\colon [(M,\partial M),(SO(3),1)] \rightarrow H^1(M,\partial M;\ZZ/2\ZZ)$ is onto.
\end{lemma}
\bp The first assertion is obvious, and the second one follows since $H^1(M,\partial M;\ZZ/2\ZZ)$ is the Poincar\'e dual of $H_2(M;\ZZ/2\ZZ)$ and since any element of $H_2(M;\ZZ/2\ZZ)$ is the class of a closed surface.
\eop

\begin{lemma}
\label{lemdeg}
The degree is a group homomorphism $$\mbox{deg} \colon [(M,\partial M),(SO(3),1)] \longrightarrow \ZZ$$
and $\mbox{deg}(\rho_M(B^3)^k)=2k$.
\end{lemma}
\bp
It is easy to see that $\mbox{deg}(fg)=\mbox{deg}(f)+\mbox{deg}(g)$ when $f$ or $g$
is a power of $[\rho_M(B^3)]$.
Let us prove that $\mbox{deg}(f^2)=2\mbox{deg}(f)$ for any $f$.
According to Lemma~\ref{lemgStau}, there is an unoriented embedded surface $S_f$ of the interior of $C$ such that $H^1(f;\ZZ/2\ZZ)=H^1(g(S_f,\tau);\ZZ/2\ZZ)$ for some trivialization $\tau$ of $TM$. Then, according to Lemma~\ref{lempreptrivun}, $fg(S_f,\tau)^{-1}$ is homotopic to some power of $\rho_M(B^3)$, and we are left with the proof that the degree of $g^2$ is $2\mbox{deg}(g)$
for $g=g(S_f,\tau)$. This can easily be done by noticing that $g^2$ is homotopic to $g(S^{(2)}_f,\tau)$ where $S^{(2)}_f$ is the boundary of the tubular neighborhood of $S_f$. In general,
$\mbox{deg}(fg)=\frac{1}2\mbox{deg}((fg)^2)=\frac{1}2\mbox{deg}(f^2g^2)=\frac{1}2\left( \mbox{deg}(f^2)+\mbox{deg}(g^2)\right)$, and the lemma is proved.
\eop

Lemmas~\ref{lempreptrivun} and \ref{lemdeg} imply the following lemma.

\begin{lemma}
\label{lemdegtwo}
The degree induces an isomorphism 

$$\mbox{deg} \colon [(M,\partial M),(SO(3),1)] \otimes_{\ZZ}\QQ \longrightarrow \QQ.$$
Any group homomorphism $\phi \colon [(M,\partial M),(SO(3),1)] \longrightarrow \QQ$ reads
$\frac12\phi(\rho_M(B^3))\mbox{deg}$.
\end{lemma}
\eop

\subsection{First homotopy groups of the groups \texorpdfstring{$SU(n)$}{SU(n)}}

Let $\KK= \RR$ or $\CC$. Let $n \in \NN$. The stabilization maps induced
by the inclusions
$$\begin{array}{llll}i: & GL(\KK^n) & \longrightarrow & GL(\KK \oplus \KK^n)\\
& g & \mapsto & \left(i(g): (x,y) \mapsto (x,g(y))\right)\end{array}$$
will be denoted by $i$. Elements of $GL(\KK^n)$ are represented by matrices whose columns contain the coordinates of the images of the basis elements, with respect to the standard basis of $\KK^n$

See $S^3$ as the unit sphere of $\CC^2$ so that its elements are the pairs $(z_1,z_2)$ of complex numbers such that $|z_1|^2 + |z_2|^2=1$.

The group $SU(2)$ is identified with $S^3$ by the homeomorphism
$$\begin{array}{llll} m^{\CC}_r \colon & S^3 &\rightarrow & SU(2)\\ 
 &(z_1,z_2) &\mapsto & \left[\begin{array}{cc} z_1&-\overline{z}_2\\z_2& \overline{z}_1\end{array} \right] \end{array}$$
so that the first non trivial homotopy group of $SU(2)$ is
$$\pi_3(SU(2))=\ZZ[m^{\CC}_r].$$

The long exact sequence associated to the fibration $$ SU(n-1) \hookfl{i} SU(n) \rightarrow S^{2n-1}$$
shows that $i^n_{\ast}: \pi_j(SU(2)) \longrightarrow \pi_j(SU(n+2))$ is an isomorphism for $j \leq 3$ and $n \geq 0$, 
and in particular, that $\pi_j(SU(4))=\{1\}$ for $j\leq 2$ and
$$\pi_3(SU(4))=\ZZ[i^2(m^{\CC}_r)]$$
where $i^2(m^{\CC}_r)$ 
is the following map
$$\begin{array}{llll} i^2(m^{\CC}_r): &(S^3 \subset \CC^2) & \longrightarrow & SU(4)\\
& (z_1,z_2) & \mapsto & \left[\begin{array}{cccc} 1&0&0&0\\0&1&0&0\\0&0&z_1&-\overline{z}_2\\0&0&z_2& \overline{z}_1\end{array} \right] \end{array}.$$

\subsection{Definition of relative Pontrjagin numbers}
\label{subdefpont}
Let $M_0$ and $M_1$ be two compact connected oriented $3$-manifolds whose boundaries  have collars that are identified by a diffeomorphism. Let $\tau_0\colon M_0 \times \RR^3 \rightarrow TM_0$ and $\tau_1\colon M_1 \times \RR^3 \rightarrow TM_1$ be two parallelizations (which respect the orientations) that agree on the collar neighborhoods of $\partial M_0=\partial M_1$. 
Then the {\em relative Pontrjagin number $p_1(\tau_0,\tau_1)$\/} is the Pontrjagin obstruction to extending the trivialization of $TW \otimes \CC$ induced by $\tau_0$ and $\tau_1$ across the interior of a signature $0$ cobordism $W$ from $M_0$ to $M_1$. Details follow.

Let $M$ be a compact connected oriented $3$-manifold. A {\em special complex trivialization\/} of $TM$ is a trivialization of $TM \otimes \CC$ that is obtained from a trivialization $\tau_M \colon M \times \RR^3 \rightarrow TM$ that induces the orientation of $M$ by composing $(\tau^{\CC}_M = \tau_M \otimes_{\RR} \CC) \colon M \times \CC^3 \rightarrow TM \otimes \CC$ by $$\begin{array}{llll} 
\psi(G): &M \times \CC^3 &\longrightarrow  &M \times \CC^3\\
&(x,y) & \mapsto &(x,G(x)(y))\end{array}$$
for a map $G\colon M \rightarrow SL(3,\CC)$. The definition and properties of relative Pontrjagin numbers, which are given with more details below, are valid for pairs of special complex trivializations.

The {\em signature\/} of a $4$-manifold is the signature of the intersection form on its $H_2(.;\RR)$ (number of positive entries minus number of negative entries in a diagonalised version of this form). Also recall that any closed oriented three-manifold bounds a compact oriented $4$-dimensional manifold whose signature may be arbitrarily changed by connected
sums with copies of $\CC P^2$ or $-\CC P^2$.
A {\em cobordism from $M_0$ to $M_1$\/} is a compact oriented $4$-dimensional manifold $W$ with corners such that
$$\partial W=-M_0\cup_{\partial M_0 \sim 0 \times \partial M_0}  (-[0,1] \times \partial M_0)  \cup_{\partial M_1 \sim 1 \times \partial M_0} M_1,$$ 
is identified with an open subspace of one of the products $[0,1[ \times M_0$ or $]0,1] \times M_1$ near $\partial W$, as the following picture suggests.

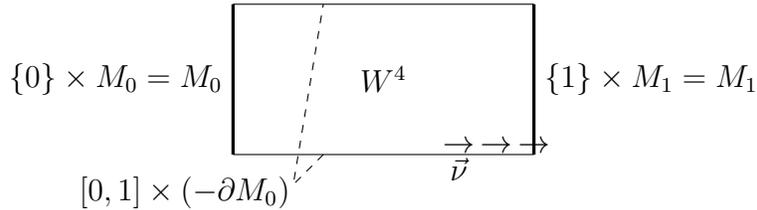
\begin{center}
\begin{tikzpicture}
\useasboundingbox (-3,-.7) rectangle (5.5,2.5);
\draw (0,0) -- (4,0) (4,2) -- (0,2) (2,1) node{$W^4$} (0,1) node[left]{$\{0\} \times M_0=M_0 $} (4,1) node[right]{$\{1\} \times M_1=M_1$};
\draw [very thick] (0,0) -- (0,2) (4,0) -- (4,2);
\draw [dashed] (1.2,0) -- (.8,-.4) -- (1.2,2);
\draw (.9,-.5) node[left]{$[0,1] \times (-\partial M_0)$};
\draw (4,.1) node{$ \rightarrow $} (3.5,.1) node{$ \rightarrow $} (3,.1) node{$ \rightarrow $} node[below](3,-.1) {$ \vec{\normnu}$};
\end{tikzpicture}
\end{center}

Let $W=W^4$ be such a cobordism from $M_0$ to $M_1$, with signature $0$.
Consider the complex $4$-bundle $TW \otimes \CC$ over $W$.
Let $\vec{\normnu}$ be the tangent vector to $[0,1] \times \{\mbox{pt}\}$ over $\partial W$ (under the identifications above), and let 
$\tau(\tau_0,\tau_1)$ denote the trivialization of $TW \otimes \CC$ over $\partial W$ that is obtained by stabilizing either $\tau_0$ or $\tau_1$ into $\vec{\normnu} \oplus \tau_0$ or  $\vec{\normnu} \oplus \tau_1$. Then the obstruction to extending this trivialization to $W$ is the relative first {\em Pontrjagin class\/} $$p_1(W;\tau(\tau_0,\tau_1))[W,\partial W] \in H^4(W,\partial W;\ZZ=\pi_3(SU(4)))=\ZZ[W,\partial W]$$ of the trivialization. 

Now, we specify our sign conventions for this Pontrjagin class. They are the same as in \cite{milnorsta}. 
In particular, $p_1$ is the opposite of the second Chern class $c_2$ of the complexified tangent bundle. See \cite[p. 174]{milnorsta}.
More precisely, equip $M_0$ and $M_1$ with Riemannian metrics that coincide near $\partial M_0$, and equip $W$ with a Riemannian metric that coincides with the orthogonal product metric of one of the products $[0,1] \times M_0$ or $[0,1] \times M_1$ near $\partial W$. 
Equip $TW \otimes \CC$ with the associated hermitian structure. The determinant bundle of $TW$
is trivial because $W$ is oriented, and $\mbox{det}(TW \otimes \CC)$ is also trivial.
Our parallelization $\tau(\tau_0,\tau_1)$ over $\partial W$ is special with respect to
the trivialization of $\mbox{det}(TW \otimes \CC)$.
Up to homotopy, assume that  $\tau(\tau_0,\tau_1)$
is unitary with respect to the hermitian structure of $TW \otimes \CC$ and the standard hermitian form of $\CC^4$.
Since $\pi_i(SU(4))=\{0\}$ when $i<3$, the trivialization $\tau(\tau_0,\tau_1)$
extends to a special unitary trivialization $\tau$ outside the interior of a $4$-ball $B^4$
and defines 
$$\tau\colon S^3 \times \CC^4 \longrightarrow (TW \otimes \CC)_{|S^3}$$
over the boundary $S^3=\partial B^4$ of this $4$-ball $B^4$.
Over this $4$-ball $B^4$, the bundle $TW \otimes \CC$ admits a trivialization
$$\tau_B\colon B^4 \times \CC^4 \longrightarrow (TW \otimes \CC)_{|B^4}.$$
Then $\tau_B^{-1} \circ \tau(v \in S^3, w \in \CC^4)=(v, \phi(v)(w))$, for a map $\phi: S^3 \longrightarrow  SU(4)$ whose homotopy class reads
$$[\phi]=-p_1(W;\tau(\tau_0,\tau_1))[i^2(m^{\CC}_r)] \in \pi_3(SU(4)).$$

Define $p_1(\tau_0,\tau_1)= p_1(W;\tau(\tau_0,\tau_1))$.

\begin{proposition}
\label{proppontdef}
The first {\em Pontrjagin number\/} $p_1(\tau_0,\tau_1)$ is well-defined by the above 
conditions.
\end{proposition}
\bp 
According to the Nokivov additivity theorem, if a closed (compact, without boundary) $4$-manifold $Y$ reads $Y= Y^+ \cup_{X} Y^-$
where $Y^+$ and $Y^-$ are two $4$-manifolds with boundary, embedded in $Y$ that intersect along a closed $3$--manifold $X$ (their common boundary, up to orientation) then $$\mbox{signature}(Y)=\mbox{signature}(Y^+)+\mbox{signature}(Y^-).$$ According to a Rohlin theorem (see \cite{Roh} or \cite[p. 18]{gm2}),
when $Y$ is a compact oriented $4$--manifold without boundary, $p_1(Y)=3\;\mbox{signature}(Y)$.

We only need to prove that $p_1(\tau_0,\tau_1)$ is independent of the signature $0$ cobordism $W$. Let $W_E$
be a $4$-manifold of signature $0$ bounded by $(-\partial W)$. Then $W \cup_{\partial W} W_E$ is a $4$-dimensional manifold without boundary whose signature is $\left(\mbox{signature}(W_E) + \mbox{signature}(W)=0\right)$ by the Novikov additivity theorem.
According to the Rohlin theorem, the first Pontrjagin class of $W \cup_{\partial W} W_E$ is also zero. On the other hand,
this first Pontrjagin class is the sum of the relative first Pontrjagin classes of $W$ and $W_E$ with respect to $\tau(\tau_0,\tau_1)$. These two relative Pontrjagin classes are opposite and therefore the relative first Pontrjagin class of $W$ with respect to $\tau(\tau_0,\tau_1)$ does not depend on $W$.
\eop

Similarly, it is easy to prove the following proposition.
\begin{proposition}
\label{proppontdefsignon}
Under the above assumptions except for the assumption on the signature of the cobordism $W$,
$$p_1(\tau_0,\tau_1)=p_1(W;\tau(\tau_0,\tau_1))-3\;\mbox{\rm signature}(W).$$
\end{proposition}

\subsection{On the groups \texorpdfstring{$SO(3)$}{SO(3)} and \texorpdfstring{$SO(4)$}{SO(4)}}

In this subsection, we describe $\pi_3(SO(4))$ and the natural maps from $\pi_3(SO(3))$ to $\pi_3(SO(4))$ and to $\pi_3(SU(4))$.

The {\em quaternion field $\HH$\/} is the vector space $\CC \oplus \CC j$ equipped with the multiplication that maps $(z_1+z_2j,z^{\prime}_1+z^{\prime}_2j)$ to $(z_1z^{\prime}_1-z_2\overline{z^{\prime}_2}) +(z_2\overline{z^{\prime}_1} +z_1z^{\prime}_2)j$, and with the conjugation that maps $(z_1+z_2j)$ to $\overline{z_1+z_2j}=\overline{z_1}-z_2j$.
The norm of $(z_1+z_2j)$ is the square root of $|z_1|^2+|z_2|^2=(z_1+z_2j)\overline{z_1+z_2j}$, it is multiplicative. Setting $k=ij$, $(1,i,j,k)$ is an orthogonal basis of $\HH$ with respect to the scalar product associated to the norm.
The unit sphere of $\HH$ is the sphere $S^3$, which is equipped with the corresponding group structure.
There are two group morphisms from $S^3$ to $SO(4)$ induced by the multiplication in $\HH$.

$$\begin{array}{llll} m_{\ell} \colon & S^3 &\rightarrow & (SO(\HH)=SO(4))\\ 
& x &\mapsto &m_{\ell}(x) \colon v \mapsto x.v\end{array}$$

$$\begin{array}{llll} \overline{m_r} \colon & S^3 &\rightarrow & SO(\HH)\\ 
& y &\mapsto & (\overline{m_r}(y)\colon v \mapsto v.\overline{y}).\end{array}$$

Together, they induce the group morphism
$$\begin{array}{lll}
S^3 \times S^3 &\rightarrow & SO(4)\\ 
(x,y) &\mapsto & (v \mapsto x.v.\overline{y}).\end{array}$$
The kernel of this group morphism is $\ZZ/2\ZZ(-1,-1)$ so that this morphism is a two-fold covering.
In particular, $$\pi_3(SO(4))=\ZZ[m_{\ell}] \oplus \ZZ[\overline{m_r}].$$

For $\KK= \RR$ or $\CC$ and $n \in \NN$, the $\KK$ (euclidean or hermitian) oriented vector space with the direct orthonormal basis $(v_1, \dots, v_n)$ is denoted by $\KK<v_1, \dots, v_n>$.
There is also the following group morphism $$\begin{array}{llll} \tilde{\rho} \colon & S^3 &\rightarrow & SO(\RR< i , j , k>)=SO(3)\\ 
& x &\mapsto &\left( v \mapsto (v \mapsto x.v.\overline{x}) \right)\end{array}$$
whose kernel is $\ZZ/2\ZZ(-1)$. This morphism $\tilde{\rho}$ is also a two-fold covering.

\begin{lemma}
 This definition of $\tilde{\rho}$ coincides with the previous one, up to homotopy.
\end{lemma}
\bp
It is clear that the two maps coincide up to homotopy, up to orientation since both classes generate $\pi_3(SO(3))=\ZZ$. We take care of the orientation using the outward normal first convention to orient boundaries, as usual.
An element of $S^3$ reads $\mbox{cos}(\theta) + \mbox{sin}(\theta) v$ for a unique $\theta \in[0,\pi]$ and a unit quaternion $v$ with real part zero, which is unique when $\theta \notin \{0,\pi\}$. 
In particular, this defines a diffeomorphism $\phi$ from $]0,\pi[\times S^2$ to $S^3\setminus \{-1,1\}$. We compute the degree of $\phi$ at $\phi(\pi/2,i)$.
The space $\HH$ is oriented as $\RR \oplus \RR< i , j , k>$, where $\RR< i , j , k>$ is oriented by the outward normal to $S^2$, which coincides with the outward normal to $S^3$ in $\RR^4$, followed by the orientation of $S^2$. In particular since $\mbox{cos}$ is an orientation-reversing diffeomorphism at $\pi/2$, the degree of $\phi$ is $1$ and $\phi$ preserves the orientation.
Now $(\mbox{cos}(\theta) + \mbox{sin}(\theta) v)w\overline{(\mbox{cos}(\theta) + \mbox{sin}(\theta) v)}=R(\theta,v)(w)$
where $R(\theta,v)$ is a rotation with axis $v$ for any $v$. Since $R(\theta,i)(j)= \mbox{cos}(2\theta)j + \mbox{sin}(2\theta) k$, the two maps $\tilde{\rho}$ are homotopic. One can check that they are actually conjugate.
\eop

Define
$$\begin{array}{llll} {m_r} \colon & S^3 &\rightarrow & (SO(\HH)=SO(4))\\ 
& y &\mapsto & ({m_r}(y)\colon v \mapsto v.{y}).\end{array}$$

\begin{lemma}
\label{lempitroissoquatrestr} In
$$\pi_3(SO(4))=\ZZ[m_{\ell}] \oplus \ZZ[\overline{m_r}],$$
$$i_{\ast}([\tilde{\rho}])=[m_{\ell}] + [\overline{m}_r]=
[m_{\ell}] - [{m}_r].$$
\end{lemma}
\bp
The $\pi_3$-product in $\pi_3(SO(4))$ coincides with the product induced by the group structure of $SO(4)$.
\eop

\begin{lemma} 
\label{lem237}
Recall that $m_r$
denotes the map from the unit sphere $S^3$ of $\HH$ to $SO(\HH)$ induced by the
right-multiplication.
Denote the inclusions $SO(n) \subset SU(n)$ by $c$. Then in $\pi_3(SU(4))$,
$$c_{\ast}([m_r])=2[i^2(m^{\CC}_r)].$$
\end{lemma}
\bp
Let $\HH + I \HH$ denote the complexification of $\RR^4= \HH= \RR<1,i,j,k>$.
Here, $\CC=\RR \oplus I\RR$.
When $x \in \HH$ and $v \in S^3$, $c(m_r)(v)(Ix)=Ix.v$, and $I^2=-1$.
Let $\varepsilon=\pm 1$, define
$$ \CC^2(\varepsilon)=\CC<\frac{\sqrt{2}}{2}(1+ \varepsilon Ii),\frac{\sqrt{2}}{2}(j+ \varepsilon Ik)>.$$
Consider the quotient $\CC^4/\CC^2(\varepsilon)$.
In this quotient, $Ii=-\varepsilon 1$, $Ik =-\varepsilon j$, and since $I^2=-1$,
$I1=\varepsilon i$ and $Ij=\varepsilon k$. Therefore this quotient is isomorphic to $\HH$ as a real vector space with its complex structure $I= \varepsilon i$.
Then it is easy to see that $c(m_r)$ maps $\CC^2(\varepsilon)$ to $0$ in this quotient. Thus $c(m_r)(\CC^2(\varepsilon)) = \CC^2(\varepsilon)$.
Now, observe that $\HH + I \HH$ is the orthogonal sum of $\CC^2(-1)$ and $\CC^2(1)$.
In particular, $\CC^2(\varepsilon)$ is isomorphic to the quotient $\CC^4/\CC^2(-\varepsilon)$, which is isomorphic to $(\HH;I= -\varepsilon i)$ and 
$c(m_r)$ acts on it by the right multiplication. Therefore,
with respect to the orthonormal basis $\frac{\sqrt{2}}{2}(1-Ii, j-Ik, 1+Ii, j+Ik )$, $c(m_r)$ reads
$$c(m_r)(z_1+z_2j) =\left[\begin{array}{cccc} 
z_1 & -\overline{z}_2 & 0 & 0\\ 
z_2 & \overline{z}_1 & 0 & 0\\
0 & 0 & \overline{z}_1=x_1-Iy_1 & -z_2\\
0 & 0 & \overline{z}_2 & z_1=x_1+Iy_1\\
\end{array} \right].$$
Therefore, the homotopy class of $c(m_r)$ 
is the sum of the homotopy classes of
$$(z_1+z_2j) \mapsto \left[\begin{array}{cc} 
m^{\CC}_r(z_1,z_2) & 0 \\ 
0 & 1
\end{array} \right] \;\; \mbox{and}\;\; (z_1+z_2j) \mapsto \left[\begin{array}{cc} 
1 & 0 \\ 
0 & m^{\CC}_r \circ \iota(z_1,z_2)
\end{array} \right] $$
where $\iota (z_1,z_2)= (\overline{z_1},\overline{z_2})$.
Since the first map is conjugate by a fixed element of $SU(4)$ to $i^2_{\ast}(m^{\CC}_r)$, it is homotopic to $i^2_{\ast}(m^{\CC}_r)$, and since $\iota$ induces the 
identity on $\pi_3(S^3)$, the second map is homotopic to $i^2_{\ast}(m^{\CC}_r)$, too.
\eop

The following lemma finishes to determine the maps $$c_{\ast}: \pi_3(SO(4)) \longrightarrow \pi_3(SU(4))$$
and $$c_{\ast}i_{\ast}: \pi_3(SO(3)) \longrightarrow \pi_3(SU(4)).$$
\begin{lemma}
\label{lempitroissoquatre}
$$c_{\ast}([\overline{m_r}])=c_{\ast}([m_{\ell}])=-2[i^2(m^{\CC}_r)].$$
$$c_{\ast}(i_{\ast}([\tilde{\rho}]))=-4[i^2(m^{\CC}_r)].$$
\end{lemma}
\bp
According to Lemma~\ref{lempitroissoquatrestr}, $i_{\ast}([\tilde{\rho}])=[m_{\ell}] + [\overline{m}_r]=
[m_{\ell}] - [{m}_r].$
Using the conjugacy of quaternions, $m_{\ell}(v)(x)=v.x=\overline{\overline{x}.\overline{v}}=\overline{\overline{m}_{r}(v)(\overline{x})}$.
Therefore $m_{\ell}$ is conjugated to $\overline{m}_{r}$ via 
the conjugacy of quaternions, which lies in $(O(4)\subset U(4))$.

Since $U(4)$ is connected, the conjugacy by an element of $U(4)$ induces the identity on $\pi_3(SU(4))$. 
Thus, $$c_{\ast}([m_{\ell}])=c_{\ast}([\overline{m}_{r}])=-c_{\ast}([{m}_{r}]),$$
and 
$$c_{\ast}(i_{\ast}([\tilde{\rho}]))=-2 c_{\ast}([{m}_{r}]).$$
\eop

\subsection{Relating the relative Pontrjagin number to the degree}

We finish proving Theorem~\ref{thmpone} by proving the following proposition.

\begin{proposition}
\label{proppontdeg}
Let $M_0$ and $M$ be two compact connected oriented $3$-manifolds whose boundaries  have collars that are identified by a diffeomorphism. Let $\tau_0\colon M_0 \times \CC^3 \rightarrow TM_0 \otimes \CC$ and $\tau\colon M \times \CC^3 \rightarrow TM \otimes \CC$ be two special complex trivializations (which respect the orientations) that coincide on the collar neighborhoods of $\partial M_0=\partial M$. Let $[(M, \partial M),(SU(3),1)]$ denote the group of homotopy classes of maps from $M$ to $SU(3)$ that map $\partial M$ to $1$. For any 
$$g: (M, \partial M) \longrightarrow (SU(3),1),$$
define
$$\begin{array}{llll} 
\psi(g): &M \times \CC^3 &\longrightarrow  &M \times \CC^3\\
&(x,y) & \mapsto &(x,g(x)(y))\end{array}$$
then 
$$p_1(\tau_0,\tau \circ \psi(g))-p_1(\tau_0,\tau)=p_1(\tau,\tau \circ \psi(g))=-p_1(\tau \circ \psi(g),\tau)=p^{\prime}_1(g)$$ is independent of $\tau_0$ and $\tau$,
$p^{\prime}_1$ induces an isomorphism from the group $[(M, \partial M),(SU(3),1)]$ to $\ZZ$,
and, if $g$ is valued in $SO(3)$, then $$p^{\prime}_1(g)=2\mbox{deg}(g).$$
\end{proposition}

In order to prove this proposition, we first prove the following lemma.

\begin{lemma}
\label{lempunind}
Under the hypotheses of Proposition~\ref{proppontdeg}, $\left(p_1(\tau_0,\tau \circ \psi(g))-p_1(\tau_0,\tau)\right)$ is independent of $\tau_0$ and $\tau$.
\end{lemma}
\bp Indeed, $\left(p_1(\tau_0,\tau \circ \psi(g))-p_1(\tau_0,\tau)\right)$ can be defined as the obstruction to extending the following trivialization of the complexified tangent bundle to
$[0,1] \times M$ restricted to the boundary. This trivialization is
$T[0,1] \oplus \tau$ on $(\{0\} \times M) \cup ([0,1] \times \partial M)$ and $T[0,1] \oplus \tau \circ \psi(g)$ on $\{1\} \times M$. But this obstruction is the obstruction
to extending the map $\tilde{g}$ from $\partial([0,1] \times M)$ to $SU(4)$ that maps  $(\{0\} \times M) \cup ([0,1] \times \partial M)$ to $1$ and that coincides with $i(g)$ on $\{1\} \times M$, regarded as a map from $\partial([0,1] \times M)$ to $SU(4)$, over $([0,1] \times M)$. This obstruction, which lies in 
$\pi_3(SU(4))$ since $\pi_i(SU(4))=0$, for $i<3$, is independent of $\tau_0$ and $\tau$.
\eop

\noindent{\sc Proof of Proposition~\ref{proppontdeg}:}
Lemma~\ref{lempunind} guarantees that $p^{\prime}_1$ defines two group homomorphisms to $\ZZ$ from $[(M, \partial M),(SU(3),1)]$ and from $[(M, \partial M),(S0(3),1)]$. 
Since $\pi_i(SU(3))$ is trivial for $i<3$ and since $\pi_3(SU(3))=\ZZ$, the group of homotopy classes  $[(M, \partial M) , (SU(3),1)]$ is generated by the class of a map that maps the complement of a $3$-ball $B$ to $1$ and that factors through a map that generates $\pi_3(SU(3))$. By definition of the Pontrjagin classes, $p^{\prime}_1$ sends such a generator to $\pm 1$ and it induces an isomorphism from $[(M, \partial M),(SU(3),1)]$ to $\ZZ$. 

According to Lemma~\ref{lempreptrivun} and to Lemma~\ref{lemdegtwo}, the restriction of $p^{\prime}_1$ to $[(M, \partial M),(SO(3),1)]$ must read
$p^{\prime}_1(\rho_M(B^3)) \frac{\mbox{deg}}{2}$, and we are left with the proof of the following lemma.

\begin{lemma}
\label{lemvarpun}
$$p^{\prime}_1(\rho_M(B^3))=4.$$
\end{lemma}
Let $g=\rho_M(B^3)$, we can extend $\tilde{g}$ (defined in the proof of Lemma~\ref{lempunind}) by the constant map with value 1 outside
$[\varepsilon, 1] \times B^3 \cong B^4$ and, in $\pi_3(SU(4))$
$$[c(\tilde{g}_{|\partial B^4})]=-p_1(\tau,\tau \circ \psi(g))[i^2(m^{\CC}_r)].$$
Since $\tilde{g}_{|\partial B^4}$ is homotopic to $c \circ i(\tilde{\rho})$, Lemma~\ref{lempitroissoquatre} allows us to conclude.
\eop

\section{Other complements}
\setcounter{equation}{0}

\subsection{More on low-dimensional manifolds}
\label{submorelowdif}

Piecewise linear (or PL) $n$--manifolds can be defined as the $C^i$-manifolds of Subsection~\ref{subbackground} by replacing $C^i$ with piecewise linear (or PL).

When $n \leq 3$, the above notion of PL-manifold coincides with the notions of smooth and topological manifold, according to the following theorem. This is not true anymore when $n>3$. See 
\cite{kui}.

\begin{theorem} 
\label{thmstructhree}
When $n \leq 3$, the category of topological $n$--manifolds is isomorphic to the category of PL $n$--manifolds and to the category of $C^r$ $n$--manifolds, for $r=1,\dots,\infty$.
\end{theorem}

For example, according to this statement, which contains several theorems (see \cite{kui}), any topological $3$--manifold has a unique $C^{\infty}$--structure. Below $n=3$.

The equivalence between the $C^i, i=1,2,\dots,\infty$-categories follows from work of Whitney in 1936 \cite{whi}. In 1934, Cairns \cite{cai1} provided a map from the $C^{1}$--category to the PL category, which shows the existence of a triangulation for $C^{1}$--manifolds, and he proved that this map is onto \cite[Theorem III]{cai2} in 1940. Moise \cite{moise} proved the equivalence between the topological category and the PL category in 1952.
This diagram was completed by Munkres \cite[Theorem 6.3]{munk} and Whitehead \cite{white} in 1960 by their independent proofs of the injectivity of the natural map from the $C^1$--category to the topological category.

\section*{Index of notations}
\addcontentsline{toc}{section}{Index of notations}
\begin{multicols}{2}
\begin{theindex}

  \item $\CA_n(C)$, \hyperpage{28}
  \item $\CA_n(\emptyset)$, \hyperpage{28}
  \item $\CA^t_n(C)$, \hyperpage{27}

  \indexspace

  \item $\beta_{\Gamma}$, \hyperpage{36}

  \indexspace

  \item $\CD^c(\RR)$, \hyperpage{48}
  \item $\CD^c_n(\RR)$, \hyperpage{43}
  \item $\CD^e_n(C)$, \hyperpage{33}

  \indexspace

  \item $I_{\theta}(K_j,\tau)$, \hyperpage{34}

  \indexspace

  \item $\omega_{S^2}$, \hyperpage{4}

  \indexspace

  \item $\psi_{\RR}$, \hyperpage{17}
  \item $p_{\tau}$, \hyperpage{14}

  \indexspace

  \item $\taust$, \hyperpage{13}
  \item $\Theta(M)$, \hyperpage{18}
  \item $\Theta(M,\tau)$, \hyperpage{15}

  \indexspace

  \item $Z$, \hyperpage{35}
  \item $Z^u$, \hyperpage{35}

\end{theindex}
\end{multicols}

\def\cprime{$'$}
\providecommand{\bysame}{\leavevmode ---\ }
\providecommand{\og}{``}
\providecommand{\fg}{''}
\providecommand{\smfandname}{and}
\providecommand{\smfedsname}{eds.}
\providecommand{\smfedname}{ed.}
\providecommand{\smfmastersthesisname}{M\'emoire}
\providecommand{\smfphdthesisname}{Th\`ese}

\end{document}